\documentclass[10pt]{amsart}
\usepackage[numbers]{natbib}
\usepackage{amsmath}
\usepackage{amsfonts}
\usepackage{bbm}
\usepackage{mathrsfs}
\usepackage{scalerel}
\usepackage{dsfont}
\usepackage{stmaryrd}
\usepackage{xifthen}
\usepackage{graphicx}
\usepackage{booktabs}
\usepackage{paralist}
\usepackage[shortlabels,inline]{enumitem}
\usepackage{animate}
\usepackage{lineno}
\usepackage{nameref}
\usepackage[hidelinks]{hyperref}  
\DeclareMathAlphabet{\Mymathbb}{U}{bbold}{m}{n}

\theoremstyle{remark}
\newtheorem{remark}{Remark}

\newcommand{\T}{\mbox{\tiny T}}

\newcommand{\timeSymbol}{t}
\newcommand{\Time}{\timeSymbol} 
\newcommand{\TimeF}{{\timeSymbol}_{\scriptscriptstyle{f}}} 
\newcommand{\tempo}{\Time} 

\newcommand{\FirstSymbol}{\mathcal{G}}
\newcommand{\First}[1][]{
  \ifthenelse{\equal{#1}{}}
  {\FirstSymbol}
  {\FirstSymbol_{#1}}
}
\newcommand{\STFirstSymbol}{\mathbb{G}}
\newcommand{\STFirst}[1][]{
  \ifthenelse{\equal{#1}{}}
  {\STFirstSymbol}
  {\STFirstSymbol_{#1}}
}

\newcommand{\FirstForm}[1][]{
  \ifthenelse{\equal{#1}{}}       
  {\operatorname{I_{\point}}}             
  {\operatorname{I_{#1}}}             
}
\newcommand{\GradSymbol}{\operatorname{\mathbf{\nabla}}}

\newcommand{\GradSurf}{\GradSymbol_{\scriptscriptstyle\First}}
\newcommand{\GradSurfST}{\GradSymbol_{\scriptscriptstyle\STFirst}}

\newcommand{\DivSurf}{\GradSurf\cdot}
\newcommand{\DivSurfST}{\GradSurfST\cdot}

\newcommand{\Der}[1][]
{
  \ifthenelse{\equal{#1}{}}
  {\partial}
  {\partial_{\scriptscriptstyle{#1}}}
}

\newcommand{\DerPar}[2]{\frac{\Der #1}{\Der #2}}

\newcommand{\DerTot}[2][t]
{
  \ifthenelse{\equal{#2}{}}
  {\frac{d #2}{dt}}
  {\frac{d #2}{d #1}}
}
\newcommand{\Diffsymbol}{\operatorname{d}\!}
\newcommand{\Diff}[2][]{
  \ifthenelse{\equal{#1}{}}
  {\Diffsymbol{#2}}
  {\Diffsymbol{#2}_{\scriptscriptstyle{#1}}}
}
\newcommand{\DirDerSymb}{D}
\newcommand{\DirDer}[2][]
{
  \ifthenelse{\equal{#1}{}}
  {\DirDerSymb^{#2}}
  {\DirDerSymb^{#2}_{#1}}
}

\newcommand{\BigO}[1]{\ensuremath{\operatorname{\mathcal{O}}\!\left(#1\right)}}

\newcommand{\REALsymbol}{\mathbb{R}}
\newcommand{\REAL}[1][]{
  \ifthenelse{\equal{#1}{}}
  {\REALsymbol}
  {{\REALsymbol}^{#1}}
}
\newcommand{\EUCLsymbol}{\mathbb E}
\newcommand{\EUCL}[1][]{
  \ifthenelse{\equal{#1}{}}
  {\EUCLsymbol}
  {{\EUCLsymbol}^{#1}}
}
\newcommand{\NATURALsymbol}{\mathbb N}
\newcommand{\NATURAL}[1][]{
  \ifthenelse{\equal{#1}{}}
  {\NATURALsymbol}
  {{\NATURALsymbol}^{#1}}
}

\newcommand{\ABS}[2][]
{
  \ifthenelse{\equal{#1}{}}
  {\left| #2 \right|}
  {\left| #2 \right|_{#1}}
}

\newcommand{\NORM}[2][]
{
  \ifthenelse{\equal{#1}{}}
  {\left\| #2 \right\|}
  {\left\| #2 \right\|_{#1}}
}

\newcommand{\SCAL}[3][]
{
  \ifthenelse{\equal{#1}{}}
  {\left\langle{#2},{#3}\right\rangle}
  {\left\langle{#2},{#3}\right\rangle_{#1}}
}
\newcommand{\SCALF}[3][]
{
  \ifthenelse{\equal{#1}{}}
  {\left({#2},{#3}\right)}
  {\left({#2},{#3}\right)_{#1}}
}

\newcommand{\scalprodSurf}[3][]
{
  \ifthenelse{\equal{#1}{}}
  {\left\langle {#2},{#3} \right\rangle_{\scriptscriptstyle\First}}
  {\left\langle {#2},{#3} \right\rangle_{{#1}}}
}
\newcommand{\scalprodSurfST}[3][]
{
  \ifthenelse{\equal{#1}{}}
  {\left\langle {#2},{#3} \right\rangle_{\scriptscriptstyle\STFirst}}
  {\left\langle {#2},{#3} \right\rangle_{{#1}}}
}

\newcommand{\Point}{\mathbf{P}}
\newcommand{\point}[1][]
{
  \ifthenelse{\equal{#1}{}}
  {\mathbf{p}}
  {\mathbf{p}_{#1}}
}

\newcommand{\midPoint}{\mathbf{m}}

\newcommand{\Dt}{\Delta\tempo}

\newcommand{\Closure}[1]{\overline{#1}}

\newcommand{\RegionSymb}{R}
\newcommand{\Region}[1][]
{
  \ifthenelse{\equal{#1}{}}
  {\RegionSymb}
  {\RegionSymb^{#1}}
}
\newcommand{\bRegion}[1][]
{
  \ifthenelse{\equal{#1}{}}
  {\overline{\RegionSymb}}
  {\overline{\RegionSymb}^{#1}}
}

\newcommand{\diam}{\operatorname{diam}}

\newcommand{\IT}[1][] 
{
  \ifthenelse{\equal{#1}{}}
  {[0,\TimeF]}
  {[\Time^{#1},\Time^{#1+1}]}
}
\newcommand{\Interval}[1][]
{
  \ifthenelse{\equal{#1}{}}
  {I}
  {I_{#1}}
}

\newcommand{\SurfDomainsymb}{\Gamma}
\newcommand{\SurfDomain}[1][]{
  \ifthenelse{\equal{#1}{}}
  {\SurfDomainsymb}
  {\SurfDomainsymb_{#1}}
}
\newcommand{\SurfDomainBndsymb}{\partial\Gamma}
\newcommand{\SurfDomainBnd}[1][]{
  \ifthenelse{\equal{#1}{}}
  {\SurfDomainBndsymb}
  {\SurfDomainBndsymb_{#1}}
}
\newcommand{\ClosedSurfDomain}[1][]{
  \ifthenelse{\equal{#1}{}}
  {\Closedsymb{\SurfDomain}}
  {\Closedsymb{\SurfDomain[#1]}}
}

\newcommand{\heightsymb}{\mathcal{B}}
\newcommand{\height}[1][]{
  \ifthenelse{\equal{#1}{}}
  {\heightsymb}
  {\heightsymb_{#1}}
}
\newcommand{\obstaclesymb}{\mathcal{Z}}
\newcommand{\obstacle}[1][]{
  \ifthenelse{\equal{#1}{}}
  {\obstaclesymb}
  {\obstaclesymb_{#1}}
}

\newcommand{\BSMsymbol}{\heightsymb}

\newcommand{\BSM}[1][]
{
  \ifthenelse{\equal{#1}{}}
  {\BSMsymbol}
  {\BSMsymbol_{#1}}
}

\newcommand{\Surf}{\mathcal{S}}

\newcommand{\SurfBndsymb}{\partial\Surf}
\newcommand{\SurfBnd}[1][]{
  \ifthenelse{\equal{#1}{}}
  {\SurfBndsymb}
  {\SurfBndsymb_{#1}}
}
\newcommand{\Closedsymb}[1]{\bar{#1}}
\newcommand{\ClosedSurf}[1][]{
  \ifthenelse{\equal{#1}{}}
  {\Closedsymb{\Surf}}
  {\Closedsymb{\Surf[#1]}}
}

\newcommand{\Vector}[1]{\mathbf{#1}}

\newcommand{\grav}{g}

\newcommand{\Depth}{\eta}
\newcommand{\density}{\rho}

\newcommand{\MapUsymb}{\phi}
\newcommand{\MapU}[1][]
{
  \ifthenelse{\equal{#1}{}}
    {\MapUsymb}
    {\MapUsymb_{#1}}
}
\newcommand{\MapVsymb}{\psi}
\newcommand{\MapV}[1][]
{
  \ifthenelse{\equal{#1}{}}
    {\MapVsymb}
    {\MapVsymb_{#1}}
}
\newcommand{\Transsymb}{\Phi}
\newcommand{\Trans}[1][]
{
  \ifthenelse{\equal{#1}{}}
    {\Transsymb} 
    {\Transsymb_{\scriptscriptstyle{#1}}}
}
\newcommand{\InvMapsymb}{\Psi}
\newcommand{\InvMap}[1][]
{
  \ifthenelse{\equal{#1}{}}
    {\InvMapsymb}
    {\InvMapsymb_{\scriptscriptstyle{#1}}}
}

\newcommand{\fsymb}{f}
\newcommand{\scalFun}[1][]
{
  \ifthenelse{\equal{#1}{}}
  {\fsymb}
  {\fsymb_{#1}}
}
\newcommand{\tscalFun}[1][]
{
  \ifthenelse{\equal{#1}{}}
  {\tilde{\fsymb}}
  {\tilde{\fsymb}_{#1}}
}
\newcommand{\gsymb}{g}
\newcommand{\scalFung}[1][]
{
  \ifthenelse{\equal{#1}{}}
  {\gsymb}
  {\gsymb_{#1}}
}
\newcommand{\Fsymb}{F}
\newcommand{\FvecFun}[1][]
{
  \ifthenelse{\equal{#1}{}}
  {\Fsymb}
  {\Fsymb_{#1}}
}
\newcommand{\tFvecFun}[1][]
{
  \ifthenelse{\equal{#1}{}}
  {\tilde{\Fsymb}}
  {\tilde{\Fsymb}_{#1}}
}
\newcommand{\Ffunc}[2][]
{
  \ifthenelse{\equal{#1}{}}
  {\ensuremath{\Fsymb_{#2}}}
  {\ensuremath{\Fsymb^{#1}_{#2}}}
}

\newcommand{\PrincipalK}[1][]
{
  \ifthenelse{\equal{#1}{}}
  {k}
  {k_{#1}}
}

\newcommand{\vecsymb}{u}
\newcommand{\vecFun}[1][]{
  \ifthenelse{\equal{#1}{}}
  {\Vector{\vecsymb}}
  {\vecsymb^{#1}}
}
\newcommand{\tvecFun}[1][]{
  \ifthenelse{\equal{#1}{}}
  {\tilde{\vecsymb}}
  {\tilde{\vecsymb}_{#1}}
}

\newcommand{\vvsymb}{v}
\newcommand{\vv}[1][]{
  \ifthenelse{\equal{#1}{}}
  {\mathbf{\vvsymb}}
  {\vvsymb^{#1}}
}
\newcommand{\wwsymb}{w}
\newcommand{\ww}[1][]
{
  \ifthenelse{\equal{#1}{}}
  {\mathbf{\wwsymb}}
  {\wwsymb^{#1}}
}
\newcommand{\uusymb}{u}
\newcommand{\uu}[1][]
{
  \ifthenelse{\equal{#1}{}}
  {\mathbf{\uusymb}}
  {\uusymb^{#1}}
}

\newcommand{\VecFieldSymbol}{X}
\newcommand{\VecField}[1][]
{
  \ifthenelse{\equal{#1}{}}
  {\VecFieldSymbol}
  {\VecFieldSymbol^{#1}}
}
\newcommand{\VecFieldYSymbol}{Y}
\newcommand{\VecFieldY}[1][]
{
  \ifthenelse{\equal{#1}{}}
  {\VecFieldYSymbol}
  {\VecFieldYSymbol^{#1}}
}

\newcommand{\xvsymb}{x}
\newcommand{\xv}[1][]
{
  \ifthenelse{\equal{#1}{}}
  {\mathbf{\xvsymb}}
  {\mathbf{\xvsymb}_{\scriptscriptstyle{#1}}}
}
\newcommand{\xvcomp}[1][]{
  \ifthenelse{\equal{#1}{}}
  {\xvsymb}
  {\xvsymb^{\scriptscriptstyle{#1}}}
}
\newcommand{\xcg}[1][]{
  \ifthenelse{\equal{#1}{}}
  {\xvcomp[1]}
  {\xvcomp[1]_{\scriptscriptstyle{#1}}}
}
\newcommand{\ycg}[1][]{
  \ifthenelse{\equal{#1}{}}
  {\xvcomp[2]}
  {\xvcomp[2]_{\scriptscriptstyle{#1}}}
}
\newcommand{\zcg}[1][]{
  \ifthenelse{\equal{#1}{}}
  {\xvcomp[3]}
  {\xvcomp[3]_{\scriptscriptstyle{#1}}}
}

\newcommand{\svsymb}{s}
\newcommand{\sv}[1][]{
  \ifthenelse{\equal{#1}{}}
  {\mathbf{\svsymb}}
  {\mathbf{\svsymb}_{\scriptscriptstyle{#1}}}
}
\newcommand{\svcomp}[1][]
{
  \ifthenelse{\equal{#1}{}}
  {\svsymb}
  {\svsymb^{\scriptscriptstyle{#1}}}
}
\newcommand{\xcl}[1][]{
  \ifthenelse{\equal{#1}{}}
   {\svcomp[1]}
   {\svcomp[1]_{\scriptscriptstyle{#1}}}
}
\newcommand{\ycl}[1][]{
  \ifthenelse{\equal{#1}{}}
   {\svcomp[2]}
   {\svcomp[2]_{\scriptscriptstyle{#1}}}
}
\newcommand{\zcl}[1][]{
  \ifthenelse{\equal{#1}{}}
   {\svcomp[3]}
   {\svcomp[3]_{\scriptscriptstyle{#1}}}
}

\newcommand{\ProjSymb}{\operatorname{\pi}}
\newcommand{\ProjFun}[2][]{
  \ifthenelse{\equal{#1}{}}
  {\ProjSymb\left(#2\right)}
  {\ProjSymb_{\scriptscriptstyle{#1}}\left(#2\right)}
}
\newcommand{\Prm}[1][]
{
  \ifthenelse{\equal{#1}{}}
  {\operatorname{pr}}
  {\operatorname{pr}_{\scriptscriptstyle{#1}}}
}
\newcommand{\TanPlane}[2][]
{
  \ifthenelse{\equal{#1}{}}
  {T_{\scriptscriptstyle{\point}}#2}
  {T_{\scriptscriptstyle{#1}}#2}
}

\newcommand{\SubsetSymbol}{\mathcal{U}}
\newcommand{\SubsetU}[1][]
{
  \ifthenelse{\equal{#1}{}}
  {{U}}
  {{U}_{#1}}
}
\newcommand{\SubsetV}[1][]
{
  \ifthenelse{\equal{#1}{}}
  {{V}}
  {{V}_{#1}}
}
\newcommand{\SubsetW}[1][]
{
  \ifthenelse{\equal{#1}{}}
  {{W}}
  {{W}_{#1}}
}
\newcommand{\NeighSymbol}{\mathcal{N}}
\newcommand{\Neigh}[1][]
{
  \ifthenelse{\equal{#1}{}}
  {\NeighSymbol_{\point}}
  {\NeighSymbol_{#1}}
}
\newcommand{\NeighSurf}[1][]
{
  \ifthenelse{\equal{#1}{}}
  {\SubsetSymbol_{\point}}
  {\SubsetSymbol_{#1}}
}
\newcommand{\NormSymb}{\mathbf{N}}
\newcommand{\normalvec}{\NormSymb}
\newcommand{\normalSurf}[1][]
{
  \ifthenelse{\equal{#1}{}}
  {\NormSymb}
  {\NormSymb(#1)}
}
\newcommand{\normalInterp}[1][]
{
  \ifthenelse{\equal{#1}{}}
   {\tilde{\NormSymb}}
   {\tilde{\NormSymb}_{\scriptscriptstyle{#1}}}
}

\newcommand{\normalEdge}{\mathbf{\nu}} 
\newcommand{\basisCC}{t}
\newcommand{\basisGC}{e}
\newcommand{\vecBaseGC}[1][]
{
  \ifthenelse{\equal{#1}{}}
  {\mathbf{\basisGC}}
  {\mathbf{\basisGC}_{#1}}
}
\newcommand{\vecBasePhys}[1][]
{
  \ifthenelse{\equal{#1}{}}
  {\mathbf{\basisGC}}
  {\mathbf{\basisGC}_{#1}}
}

\newcommand{\vecBaseCCcv}[1][]
{
  \ifthenelse{\equal{#1}{}}
  {\mathbf{\basisCC}}
  {\mathbf{\basisCC}_{#1}}
}

\newcommand{\tvecBaseCCcv}[1][]
{
  \ifthenelse{\equal{#1}{}}
  {\tilde{\mathbf{\basisCC}}}
  {\tilde{\mathbf{\basisCC}}_{#1}}
}
\newcommand{\hvecBaseCCcv}[1][]
{
  \ifthenelse{\equal{#1}{}}
  {\hat{\mathbf{\basisCC}}}
  {\hat{\mathbf{\basisCC}}_{#1}}
}
\newcommand{\vecBaseCCctrv}[1][]
{
  \ifthenelse{\equal{#1}{}}
  {\mathbf{\basisCC}}
  {\mathbf{\basisCC}^{#1}}
}

\newcommand{\metrcoefsymbol}{h}
\newcommand{\metrcoefH}[1]{\metrcoefsymbol_{\mbox{\tiny{(#1)}}}}

\newcommand{\first}[1]{
  \IfEqCase{#1}{
    {1}{\operatorname{E}}
    {2}{\operatorname{F}}
    {3}{\operatorname{G}}
  }
  [\PackageError{first}{Undefined option to first: #1}{}]%
}
\newcommand{\SecondFormSymbol}{\ensuremath{\operatorname{II}}}
\newcommand{\SecondForm}[1][]
{
  \ifthenelse{\equal{#1}{}}
  {\SecondFormSymbol_{\point}}
  {\SecondFormSymbol_{#1}}
}
\newcommand{\second}[1]{
  \IfEqCase{#1}{
    {1}{\operatorname{e}}
    {2}{\operatorname{f}}
    {3}{\operatorname{g}}
  }
  [\PackageError{first}{Undefined option to first: #1}{}]
}
\newcommand{\WeigSymbol}{\mathcal{W}}
\newcommand{\Weig}[1][]
{
  \ifthenelse{\equal{#1}{}}
  {\WeigSymbol}
  {\WeigSymbol_{#1}}
}

\newcommand{\velSymbol}{u}
\newcommand{\vectvel}[1][]
{
   \ifthenelse{\equal{#1}{}}
   {\vec{\velSymbol}}
   {\vec{\velSymbol}(#1)}
}

\newcommand{\velcompContr}[2][i]
{
   \ifthenelse{\equal{#2}{}}
   {\velSymbol^{#1}}
   {\velSymbol^{#1}(#2)}}
\newcommand{\velcompPhys}[2][i]
{
   \ifthenelse{\equal{#2}{}}
   {\velSymbol_{(#1)}}
   {\velSymbol_{(#1)}(#2)}}

\newcommand{\velSymbolRP}{v}
\newcommand{\velcompContrRP}[2][i]
{
   \ifthenelse{\equal{#2}{}}
   {\velSymbolRP^{#1}}
   {\velSymbolRP^{#1}(#2)}}
\newcommand{\velRP}[1][]
{
  \ifthenelse{\equal{#1}{}}
  {\velSymbolRP}
  {\velSymbolRP_{#1}}
}

\newcommand{\velcompApprox}[2][i]
{
   \ifthenelse{\equal{#2}{}}
   {\velSymbol^{#1)}}
   {\velSymbol^{#1}_{(#2)}}
}

\newcommand{\VelSymbol}{U}
\newcommand{\vectVel}[1][]
{
   \ifthenelse{\equal{#1}{}}
   {\vec{\VelSymbol}}
   {\vec{\VelSymbol}(#1)}
}
\newcommand{\Velcomp}[2][i]
{
   \ifthenelse{\equal{#2}{}}
   {\VelSymbol^{#1}}
   {\VelSymbol^{#1}(#2)}
}
\newcommand{\VprimoSymbol}{\tilde{u}}
\newcommand{\Vprimo}[1][]
{
   \ifthenelse{\equal{#1}{}}
   {\VprimoSymbol}
   {\VprimoSymbol(#1)}
}
\newcommand{\VprimoComp}[2][i]
{
   \ifthenelse{\equal{#2}{}}
   {\VprimoSymbol^{#1}}
   {\VprimoSymbol^{#1}(#2)}
}
\newcommand{\ttvelSymbol}{\tilde{\Mymathbb{u}}}
\newcommand{\ttvel}[1][]
{
   \ifthenelse{\equal{#1}{}}
   {\mathbf{\ttvelSymbol}}
   {\mathbf{\ttvelSymbol}(#1)}
}
\newcommand{\ttvelComp}[2][i]
{
   \ifthenelse{\equal{#2}{}}
   {\ttvelSymbol^{#1}}
   {\ttvelSymbol^{#1}(#2)}
}

\newcommand{\MatAlphaSymbol}{\mathbb{A}}
\newcommand{\MatAlpha}[1][]{%
  \ifthenelse{\equal{#1}{}}
  {\MatAlphaSymbol}
  {\MatAlphaSymbol_{#1}}
}

\newcommand{\QSymbol}{q}
\newcommand{\Qdisch}[1][]
{
   \ifthenelse{\equal{#1}{}}
   {\mathbf{\QSymbol}}
   {\mathbf{\QSymbol}(#1)}
}
\newcommand{\Qcomp}[2][i]
{
   \ifthenelse{\equal{#2}{}}
   {\QSymbol^{#1}}
   {\QSymbol^{#1}(#2)}
}
\newcommand{\Qvect}[1][]
{
   \ifthenelse{\equal{#1}{}}
   {\vec{\QSymbol}}
   {\vec{\QSymbol}(#1)}
}
\newcommand{\FricSymbol}{f}
\newcommand{\vectFric}[1][]
{
   \ifthenelse{\equal{#1}{}}
   {\mathbf{\FricSymbol}}
   {\mathbf{\FricSymbol}_{\scriptscriptstyle{#1}}}
}
\newcommand{\Friccomp}[2][i]
{
   \ifthenelse{\equal{#2}{}}
   {\FricSymbol_{#1}}
   {\FricSymbol_{#1}(#2)}}
\newcommand{\BFsymbol}{\tau}
\newcommand{\BottomFriction}[1][]{
  \ifthenelse{\equal{#1}{}}
  {\BFsymbol_{b}}
  {\BFsymbol_{b}^{#1}}
}

\newcommand{\ProjMatSymb}{\mathbb{P}}
\newcommand{\ProjMat}[1][]{
  \ifthenelse{\equal{#1}{}}
  {\ProjMatSymb}
  {\ProjMatSymb_{#1}}
}
\newcommand{\IDSymbol}{\mathbb{I}}
\newcommand{\IDtens}[1][]{
  \ifthenelse{\equal{#1}{}}
  {\IDSymbol}
  {\IDSymbol(#1)}
}

\newcommand{\tenscomp}[2][ij]
{
  \ifthenelse{\equal{#2}{}}
  {\tenscompSymbol^{#1}}
  {\tenscompSymbol^{#1}(#2)}
}
\newcommand{\tensrow}[2][i]
{
  \ifthenelse{\equal{#2}{}}
  {\tensSymbol^{(#1)}}
  {\tensSymbol^{(#1)}(#2)}
}
\newcommand{\TensSymbol}{\mathbf{T}}
\newcommand{\Tens}[1][]{
  \ifthenelse{\equal{#1}{}}
  {\TensSymbol}
  {\TensSymbol_{#1}}
}

\newcommand{\TensCompSymbol}{\TensSymbol}
\newcommand{\TensComp}[2][ij]
{
  \ifthenelse{\equal{#2}{}}
  {\TensCompSymbol^{#1}}
  {\TensCompSymbol^{#1}(#2)}
}

\newcommand{\tensPrimoSymbol}{\tilde{\mathbf{\tau}}}
\newcommand{\tensPrimo}[1][]{
  \ifthenelse{\equal{#1}{}}
  {\tensPrimoSymbol}
  {\tensPrimoSymbol(#1)}
}
\newcommand{\tensPrimoCompSymbol}{\tensPrimoSymbol}
\newcommand{\tensPrimoComp}[2][ij]
{
  \ifthenelse{\equal{#2}{}}
  {\tensPrimoCompSymbol^{#1}}
  {\tensPrimoCompSymbol^{#1}(#2)}
}
\newcommand{\MCxl}[1][]{\ifthenelse{\equal{#1}{}}{h_{(1)}}{h_{(1),#1}}} 
\newcommand{\MCyl}[1][]{\ifthenelse{\equal{#1}{}}{h_{(2)}}{h_{(2),#1}}} 
\newcommand{\MCzl}[1][]{\ifthenelse{\equal{#1}{}}{h_{(3)}}{h_{(3),#1}}}

\newcommand{\MPsymb}{h}

\newcommand{\smallestH}[1][]
{
  \ifthenelse{\equal{#1}{}}
    {{l}}
    {{l}_{#1}}
}
\newcommand{\meshparam}[1][]
{
  \ifthenelse{\equal{#1}{}}
    {\MPsymb}
    {\MPsymb_{#1}}
}
\newcommand{\InradiusSymbol}{r}
\newcommand{\Inradius}[1][]
{
  \ifthenelse{\equal{#1}{}}
  {\InradiusSymbol}
  {\InradiusSymbol_{\scriptscriptstyle{#1}}}
}
\newcommand{\Tsymb}{\mathcal{T}}
\newcommand{\Triang}[1][]
{
  \ifthenelse{\equal{#1}{}}
    {\Tsymb}
    {\Tsymb_{#1}}
}
\newcommand{\TriangH}[1][]
{
  \ifthenelse{\equal{#1}{}}
    {\Tsymb_{\meshparam}}
    {\Tsymb_{#1}}
}
\newcommand{\Edgesymb}{e}
\newcommand{\Edge}[1][]{
  \ifthenelse{\equal{#1}{}}
    {\Edgesymb}
    {\Edgesymb_{#1}}
}
\newcommand{\EdgeH}[1][]{
  \ifthenelse{\equal{#1}{}}
    {\Edgesymb_{\meshparam}}
    {\Edgesymb_{\meshparam,#1}}
}
\newcommand{\NEdge}[1][]{
  \ifthenelse{\equal{#1}{}}
    {N_{\Edgesymb}}
    {N_{\Edgesymb({#1})}}
}
\newcommand{\Cellsymb}{K}
\newcommand{\Cell}[1][]{
  \ifthenelse{\equal{#1}{}}
    {\Cellsymb}
    {\Cellsymb_{#1}}
}
\newcommand{\CellH}[1][]{
  \ifthenelse{\equal{#1}{}}
    {\Cellsymb_{\meshparam}}
    {\Cellsymb_{\meshparam,#1}}
}
\newcommand{\areaSymb}{\mathcal{A}}
\newcommand{\CellArea}[1][]
{
  \ifthenelse{\equal{#1}{}}
    {\areaSymb_{\Cell}}
    {{A\left({#1}\right)}}
}
\newcommand{\CellHArea}[1][]
{
  \ifthenelse{\equal{#1}{}}
    {\areaSymb_{\CellH}}
    {\areaSymb_{\meshparam,#1}}
}

\newcommand{\NCell}[1][]{
  \ifthenelse{\equal{#1}{}}
    {N_{\Cellsymb}}
    {N_{\Cellsymb({#1})}}
}

\newcommand{\edgeLength}[1][]
{
  \ifthenelse{\equal{#1}{}}
  {\ABS{\Edge}}
  {\ABS{\Edge_{#1}}}
}
 
\newcommand{\edgeHLength}[1][]
{
  \ifthenelse{\equal{#1}{}}
  {\lengthSymb_{\EdgeH}}
  {\lengthSymb_{\meshparam,#1}}
}
\newcommand{\subVolume}[2]{D_{#1}^{#2}}

\newcommand{\Sourcesymb}{\mathbf{S}}
\newcommand{\Source}[1][]
{
  \ifthenelse{\equal{#1}{}}
    {\Sourcesymb}
    {\Sourcesymb_{#1}}
}
\newcommand{\bSource}[1][]
{
  \ifthenelse{\equal{#1}{}}
    {\overline{\Sourcesymb}}
    {\overline{\Sourcesymb}_{#1}}
}

\newcommand{\SourceEdge}[1][]{
  \ifthenelse{\equal{#1}{}}
    {\Sourcesymb_{ij}}
    {\Sourcesymb_{#1}}
}

\newcommand{\ConservVar}[2]{\mathbf{U}_{#1}^{{#2}}}

\newcommand{\bConservVar}[2]{\overline{\mathbf{U}}_{#1}^{{#2}}}
\newcommand{\FluxSWE}{\mathbf{F}}
\newcommand{\Flux}{\mathcal{F}}
\newcommand{\FluxLagrangian}{\underline{\underline{F}}}
\newcommand{\FluxEulerian}{\underline{\underline{G}}}

\newcommand{\FluxEdgesymb}{\mathbf{F}}
\newcommand{\FluxEdge}[1][]{
  \ifthenelse{\equal{#1}{}}
    {\FluxEdgesymb_{ij}}
    {\FluxEdgesymb_{#1}}
}

\newcommand{\numFlux}[1][]
{
  \ifthenelse{\equal{#1}{}}
  {\tilde{\fsymb}}
  {\tilde{\fsymb}_{#1}}
}

\newcommand{\FluxFuncNormSymbol}{\mathbf{F}}
\newcommand{\FluxFuncNorm}[1][]{
  \ifthenelse{\equal{#1}{}}
    {\FluxFuncNormSymbol^{\normalEdge}}
    {\FluxFuncNormSymbol^{\normalEdge}_{#1}}    
}

\newcommand{\JacobianSymbol}{\mathbf{A}}
\newcommand{\Jacobian}[1][]{
  \ifthenelse{\equal{#1}{}}
    {\JacobianSymbol}
    {\JacobianSymbol_{#1}}    
}
\newcommand{\EValSymbol}{\lambda}
\newcommand{\EVal}[1][]{
  \ifthenelse{\equal{#1}{}}
  {\EValSymbol}
  {\EValSymbol_{#1}}
}
\newcommand{\EVecSymbol}{\mathbf{r}}
\newcommand{\EVec}[2][]{
  \ifthenelse{\equal{#1}{}}
  {\EVecSymbol^{(#2)}}
  {\EVecSymbol^{(#2)}_{#1}}
}

\newcommand{\midPointEdge}[1][]
{
  \ifthenelse{\equal{#1}{}}
  {\midPoint_{\scriptscriptstyle\Edge}}
  {\midPoint_{\scriptscriptstyle\Edge[#1]}}
}
\newcommand{\gpPointEdgeDG}[1][]
{
  \ifthenelse{\equal{#1}{}}
  {\point_{\scriptscriptstyle\Edge}}
  {\point_{\scriptscriptstyle\Edge,#1}}
}
\newcommand{\midPointCell}[1][]
{
  \ifthenelse{\equal{#1}{}}
  {\midPoint_{\scriptscriptstyle\Cell}}
  {\midPoint_{\scriptscriptstyle\Cell[#1]}}
}

\newcommand{\speedRS}[2][]
{
  \ifthenelse{\equal{#2}{}}
  {S_{#2}}
  {S_{#2}^{#1}}
}

\newcommand{\ContSymbol}{C}
\newcommand{\Cont}[1][]{
  \ifthenelse{\equal{#1}{}}
  {\ContSymbol^{0}}
  {\ContSymbol^{#1}}
}
\newcommand{\Cinf}[1][]{
  \ifthenelse{\equal{#1}{}}
  {\ContSymbol^{\infty}}
  {\ContSymbol^{\infty}(#1)}
}
\newcommand{\SobSymbol}{W}
\newcommand{\HilbSymbol}{H}
\newcommand{\Sob}[2][]{
  \ifthenelse{\equal{#1}{}}
  {\HilbSymbol^{#2}}
  {\SobSymbol^{#2,#1}}  
}

\newcommand{\LspaceSymb}{L}
\newcommand{\Lspace}[1][]{
  \ifthenelse{\equal{#1}{}}
  {\LspaceSymb^{2}}
  {\LspaceSymb^{#1}}  
}

\newcommand{\TestSpSymbol}{\mathcal{V}}
\newcommand{\TestSpace}[1][]{
  \ifthenelse{\equal{#1}{}}
  {\TestSpSymbol({\SurfDomain})}
  {\TestSpSymbol_{#1}(\SurfDomain)}
}
\newcommand{\TestSpaceTime}[1][]{
  \ifthenelse{\equal{#1}{}}
  {\TestSpSymbol(\SurfDomain{}(\tempo))}
  {\TestSpSymbol_{#1}(\SurfDomain{}(\tempo))}
}

\newcommand{\AngleSymbol}{\theta}
\newcommand{\DevAngle}[1][]
{
  \ifthenelse{\equal{#1}{}}
  {\AngleSymbol}
  {\AngleSymbol_{\scriptscriptstyle{#1}}}
}
\newcommand{\relheightSymb}{\pi}
\newcommand{\relheight}[1][]
{
  \ifthenelse{\equal{#1}{}}
  {\relheightSymb_{\scriptscriptstyle{\SurfDomain}}}
  {\relheightSymb_{\scriptscriptstyle{#1}}}
}

\newcommand{\SolSymbol}{u}
\newcommand{\Sol}{\SolSymbol}

\newcommand{\BilinearStiffSymbol}{a}
\newcommand{\BilinearStiff}[3][]
{
  \ifthenelse{\equal{#1}{}}
  {\BilinearStiffSymbol(#2,#3)}    
  {\BilinearStiffSymbol_{#1}(#2,#3)}    
}
\newcommand{\BilinearAdvSymbol}{b}
\newcommand{\BilinearAdv}[3][]
{
  \ifthenelse{\equal{#1}{}}
  {\BilinearAdvSymbol(#2,#3)}    
  {\BilinearAdvSymbol_{#1}(#2,#3)}    
}
\newcommand{\BilinearMassSymbol}{m}
\newcommand{\BilinearMass}[3][]
{
  \ifthenelse{\equal{#1}{}}
  {\BilinearMassSymbol(#2,#3)}    
  {\BilinearMassSymbol_{#1}(#2,#3)}    
}
\newcommand{\BilinearReactSymbol}{c}
\newcommand{\BilinearReact}[3][]
{
  \ifthenelse{\equal{#1}{}}
  {\BilinearReactSymbol(#2,#3)}    
  {\BilinearReactSymbol_{#1}(#2,#3)}    
}

\newcommand{\TestSymbol}{v}
\newcommand{\Test}[1][]
{
  \ifthenelse{\equal{#1}{}}
  {\TestSymbol}  
  {\TestSymbol_{\scriptscriptstyle{#1}}}
}

\newcommand{\coeffSol}[2][]
{
  \ifthenelse{\equal{#1}{}}
  {\Sol_{#2}}  
  {\Sol_{#2}^{\scriptscriptstyle{#1}}}
}

\newcommand{\nodalBasis}[2][]
{
  \ifthenelse{\equal{#1}{}}
  {\varphi_{\scriptscriptstyle{#2}}}  
  {\varphi_{\scriptscriptstyle{#2}}^{\scriptscriptstyle{#1}}}
}

\newcommand{\nNodes}[1][]
{
  \ifthenelse{\equal{#1}{}}
  {N}
  {N_{\scriptscriptstyle{#1}}}
}

\newcommand{\TestApprox}[1][]
{
  \ifthenelse{\equal{#1}{}}
  {\TestSymbol_{\scriptscriptstyle{\meshparam}}}  
  {\TestSymbol_{\scriptscriptstyle{\meshparam,#1}}}
}

\newcommand{\viscositySymb}{\nu}
\newcommand{\viscosity}[1][]
{
  \ifthenelse{\equal{#1}{}}
  {\viscositySymb}
  {\viscositySymb_{\scriptscriptstyle{#1}}}
}

\newcommand{\ResidualSymbol}{R}
\newcommand{\Residual}[1][]
{
  \ifthenelse{\equal{#1}{}}
  {\ResidualSymbol}
  {\ResidualSymbol_{#1}}
}

\renewcommand{\xcg}{x}
\renewcommand{\ycg}{y}
\renewcommand{\zcg}{z}
\setlist[description]{font=\normalfont\bfseries}
\renewcommand{\paragraph}{\subsubsection*}

\usepackage{setspace}
\usepackage[a4paper,left=2.cm,right=2.cm, bottom=2.5cm, top=2.5cm]{geometry}
\title[Lagrangian-Eulerian FV for ISWE on locally discontinuous
  beds]{A robust Lagrangian-Eulerian Finite Volume scheme for intrinsic
  shallow water equations with wetting-and-drying transitions on
  locally discontinuous beds}

\author{Eduardo Abreu$^1$ \and Elena Bachini$^2$ \and John P\'erez$^3$ \and Mario Putti$^2$}
\address{
  $^1$ Department of Applied Mathematics,
              University of Campinas - IMECC, Brazil\\
  $^2$ Department of Agronomy, Food, Natural
    resources, Animals and Environment, University of Padua, Italy \\
    $^3$ ITM-University Institution, Medellin,
           Colombia
}

\begin{document}

\maketitle

\begin{abstract}
      In this work, we extend the Intrinsic Lagrangian–Eulerian Finite
    Volume method for the Intrinsic Shallow Water Equations (ISWE) to
    account for localized topographic discontinuities that influence 
    flow patterns. The ISWE formulation is modified to incorporate these
    bottom perturbations without altering the global structure of the
    system. The numerical discretization is based on a Lagrangian–Eulerian
    Finite Volume scheme within an intrinsic geometric framework and
    consists of a Lagrangian evolution step followed by an Eulerian
    remapping onto a fixed reference grid.
    The combination of the intrinsic nature of the ISWE formulation
    and the dual Lagrangian–Eulerian structure of the proposed
    discretization yields a numerical method that naturally satisfies
    the three key properties essential for stability, accuracy, and
    robustness: well-balancing for preservation of the
    ``lake-at-rest'' equilibrium over irregular bathymetry;
    wetting-and-drying, enabling robust tracking of moving fluid
    boundaries; and hyperbolicity guaranteeing consistency with the
    system wave propagation and stability constraints. Numerical
    experiments on curved two-dimensional surfaces with localized
    bottom discontinuities demonstrate the accuracy and robustness of
    the method, including accurate long-time velocity computations in
    near-dry conditions and transitions between wet and dry states.
\end{abstract}

  

\section{Introduction}
\label{sec:intro}

The Shallow Water Equations (SWE) are a reduced dimensional model of
the Navier-Stokes equations derived under the assumption that
free-surface fluid motion develops essentially on a topologically
two-dimensional domain. This situation occurs typically when
perturbations of the fluid depth are small with respect to the
wavelengths of the predominant flow. SWE are used regularly for
small bed slopes for which the direction of gravity can be considered
as the normal to the bottom surface~\cite{book:Toro2001,book:Vreugdenhil2013}.
When fluid flow takes place on a general three-dimensional surface,
this approximation is no longer valid and the bottom normal must be
used to measure fluid depth. In this case, a more geometrical
viewpoint allows for automatic consideration of surface
characteristics such as slope and curvatures, and overcome classical
slope limitations of the more common methodologies. This geometric
approach, pioneered in the field of SWE
by~\cite{art:SavageHutter1989}, has been formalized
by~\cite{art:Bouchut2003,art:Bouchut2004} using an extrinsic point of
view, and, more recently, by \cite{art:BP20} (see also
\cite{art:ABREU2023iswe}) employing a fully intrinsic formulation,
called ISWE (Intrinsic SWE). ISWE constitutes the starting point of
our work.

From the numerical perspective, one of the major problems of the ISWE
approach is that the ensuing hyperbolic fluxes are non-autonomous,
thus leading to well-posedness difficulties in the associated
Riemann-problems, characterized by discontinuous
fluxes~\cite{art:Audusse2004,art:andreianov2010,art:andreianov2011,art:andreianov2015}. To resolve this obstruction, recently we
proposed the use of a Lagrangian-Eulerian approach for the
finite-volume discretization of ISWE~\cite{art:ABREU2023iswe}. In
essence, the Lagrangian-Eulerian framework consists in a first
Lagrangian step, where the solution evolves following a Lagrangian
path, followed by a second Eulerian step, where the conserved
quantities are remapped back to the original finite volume mesh using
locally conservative cell
averaging~\cite{art:abreu17,abreu2018weak,art:abreu19,art:abreu20,art:abreu22,zbMATH07987763}. These
two building blocks do not require the formulation of any (local or
global) Riemann problem and thus are not subjected to the above
mentioned difficulties that arise in the case of non-autonomous
fluxes.

One of the main challenges in the modeling of gravity driven motion of
fluids is the accurate approximation of complex flow patterns that
arise in the presence of localized obstacles and discontinuities in
the topography.
In general, more complex and non-localized discontinuities either do
not fall within the shallow water hypotheses or can be approximated
with sufficient accuracy by regularizing the bottom geometry. Thus, we
focus on the localized cases.
Theoretical and numerical difficulties arise in the presence of
locally discontinuous bottoms: preservation of the steady-state of
still water, entropy inequalities and uniqueness of solution,
preservation of the non-negativity of the water depth, and loss of
hyperbolicity (see, e.g.,
\cite{zbMATH01846809,zbMATH07935471,diaz2023discussion,krvavica2018analytical,ALEKSYUK2022110801,art:xia20,abreu2024riemann}
and the references cited therein).  Indeed, the source term of the
SWEs contains a nonconservative product of two possibly discontinuous
functions, making the problem non-strictly hyperbolic and
nonconservative. The typical solution is to augment the system by
imposing the artificial condition that the time derivative of the
bottom is zero, thus allowing the emergence of the full wave
structure, including the linearly degenerate characteristics arising
from the nonconservative product~\cite{book:Bouchut2004,GUERRA2004438,lefloch2007riemann,PARES2019344,ALEKSYUK2022110801,VALIANI2024104735}. In addition, these discontinuous forcings cause
the failure of the $C$-property, or well-balance, in the scheme,
providing difficulties in the definition of appropriate Riemann
solvers and thus requiring ad-hoc well-balance techniques and
path-conservative schemes~\cite{book:Bouchut2004,dumbser2011simple}.

In this work, we overcome these difficulties by careful adaptation of
both the continuous ISWE formulation and the numerical
Lagrangian-Eulerian approach. The set of adaptations constitute the
major advancements we bring to our ISWE-based Lagrangian-Eulerian
approach presented in~\cite{art:ABREU2023iswe}.
In more details, first, we relax the limitation of the geometric
formulation by allowing localized possibly discontinuous
perturbations, i.e., obstacles defined on top of the (global) bottom
topography.
Topographic discontinuities are incorporated by local numerical
smoothing of the bottom surface to provide sound mathematical
connotation to the nonconservative product. This allows for a
straight-forward implementation into the Lagrangian-Eulerian framework
of these local perturbations, which are directly included in the
source terms without modifying the global structure of the
(hyperbolic) system.
Second, by exploiting the properties of ISWE and by a careful
decomposition of the flux function between the Lagrangian and Eulerian
steps, we guarantee well-balance and the satisfaction of the
C-property of the solution. In fact, the Lagrangian evolution is
defined only by the advective terms, while the remaining flux terms
are conserved by means of the Eulerian projection, thus allowing
automatic well-balance and wetting-and-drying.
Finally, we provide a novel and more thorough interpretation of our
previous results presented in~\cite{art:ABREU2023iswe} from the
perspective of ``control of hyperbolicity'', by careful geometric
identification of entropy solutions within a Jacobian-free and
Riemann-solver-free framework. Note that, these points are not mere
re-descriptions but represent substantial key features of our
intrinsic Lagrangian-Eulerian approach.
We would like to specifically summarize here these significant
novelties as follows:
\begin{description}
\item[Jacobian-Free and Riemann-Solver Free Framework] The scheme
  operates on an \emph{unconventional paradigm} of no-flow curves,
  which serve as geometric alternatives to classical
  characteristics. The evolution step is performed without the use of
  spectral information (eigenpairs) from the flux function Jacobian
  matrix. Consequently, it bypasses the need for costly Riemann
  solvers, making it highly efficient and accurate for the ISWE
  non-autonomous fluxes.  Stability (CFL) conditions and system
  positivity properties are easily shown.
\item[Geometric Identification of Entropy Solutions] Instead of
  relying on field-by-field characteristic decompositions, our
  framework based on the no-flow curves is able to identify the
  entropy solution by tracking local wave speeds along space-time
  trajectories. This ensures physical consistency through the direct
  interaction between the flux function and conserved variables,
  resulting in a self-contained solver with a compact spatial
  footprint.
\item[Inherent Well-Balancing and ``Lake-at-Rest'' Stability] The
  scheme is specifically designed to handle irregular bottom surfaces
  by naturally maintaining the ``lake-at-rest'' equilibrium. This is
  achieved through a precise decomposition of the flux function during
  the two-step building block process, the Lagrangian evolution and
  the Eulerian remapping, as enabled by the intrinsic approach of
  ISWE.
\item[Robust Handling of Local Discontinuities and Moving Fronts] The
  formulation is uniquely equipped to manage localized obstacles and
  topographic jumps via local numerical smoothing. Furthermore, ISWE
  allows to naturally capture ``wetting-and-drying'' moving boundaries
  at the fluid interface, thus eliminating the need for ad-hoc
  front-tracking or thin-layer algorithms. This is the fundamental
  ingredient that allows the inclusion of localized bottom
  perturbations without requiring major modifications to the global
  system.
\end{description}
We conclude the paper by presenting in-depth numerical results to
support the accuracy and robustness of our strategy, leaving a more
rigorous and in-depth mathematical analysis of the theoretical
convergence properties to future studies.
Finally, we would like to mention that the approach can be extended to
general polygons, with particular reference to quadtree meshes, which
are typically used for their efficiency and parallelization
properties~\cite{gatti2023,GLF21}.

The paper is organized as follows. In \autoref{sec:iswe}
and \autoref{sec:LEfv}, we present the extended ISWE model in the presence
of a localized discontinuous topography and we summarize the
Lagrangian-Eulerian Finite Volume scheme. %
\autoref{sec:results} reports numerical results on various
surfaces considering the presence of local discontinuities. Extra
details are presented in the appendix sections: \autoref{app:quadrule}
contains the calculations of the quadrature rules and \autoref{app:ref}
the experimental convergence assessment on mesh refinements.

\section{Extended model: ISWE in the presence  of  localized
  discontinuities}
\label{sec:iswe}

Consider a regular surface $\SurfDomain\subset\REAL^3$ equipped
with a set of local coordinates $\sv=(\xcl,\ycl)$ and the associated
metric tensor $\First$. We assume that the local reference frame
$\{\vecBaseCCcv[1](\Point),\vecBaseCCcv[2](\Point)\}$ at each point $\Point\in\SurfDomain$
is orthogonal, and thus the metric tensor is diagonal and can be
defined by:
\begin{equation*}
  \First = \begin{pmatrix}
    \NORM{\vecBaseCCcv[1](\Point)}^2 & 0\\
    0 &\NORM{\vecBaseCCcv[2](\Point)}^2 \\
  \end{pmatrix}
  =
  \begin{pmatrix}
    \metrcoefH{1}^2 & 0\\
    0 &\metrcoefH{2}^2  \\
  \end{pmatrix}\,.
\end{equation*}
Note that, up to a local change of reference frame, these assumptions
are satisfied and do not constitute a limitation to the problem. The
definition of a local reference frame relies on the further assumption
that the surface in regular ``enough'' and that an atlas of local
parametrizations exists (the surface needs to have $\Cont[1,1]$
regularity).

In this setting, the system of intrinsic shallow water equations
(ISWE)~\cite{art:BP20} is written as:
\begin{equation}
  \DerPar{\ConservVar{}{}}{\Time}
  + \DivSurf\FluxSWE(\sv,\ConservVar{}{}) +
  \Source(\sv,\ConservVar{}{}) =0
  \qquad \mbox{ in } \IT\times\SurfDomain\,,
  \label{eq1swe}
\end{equation}
with initial condition $\ConservVar{}{}(0,\sv)=\ConservVar{0}{}$,
where
$\ConservVar{}{} =\left[ \Depth,\,\Qcomp[1]{},\,\Qcomp[2]{}
\right]^{\T}$ is the conservative variable, with $\Depth$ the fluid
depth along the normal direction to the bottom surface and
$\Qcomp[i]{}$, $i=1,2$, the averaged discharges in the two tangential
directions.  The symbol $\DivSurf$ indicates the intrinsic divergence
operator, defined with respect to the metric tensor $\First$.
The hyperbolic flux tensor $\FluxSWE$ and the source term $\Source$ 
are defined by:
\begin{equation*}
  \label{eq:flux}
  \FluxSWE(\sv,\ConservVar{}{}) =
  \begin{bmatrix}
    \Qcomp[1]{}   &  \Qcomp[2]{} \\[0.7em]
    \displaystyle\DerPar{(\Qcomp[1]{})^2}{\Depth} +
    \frac{\grav \Depth^2}{2 \metrcoefH{1}^2} \DerPar{\zcg}{\zcl}
    &
    \displaystyle\frac{\Qcomp[1]{}\Qcomp[2]{}}{\Depth}  \\[0.7em]
    \displaystyle\frac{\Qcomp[1]{}\Qcomp[2]{}}{\Depth}
    &
    \displaystyle\frac{(\Qcomp[2]{})^2}{\Depth} +
    \frac{\grav \Depth^2}{2 \metrcoefH{2}^2} \DerPar{\zcg}{\zcl} 
  \end{bmatrix}
\end{equation*}
and
\begin{equation*}
  \label{eq:source}
  \Source(\sv,\ConservVar{}{}) =
  \begin{bmatrix}
    0 \\[0.7em]
    \displaystyle
    \frac{\grav \Depth^2}{2 \metrcoefH{1}^2}\DerPar{}{\xcl}
    \left(\DerPar{\zcg}{\zcl}\right)
    +\frac{\grav\Depth}{\metrcoefH{1}^2}\DerPar{\zcg}{\xcl}
    - \frac{1}{\density}[\DivSurf \Tens[sw]]^{(1,\cdot)}
    -\frac{\BottomFriction[1]}{\density} \\[0.7em]
    \displaystyle
    \frac{\grav \Depth^2}{2
      \metrcoefH{2}^2}\DerPar{}{\ycl}
    \left(\DerPar{\zcg}{\zcl}\right)
    +  \frac{\grav\Depth}{\metrcoefH{2}^2}\DerPar{\zcg}{\ycl}
    - 	\frac{1}{\density}[\DivSurf  \Tens[sw]]^{(2,\cdot)}
    - \frac{\BottomFriction[2]}{\density} 
  \end{bmatrix},
\end{equation*}

As stated above, the definition of a local reference frame relies on
the assumption that the surface in regular ``enough''. Thus, in
principle, it is not possible to directly consider a discontinuous
surface $\SurfDomain$ as domain of the ISWE model.
Our strategy is to consider localized discontinuities as local
perturbations of the bottom topography and to take these into account,
without changing the geometric regularity assumption of the underlying
surface, exploiting the fact that our Lagrangian evolution does not
involve Riemann problems. In other words, the ISWE framework rests on
a globally regular bottom surface, with the effects of the local
discontinuity confined in a modification of the source term. This
approach does not introduce extra equations and does not modify the
hyperbolic nature of the system, contrary to the approach typically
adopted in these cases (e.g.,~\cite{COZZOLINO20111413,art:Pares2006}).
Explicitly, starting from the system~\eqref{eq1swe}, we add a new term
$\mathbf{Z}$ containing the gradient of the localized extension to the
bottom $\obstacle$, i.e., the obstacle on top of the regular bottom
$\SurfDomain$:
\begin{equation*}
  \DerPar{\ConservVar{}{}}{\Time} +
  \DivSurf \FluxSWE(\sv,\ConservVar{}{}) +
  \Source(\sv,\ConservVar{}{}) + \mathbf{Z} = 0\,.
  \label{eq1}
\end{equation*}
It is then straightforward to define a ``new'' source term $\bSource$
putting together the contributions of $\Source$ and $\mathbf{Z}$:
\begin{multline}
  \bSource(\sv,\ConservVar{}{}) =  \Source(\sv,\ConservVar{}{}) +
  \mathbf{Z}
  =\begin{bmatrix}
  0 \\[0.7em]
  \displaystyle
  \frac{\grav \Depth^2}{2 \metrcoefH{1}^2}\DerPar{}{\xcl}
  \left( \DerPar{\zcg}{\zcl}\right)
  + \frac{\grav \Depth}{\metrcoefH{1}^2}\DerPar{\zcg}{\xcl}
  -\frac{1}{\density}[\DivSurf \Tens[sw]]^{(1,\cdot)}
  -\frac{\BottomFriction[1]}{\density} \\[1em]
  \displaystyle
  \frac{\grav \Depth^2}{2 \metrcoefH{2}^2}\DerPar{}{\ycl}
  \left(\DerPar{\zcg}{\zcl}\right)
  +\frac{\grav \Depth}{\metrcoefH{2}^2}\DerPar{\zcg}{\ycl}
  -\frac{1}{\density}[\DivSurf \Tens[sw]]^{(2,\cdot)}
  -\frac{\BottomFriction[2]}{\density} 
  \end{bmatrix} +
  \begin{bmatrix}
    0 \\[0.7em]
    \displaystyle\frac{\grav \Depth}{\metrcoefH{1}^2} \DerPar{\obstacle}{\xcl} \\[1.1em]
    \displaystyle \frac{\grav \Depth}{\metrcoefH{2}^2}\DerPar{\obstacle}{\ycl}
  \end{bmatrix}
  \label{IBDMS}
\end{multline}
and to re-define the ``extended'' ISWE model by:
\begin{equation}
  \DerPar{\ConservVar{}{}}{\Time} + \DivSurf \FluxSWE(\sv,\ConservVar{}{}) + \bSource(\sv,\ConservVar{}{})  = 0\,.
  \label{eq:extendedISWE}
\end{equation}
We start from this system and present our Lagrangian-Eulerian Finite
Volume scheme, already introduced in ~\cite{art:ABREU2023iswe}.

\section{Lagrangian-Eulerian Finite Volume scheme with flux
  decomposition}
\label{sec:LEfv}

In this section, we report a summary of our intrinsic
Lagrangian-Eulerian approach, introduced in~\cite{art:ABREU2023iswe},
and applied to the extended ISWE model~\eqref{eq:extendedISWE}. We
will see that the extra term $\mathbf{Z}$ introduced to handle local
discontinuities does not present any additional difficulty in our
Lagrangian-Eulerian setting. We also add a flux decomposition strategy
that guarantees the well-balance property at the discrete level. We
devise specific numerical example to support our procedure.

Our starting point is the definition of a \emph{generalized time-space
  divergence form} for the balance law~\eqref{eq:extendedISWE}. To do
so, we define an extension of the metric tensor in time-space
coordinates, i.e.:
\begin{equation*}
  \mathbb{G} :=
  \left(
  \begin{array}{rr}
    1\; \vline&  \\
    \hline
    \vline& \First 
  \end{array}
  \right)\,,
\end{equation*}
with the associated reference frame
$\langle\normalvec, \vecBaseCCcv[1],\vecBaseCCcv[2]\rangle$, where
$\vecBaseCCcv[1],\vecBaseCCcv[2]$ span the tangent plane and are used
to represent the spatial local coordiantes $\sv$, while $\normalvec$
is the local normal to the surface used to extend the tangent frame to
a basis of $\REAL^3$ and represents the time coordinate.
We denote these extended coordinates, i.e. $(\tempo,\sv)$, the
Lagrangian-Eulerian local coordinates.

In this time-space framework, the system~\eqref{eq:extendedISWE} can be
re-written in divergence form as:
\begin{align}
  \label{eq:eqSTdiv}
  \DivSurfST\Flux=-\bSource,
\end{align}
where $\DivSurfST$ is the intrinsic divergence with respect to the
extended metric $\mathbb{G}$, and $\Flux$ is an extended (symmetric)
flux function defined as: 
\begin{equation}
  \label{eq:fluxtensor}
  \Flux(\sv,\ConservVar{}{}) = 
  \FluxLagrangian +   \FluxEulerian = 
  \begin{bmatrix} 
    \Depth &\Qcomp[1]{} & \Qcomp[2]{} \\[0.7em]
    \Qcomp[1]{}
    & \displaystyle
    \frac{(\Qcomp[1]{})^2}{\Depth}
    & \displaystyle
    \frac{\Qcomp[1]{}\Qcomp[2]{}}{\Depth}  \\[0.7em]
    \Qcomp[2]{}
    & \displaystyle
    \frac{\Qcomp[1]{}\Qcomp[2]{}}{\Depth}
    & \displaystyle
    \frac{(\Qcomp[2]{})^2}{\Depth}
  \end{bmatrix} +
  \begin{bmatrix} 
    0 & 0 & 0 \\[0.7em]
    0
    & \displaystyle
    \frac{\grav
      \Depth^2}{2\metrcoefH{1}^2}\DerPar{\zcg}{\zcl}
    &     0 \\[0.7em]
    0
    &   0
    & \displaystyle
    \frac{\grav
      \Depth^2}{2\metrcoefH{2}^2}\DerPar{\zcg}{\zcl}
  \end{bmatrix}.
\end{equation}
Note that, we split the extended flux function $\Flux$ in the sum of
the ``Lagrangian flux'' $\FluxLagrangian$ and the ``Eulerian flux''
$\FluxEulerian$. This physics-based Lagrangian-Eulerian decomposition,
the second main novelty of this work, helps in maintaining the
well-balance at the discrete level.

The Lagrangian-Eulerian numerical solver, as adapted to ISWE
in~\cite{art:ABREU2023iswe}, can be summarized as follows. First, the
bottom surface $\SurfDomain$ is partitioned into non-overlapping
regions or cells, $\Region$, such that their union is exactly
$\Closure{\SurfDomain}$.
Then, the discretization proceeds by testing \autoref{eq:eqSTdiv} using
the characteristic function of the so-called space-time \emph{Discrete
control volumes} (DCV).  A discrete control volume is defined as the
time evolution of $\Region$ along no-flow surfaces,
and it is represented by the region
$\subVolume{}{n}\subset\IT\times\SurfDomain$ with boundary
$\partial \subVolume{}{} = \Region[n] \cup \mathcal{LS}^n\cup
\bRegion[n+1]$ given by the union of $\Region[n]$ and $\bRegion[n+1]$,
which constitute the \emph{Discrete Entry} and \emph{Discrete Exit}
surfaces of the Lagrangian flux $\FluxLagrangian$ at times $\tempo^n$
and $\tempo^{n+1}$, and the lateral surface $\mathcal{LS}^n$, the
so-called \emph{Discrete no-flow surface}. The latter is defined by
the condition that the normal Lagrangian flux be zero, i.e.:
$\scalprodSurfST{\FluxLagrangian}{\vec{n}}=
\FluxLagrangian\STFirst\vec{n}=0$, where here $\vec{n}$ denotes the
outer normal to the boundary $\partial\mathcal{LS}^n$.
Then, applying the divergence theorem to \autoref{eq:eqSTdiv}, we
obtain for every $DCV$ $\Region[n]$: 
\begin{align}
  \int_{\partial \subVolume{}{n}} \scalprodSurfST{\Flux}{\vec{n}}
  &=-\int_{\subVolume{}{n}} \Source - \int_{\subVolume{}{n}} \mathbf{Z}  
  =: I(\bSource) \,.
  \label{eq4}
\end{align}
Note that, $\vec{n}=-\normalvec$ in the entry region $\Region[n]$,
and $\vec{n}=\normalvec$ in the exit region
$\bRegion[n+1]$. Moreover, when written in the time-space local
coordinate system, $\normalvec=[1,0,0]^{\T}$. As a consequence, the
left-hand-side of equation~\eqref{eq4} can be written as:
\begin{align*}
  -\int_{\Region[n]} \ConservVar{}{}(\tempo^n,\xcl,\ycl)
  + \int_{\mathcal{LS}^n} \scalprodSurfST{\FluxLagrangian+\FluxEulerian}{\normalEdge} 
  + \int_{\bRegion[n+1]}\ConservVar{}{}(t^{n+1},\xcl,\ycl) &=
  I(\bSource)\,,
\end{align*}
where $\normalEdge$ is the normal to $\mathcal{LS}^n$ written in the
local coordinates $\sv=(\xcl,\ycl)$ at each time $t$. The term
containing $\FluxLagrangian$ will defines our Lagrangian no-flow
condition (see \autoref{hypnoflow}), while the term containing
$\FluxEulerian$ will be computed together with $I(\bSource)$ in the
Eulerian step. This decomposition guarantees a discrete well-balance
property and hinders blow-ups and oscillations that may occur in
wetting-and-drying areas for small or zero depths.

The no-flow condition through the lateral surface $\mathcal{LS}^n$ for
the Lagrangian flux can can be enforced by imposing:
\begin{equation}
  \int_{\mathcal{LS}^n} \scalprodSurfST{\FluxLagrangian}{\normalEdge} = 0\,.
  \label{hypnoflow}
\end{equation}
Hence, the evolution relation for the state variable is given by:
\begin{equation}
  \int_{\bRegion[n+1]} 
  \ConservVar{}{}(t^{n+1},\xcl,\ycl) 
  \, = \int_{\Region[n]} 
  \ConservVar{}{}(t^{n},\xcl,\ycl) 
  \, + I(\bSource) -  \int_{\mathcal{LS}^n} \scalprodSurfST{\FluxEulerian}{\normalEdge}\,.
  \label{eq6}
\end{equation}
Putting these two conditions together, we obtain an equivalent
formulation for our initial system~\eqref{eq1swe}.
In summary, in $[\tempo^{n},\tempo^{n+1}]$, the following two problems
are equivalent in the weak (distributional) sense:
\begin{equation*}
  {\bf (I)}
  \quad
  \ConservVar{}{}(\tempo^{n+1},\sv)=
  \ConservVar{}{}(\tempo^{n},\sv)
  -
  \int_{\tempo^{n}}^{\tempo^{n+1}}
  \left[
    \DivSurf\FluxSWE(\sv,\ConservVar{}{})
    + \bSource(\sv,\ConservVar{}{})
  \right]\Diff\tempo
\end{equation*}
and
\begin{equation*}
  {\bf (II)}
  \begin{cases}
    \displaystyle \int_{\mathcal{LS}^n} \scalprodSurfST{\FluxLagrangian}{\normalEdge} = 0 \\[1em]
    \displaystyle  \frac{1}{\CellArea[{{\bRegion[{n+1}]}}]}
    \int_{\bRegion[n+1]} 
    \ConservVar{}{}(t^{n+1},\xcl,\ycl) 
    = \frac{1}{\CellArea[{{\bRegion[{n+1}]}}]} \int_{\Region[n]} 
    \ConservVar{}{}(t^{n},\xcl,\ycl) 
    \,\\
    \phantom{ \frac{1}{A} \int_{\bRegion[n+1]} U }
    \displaystyle
    - \frac{1}{\CellArea[{{\bRegion[{n+1}]}}]}
    \left[
    \int_{\subVolume{}{n}} (\Source +  \mathbf{Z})
    + \int_{\mathcal{LS}^n}
    \scalprodSurfST{\FluxEulerian}{\normalEdge}
    \right]
    \\[0.7em]
  \end{cases}
\end{equation*}
where $\CellArea[{\bRegion[n+1]}]$ denotes the area of the evolved
region $\bRegion[n+1]$.
Problem \textbf{(II)} represents the geometrically intrinsic
Lagrangian-Eulerian building block, which will form the basis for our
numerical scheme.
Recalling that the regions $\Region[n],\bRegion[n+1]$,
$\subVolume{}{n}$, and $\mathcal{LS}^n$ are defined in a time-space
reference frame, we can understand the Cauchy problem \textbf{(II)} as
the Lagrangian evolution step followed by the remapping onto the
original region $\Region[n]$ obtained by imposing in an Eulerian sense
the local mass conservation (Eulerian remapping step).

\paragraph{Discrete Lagrangian-Eulerian scheme.}
We consider a fixed (in time) surface triangulation $\Triang(\SurfDomain)$
of the bottom surface $\SurfDomain$, formed by the
union of non-intersecting geodesic cells $\Cell[\ell]$ with vertices
on $\SurfDomain$ and such that
$\Triang(\SurfDomain)=\cup_{\ell=1}^{\NCell}\Cell[\ell] =\Closure{\SurfDomain}$
and $\Edge[\ell m]=\Cell[\ell]\cap\Cell[m]$ is a geodesic edge.
The triangulation is characterized by a mesh parameter
$\meshparam=\max_\ell \{\diam(\Cell[\ell])\}$.
We assume that there is a one-to-one correspondence between a cell
$\Cell[\ell]$ and the flat cell defined by projecting its vertices on
the tangent plane anchored at the cell mid-point
$\TanPlane[\midPoint_\ell]\SurfDomain$. Note that this is always
possible, up to a rescaling of the cell $\Cell[\ell]$. This piecewise
linear approximation of the surface enables us to work on (local) flat
cells with a geometrical error proportional to $\meshparam^2$.
Together with the space discretization, we consider a time
discretization $\{\tempo^n\}$ of the time interval $\IT$, where we define a
non-uniform time-step $\Dt=\tempo^{n+1}-\tempo^n$.
At each time $\tempo^n$, the cell $\Cell[\ell]$ coincides with the inflow
region $\Region[n]_\ell$, namely the triangulation $\Triang(\SurfDomain)$
is our working grid.

On each cell $\Region[n]_\ell$, we define the approximate 
solution  $\ConservVar{\ell}{n}$ at time $\tempo^n$ by taking the
mean value over the region $\Region[n]_{\ell}$, i.e.:
\begin{equation*}
  \ConservVar{\ell}{n} := \frac{1}{\CellArea[{\Region[n]_\ell}]} 
  \int_{\Region[n]_\ell} \ConservVar{}{}(\tempo^n,\sv)\,.
\end{equation*}
Then, $\ConservVar{\ell}{n+1}$ over the region
$\bRegion[n+1]_\ell$, is obtained after a remapping procedure
(Eulerian remapping step) over the original grid
$\Cell[\ell]=\Region[n]_\ell$ in the corresponding control volume.
If we denote by $\bConservVar{\ell}{n+1}$ the approximate value at time
$\Time^{n+1}$ obtained from the Lagrangian step~\eqref{eq6}
(intermediate result), we can write:
\begin{align*}
  \bConservVar{\ell}{n+1} &= \frac{1}{\CellArea[{\bRegion[n+1]_\ell}]}
  \int_{\bRegion[n+1]_\ell} \ConservVar{}{}(t^{n+1},\sv) \nonumber\\
  &= 
  \frac{\CellArea[{\Region[n]_\ell}]}{\CellArea[{\bRegion[n+1]_\ell}]} 
  \left[\ConservVar{\ell}{n} + \frac{1}{\CellArea[{\Region[n]_\ell}]}
    I(\bSource)
    - \int_{\mathcal{LS}^n} \scalprodSurfST{\FluxEulerian}{\normalEdge}\right]\,,
  \label{eq7}
\end{align*}
and the final updated variable $\ConservVar{\ell}{n+1}$ is obtained by
a simple geometric projection of $\bConservVar{\ell}{n+1}$ onto the
surface cell $\Region[n]_\ell$ (for more information
see~\cite{art:ABREU2023iswe}).

To complete our scheme we need a discrete version of the first
equation of problem \textbf{(II)}, i.e., we need to approximate the
no-flow condition~\eqref{hypnoflow}.
For any time $\tau\in\IT[n]$, the equation can be re-written
separating space and time and evaluating the integral numerically
using the mid-point rule on each edge of $\Region[n]_\ell$ and summing
over all edges. Note that a generic face of $\mathcal{LS}^n$ can be
seen as the evolution in time of a cell edge $\Edge[\ell m]$. Our
mid-point quadrature rule approximates the integral over this face by
a line integral in the time direction of the edge midpoint
$\midPointEdge$. This yields:
\begin{multline*}
  0=
  \int_{\tempo^n}^{\tau} \int_{\Edge[\ell m](\tempo)}
  \scalprodSurfST{\FluxLagrangian(\sv,\ConservVar{}{}(\tempo,\sv))}{\normalEdge(\tempo,\sv)}
  \Diff \sv\;\Diff \tempo
  \approx \int_{\tempo^n}^{\tau} \edgeLength[\ell m](\tempo)
  \scalprodSurfST{\FluxLagrangian(\midPointEdge,\ConservVar{}{}(\tempo,\midPointEdge))}{\normalEdge(\tempo,\midPointEdge)}
  \Diff \tempo\,,
\end{multline*}
where $\edgeLength[\ell m]$ and $\midPointEdge$ are the length and the
mid-point of $\Edge[\ell m](\tempo)$, respectively.
The time integral is approximated by means of explicit Euler, i.e.,
using the known quantities at time $\tempo^n$. This allows an explicit
characterization of the approximate no-flow curves
$\sigma_{\ell m}\subset \mathcal{LS}_{\ell m}^n$ that starts from
$\midPointEdge$ on each edge.
The extra condition we need to impose is that $\sigma_{\ell m}$ must
be orthogonal to the local normal to the cell edge
$\normalEdge_{\ell m}(\tempo)=\normalEdge(\tempo,\midPointEdge)$.
This procedure yields a nonlinear system of ODEs that describes
$\sigma'_{\ell m}$ as a function of the ISWE flux, with initial
condition given by
$\sigma_{\ell m}(\tempo^n) =
[\tempo^n,\midPointEdge^1(\tempo^n),\midPointEdge^2(\tempo^n)]^{\T}$
and flow defined implicitly at time $\tempo^n$ by
$\scalprodSurfST{\normalEdge_{\ell m}(\tempo^n)}{\sigma_{\ell m}'(\tempo^n)}=0$.
Following this approach, the position of the evolved edge mid-point at
$\tempo\in\IT[n]$ can be calculated as:
\begin{equation*}
  \sigma_{\ell m}(\tempo) = \sigma_{\ell m}(\tempo^n) +
  (\tempo-\tempo^n)\sigma_{\ell m}'(\tempo^n)\,,
\end{equation*}
with
$\sigma_{\ell m}'=[1,(\sigma^1_{\ell m})',(\sigma^2_{\ell m})']^{\T}$
obtained by
solving the following linear system:
\begin{equation*}
  \begin{cases}
    \scalprodSurfST{\normalEdge_{\ell m}(\tempo^n)}{\sigma_{\ell
        m}'(\tempo^n)}=0 \\[0.7em]
    \scalprodSurfST{\FluxLagrangian(\midPointEdge,\ConservVar{\midPointEdge}{n})}{\normalEdge(\tempo^n,\midPointEdge)}=\mathbf{0}
  \end{cases}
\end{equation*}
We refer again to~\cite{art:ABREU2023iswe} for more details.

\begin{remark}[Example of no-flow curve in the case of a Cartesian
    grid]
  In the special case of an edge with unit tangent vector in
  $\midPointEdge(\tempo^n)$ that is aligned with $\vecBaseCCcv[2]$ and normal vector 
  $\normalEdge=\left[-1,\displaystyle\frac{1}{(\sigma^{1})'(\tempo^n)},0\right]^{\T}$, we 
  obtain explicitly $\sigma'=[1,(\sigma^1)',0]^{\T}$ by
  solving the following linear system:
  \begin{equation*}
    0 = \scalprodSurfST{\FluxLagrangian(\midPointEdge)}{\normalEdge(\tempo^n,\midPointEdge)}\Big|^{\beta}=
    \FluxLagrangian^{\beta}(\midPointEdge)\, \STFirst(\midPointEdge)\,\begin{bmatrix}-1\\
    \displaystyle\frac{1}{(\sigma^{1\beta})'(\tempo^n)}\\
    0
    \end{bmatrix}
    \,,
    \label{odecp1}
  \end{equation*}
  for each component $\beta=1,2,3$.
  Thus:
  \begin{equation*}
    (\sigma^{1\beta})'(\tempo^n) =
    \FluxLagrangian^{(\beta,2)}\metrcoefH{1}^2/\ConservVar{}{\beta, n} \,,
    \qquad \beta = 1,2,3\,.
  \end{equation*}
\end{remark}

Summarizing, the Lagrangian evolution step consists in obtaining the
no-flow curves and use these directions to evolve the states:
\begin{equation}
  \begin{cases}
    \sigma_{\ell m}(\tempo^{n+1}) = \sigma_{\ell m}(\tempo^n) +
    \Dt\; \sigma_{\ell m}'(\tempo^n) 
    \quad\mbox{ s.t. } 
    \begin{cases}
      \scalprodSurfST{\normalEdge_{\ell m}(\tempo^n)}{\sigma_{\ell
          m}'(\tempo^n)}=0 \\
      \scalprodSurfST{\FluxLagrangian(\midPointEdge,\ConservVar{\midPointEdge}{n})}{\normalEdge(\tempo^n,\midPointEdge)}=\mathbf{0}\\
      \sigma_{\Edge[\ell m]}(\tempo^n) =
            [\tempo^n,\midPointEdge^1,\midPointEdge^2]^{\T}  
    \end{cases}\qquad\\[3em]
    \displaystyle \bConservVar{\ell}{n+1} = \frac{\CellArea[{\Region[n]_\ell}]}{\CellArea[{\bRegion[n+1]_\ell}]} 
    \left[\ConservVar{\ell}{n} +
      \frac{{\Dt}}{\CellArea[{\Region[n]_\ell}]}I({\bSource})
      - \Dt\int_{\mathcal{LS}^n} \scalprodSurfST{\FluxEulerian}{\normalEdge}
      \right]
    \\[2em]
    \ConservVar{\ell}{}(0,\midPoint_\ell)=\ConservVar{0}{}
  \end{cases}
  \label{DisLagEul}
\end{equation}
The numerical formulation is completed by the Eulerian remapping step,
in which the no-flow curves are followed back to remap the solution on
the computational grid.
The source term $\bSource$ can be integrated with any consistent
quadrature formula. Here we use either the mid-point or the
trapezoidal rules, whose derivation is provided in
\autoref{app:quadrule}.
Note that, it is possible to extend the discretization to other
quadrature rules~\cite{art:abreu19,art:abreu17,art:abreu21,art:abreu22b} to be able to reach high order in space.

\subsection{Stability and well-balance}

\paragraph{The CFL condition.}
Following~\cite{art:ABREU2023iswe,ABREU2022114011,ABREU2024115465}, the scheme is stable
under the following CFL condition:
\begin{equation}
  \label{eq:cflgeneral}
  \max_{\ell,m}\left\{
  \left|\frac{\FluxLagrangian^{^{\Depth,\normalEdge}}_{\Edge[\ell
        m]}}{\Depth}\right|,
  \left|\frac{\FluxLagrangian^{^{\Qcomp[1]{},\normalEdge}}_{\Edge[\ell
        m]}}{\Qcomp[1]{}}\right|,
  \left|\frac{\FluxLagrangian^{^{\Qcomp[2]{},\normalEdge}}_{\Edge[\ell
        m]}}{\Qcomp[2]{}}\right|
  \right\}
  \frac{\Dt}{\meshparam}< C\,,
\end{equation}
where $\FluxLagrangian^{\normalEdge}$ is the flux projection along the
normal direction to the edge $\Edge[\ell m]$, $\ell, m$ are the
(global) cell indices of neighboring cells, and $C<1/2$ is the CFL
constant.
We want to stress here that, due to the Lagrangian-Eulerian
formulation, it is not necessary to compute the flux Jacobian and this
is one of the main advantage of the approach. Moreover, in our case,
only the ``advective'' part $\FluxLagrangian$ of the flux function
enters into the CFL condition. Moreover, since the effects of the
additional local discontinuity are included in the source term, the
CFL condition is not affected by the extra terms.

\paragraph{Local well-balance and wetting-and-drying.}

The wet-dry transition may lead to large numerical errors when the
depth tends to be small and needs ad-hoc techniques to ensure
stability and conservation, in particular for non-planar topography.
Standard wetting-and-drying algorithms either perform
front-tracking of the wet-dry interface by mesh adaptation, or impose
a thin water layer in each cell e.g. by specifying that the water mass
does not move when the depth is below a specified threshold.
The latter approach may create oscillations in the case of non-flat
topography, due to topographic gradients erroneously propagating to
the solution~\cite{art:bunya2009wetting}, possibly yielding
negative depths.
These issues affect the solution of the Riemann problem in purely
Eulerian schemes, but is not detrimental to the accuracy of our
Lagrangian-Eulerian approach.
Indeed, in our approach, two important features contribute to
eliminate these problems altogether. On one hand, the geometrically
intrinsic nature of the proposed formulation automatically prevents
the possibility of producing negative water depths. In addition, its
combination with the Lagrangian-Eulerian steps eliminates potential
oscillations occurring when water depth, which appears at the
denominator of the flux function, tends to zero.
Indeed, the ratio $\Qdisch/\Depth$ concurs to the definition of the
no-flow surfaces and their Lagrangian time integration inherently
preserves non-negativity of the
depth~\cite{art:ABREU2023iswe,ABREU2022114011,ABREU2024115465}. Note
that, this ratio, which represents the water velocity, intervenes in
the definition of the local characteristic curves. 
For these same reasons, the proposed formulation is automatically
locally well-balanced.

Numerical verification of these properties can be done by looking at
the evolution in time of the ratio $\Qdisch/\Depth$, which should not
increase indefinitely, corresponding to the fact that $\Depth$ and
$\Qdisch$ go to zero with the same speed. Local well-balance is tested
by looking at the long-time behavior of the solution, in particular by
checking that the value of $I(\bSource)$ goes to zero.

\subsection{The numerical integration of the source term on a
  Cartesian grid}

Because of HPC efficiency and simplicity of implementation of our
Lagrangian-Eulerian intrinsic formulation on quadtree data
structures~\cite{gatti2023,GLF21,pochet2017new,tabarraei2005adaptive},
we decided to run our test cases on a quadrilateral mesh.  For this
reason we conclude this section by specializing the quadrature rule to
compute $I(\bSource)$ on a Cartesian grid.
As mentioned above, any quadrature rule and nonlinear reconstruction
techniques (see, e.g., ~\cite{DCET08,DET08,TACHTR2012,TACSH10,HAQ12})
can be employed to numerically compute the value $I(\bSource)$ in
\eqref{DisLagEul}.  In the following, we use standard first order
quadrature rules to account for the contribution of our
Lagrangian-Eulerian procedure. This procedure ensures the appropriate
balancing between the computation of the numerical flux and the source
term, in this naturally upwind Riemann-solver-free scheme.

Consider the region $D^n$ from time level $\tempo^{n}$ to
$\tempo^{n+1}$ such that $\tempo^{n+1}=\tempo^{n}+\Dt$, and enclosed
between the following no-flow surfaces (here $\beta$ identifies the
$\beta$-th equation in the system):
\begin{align*}
  \sigma_{i-\frac{1}{2},j}^{\beta}(\tempo)
  &\approx {x_{i-\frac{1}{2}}} + (\tempo -
  \tempo^n){\frac{\FluxLagrangian_{W}^{(\beta,2)}\metrcoefH{1}^2}{\ConservVar{}{\beta,n}}}\,,
  \quad
  \sigma_{i+\frac{1}{2},j}^{\beta}(\tempo) \approx x_{i+\frac{1}{2}} + (\tempo -
  \tempo^n){\frac{\FluxLagrangian_{E}^{(\beta,2)}\metrcoefH{1}^2}{\ConservVar{}{\beta,n}}}\,,\\ \sigma_{i,j-\frac{1}{2}}^{\beta}(\tempo)
  &\approx y_{j-\frac{1}{2}} + (\tempo -
  \tempo^n){\frac{\FluxLagrangian_{S}^{(\beta,3)}\metrcoefH{2}^2}{\ConservVar{}{\beta,n}}}\,,
  \quad
  \sigma_{i,j+\frac{1}{2}}^{\beta}(\tempo)
  \approx y_{j+\frac{1}{2}} + (\tempo -
  \tempo^n){\frac{\FluxLagrangian_{N}^{(\beta,3)}\metrcoefH{2}^2}{\ConservVar{}{\beta,n}}}\,.
\end{align*}
The final expression of $I(\bSource)$ after quadrature by midpoint
rule is:
\begin{multline}
  \label{eq:midpoint-final}
  I( {\bSource} )
  \approx \left\{
  \Dt\,\meshparam^2 \,
  \bSource\left(x_j^n+\left[\frac{\FluxLagrangian_{W}^{(\beta,2)}\metrcoefH{1}^2}{\ConservVar{}{\beta,n}}
    -
    \frac{\FluxLagrangian_{E}^{(\beta,2)}\metrcoefH{1}^2}{\ConservVar{}{\beta,n}}\right]\frac{\Dt}{2} + \frac{\meshparam}{2},\right.\right.\\
  \left.
  \left.
  y_j^n+\left[\frac{\FluxLagrangian_{N}^{(\beta,3)}\metrcoefH{2}^2}{\ConservVar{}{\beta,n}} - \frac{\FluxLagrangian_{S}^{(\beta,3)}\metrcoefH{2}^2}{\ConservVar{}{\beta,n}}\right]\frac{\Dt}{2} + \frac{\meshparam}{2},
  \tempo^n + \frac{\Dt}{2}\right)\right\}_{\beta=1,2,3}\,.
\end{multline}
See \autoref{app:quadrule} for the detailed computations and the
extension to the trapezoidal rule.

\section{Numerical simulations}
\label{sec:results}

The simulations presented in this section are devised in order to
validate the main properties of our proposed approach.
We specifically focus on the ability of the method to circumvent the
numerical difficulties associated with:
\begin{enumerate*}[label={\roman*)}]
  \item  discrete well-balance;
  \item wetting-and-drying transitions and moving fronts;
  \item  preservation of hyperbolicity and stability limits in long-time
    simulations; and,
  \item handling local topographic discontinuities.
\end{enumerate*}
We highlight four key features distinguishing our framework from
standard methods:
\begin{description}
\item[Capture of Non-Trivial Equilibria:] the scheme naturally
  recovers non-trivial equilibrium solutions, essential for
  hyperbolic balance laws where the source term and flux
  gradient must cancel exactly.

\item[Flux-Source Interaction:] the scheme accurately captures the
  interplay between convective fluxes and source terms, preventing
  spurious oscillations when they counteract each other.

\item[Long-Term Asymptotic Stability:] the method maintains
  equilibrium states as time tends to infinity, avoiding cumulative
  drift.

\item[Implicit Well-Balance:] the method achieves well-balance
  without external root-finding procedures. This property emerges
  naturally for the Lagrangian evolution and the geometric
  remapping, extending seamlessly to the intrinsic geometric
  setting.
\end{description}
These properties are crucial for real-world
applications, particularly in ``lake-at-rest'' scenarios with
vanishing flow velocity. Our numerical experiments demonstrate that
the scheme achieves excellent resolution for both unperturbed and
slightly perturbed equilibria, naturally preserving steady-states at
rest and complex wet-dry transitions. 

All the test cases consider a gravity-driven flow, assuming an
inviscid flow, i.e. no energy losses, $\Tens[sw]=0$, and
$\BottomFriction=0$.  In each simulation the time step size is defined
by the CFL stability condition in \autoref{eq:cflgeneral}.
Obviously, due to the absence of significant tests with exact
solutions, convergence can be verified only qualitatively, or against
a highly refined solution.
The different simulations that follow are represented by means of
three-dimensional (axonometry) renderings of the normal depth $\Depth$
at different time instances. A single viewport is used to
better appreciate the dynamics of the important features of the
solution thus allowing the evaluation of its properties.
The specific times and the viewports are chosen to highlight the
qualitative behavior of the symmetries in the obtained solution.
We would like to stress here that in some test cases the choice of a 
single viewport may hide specific features of the solution at certain
times, in which case symmetries may not be fully appreciated.
Also, we include normal depth contour lines with constant jumps
$\Delta\eta$ between successive contours. The values of $\Delta\eta$
changes at each time for a better visual rendering of symmetries of
the normal depth distribution and to identify even the tiniest
waves appearing in the numerical solution.
In addition, for specific purely gravity-driven test cases, we show
the time evolution of the $\Lspace[1]$-norm of the forcing
$I(\bSource)$ (see eq.~\eqref{eq4}) and of the $\Lspace[1]$- and
$\Lspace[\infty]$-norm of the fluid depth and velocity
components. These plots are used to indicate achievement and
time-persistence of steady state conditions.
\subsection{Well-balance and literature benchmarks}

A critical validation for well-balanced schemes is equilibrium
preservation over extended temporal scales. We consider three
literature benchmarks demonstrating the scheme ability to reach and
maintain equilibrium as $\tempo \to \infty$, from both steady-state and
slightly perturbed initial conditions. In these test cases, the
effects of the surface topography is not considered
(i.e. $\SurfDomain$ is a flat region).

\paragraph{Test case 1.1: Long-time equilibrium test.}
The simulation is conducted on a closed (no-flow boundary conditions)
domain $\SubsetU=[-3,3] \times [-3,3]\subset\REAL[2]$ with constant
initial water depth in static equilibrium, where the pressure gradient
and gravitational source terms must cancel exactly.  The mesh is
formed by $256\times 256$ cells and the simulation uses CFL=0.1.
Our method successfully maintains this equilibrium state over
long-time simulations, a task that is traditionally challenging due to
the accumulation of truncation errors. This is fundamentally due to
the \emph{no-flow curves} in our Lagrangian-Eulerian formulation.
In fact, the no-flow curve ensures that the union of local regions,
$\bigcup_{\ell=1}^{n} \Region_{\ell}$, forms a conservative line of
flow, allowing the discrete evolution to match the physical stable
state.
Mathematically,  the balance 
\begin{equation*}
  \DivSurf \FluxSWE(\ConservVar{}{}) = \Source(\ConservVar{}{})
\end{equation*}
is satisfied at the discrete level within each local region
$\Region_{\ell}$, preventing the generation of spurious waves.

\begin{figure}
  \centering
  \includegraphics[width=\textwidth]{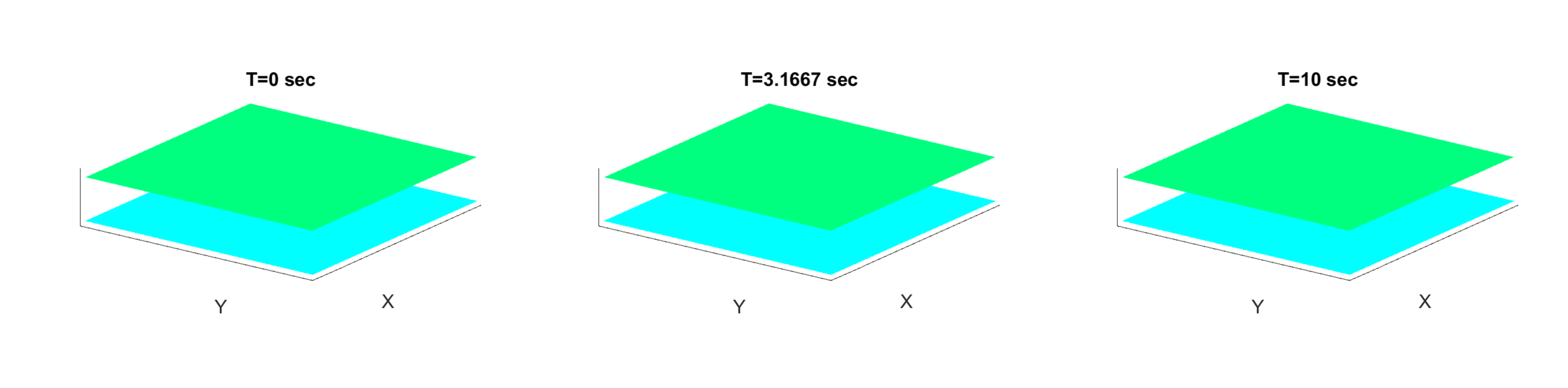}
  \caption{Test case 1.1 - Long-term preservation of the quiescent
    steady-state. The three panels show the evolution at different
    times: $T=0$~s (left), $T=3.1667$~s (central), and $T=10$~s
    (right). The static equilibrium is maintained without the
    emergence of spurious oscillations or numerical instabilities.}
  \label{WB1}
\end{figure}

\autoref{WB1} shows the numerical results: the steady-state
equilibrium is maintained along the simulation, reaching a long-time
state ($T=10$~s) without the emergence of spurious oscillations or
numerical instabilities.

\paragraph{Test case 1.2: Small perturbations to the ``lake-at-rest''.}
This benchmark, discussed in~\cite{hernandez2011shallow}, considers
the propagation of small-amplitude waves over a stagnant background,
and requires the scheme to resolve subtle fluctuations without
disrupting the underlying hydrostatic equilibrium.
The initial state consists in $\Depth=1.0$~m in the domain $\SubsetU=
[-4,4]\times[-10,10]$, with a localized small-perturbation of this
steady-state within the subset $[-4,4]\times[-0.4,0.4]$, where
$\Depth=1.5$~m. 
Outflow boundary conditions are imposed at the left and right
boundaries. In this case, the mesh is formed by $256\times 640$ cells
and the simulation consider CFL=0.15.

\begin{figure}
  \centering
  \includegraphics[width=\textwidth]{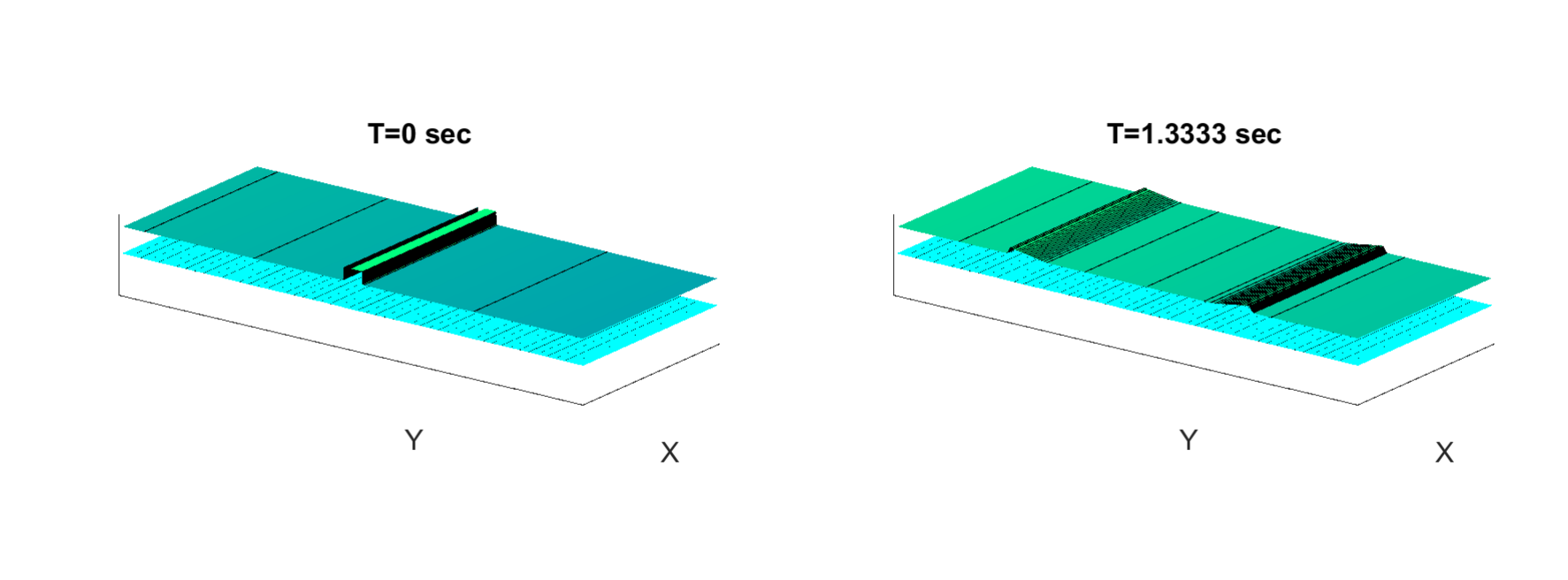}\\[-2.5em]
  \includegraphics[width=\textwidth]{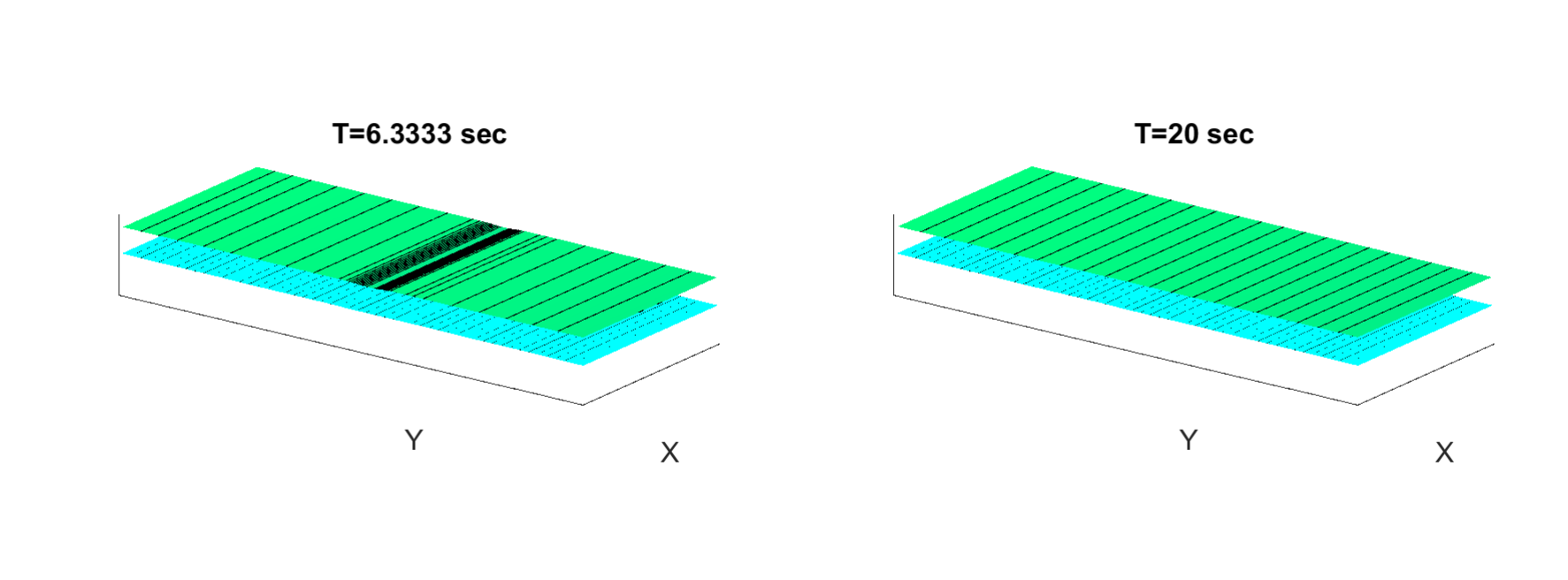}\\[-3em]
  \caption{Test case 1.2 - Evolution of an initial perturbation toward
    a quiescent steady-state. The initial small-amplitude perturbation
    superimposed on a fluid at equilibrium (top-left panel,
    $T=0$~s). The wave bifurcates into two counter-propagating gravity
    waves that traverse the domain (top-right panel) and later
    reconnect (bottom-left panel), before completely exiting the
    domain (bottom-right panel).}
  \label{WB2}
\end{figure}  

As shown in \autoref{WB2}, the perturbation splits into two
counter-propagating waves that travel across the domain.  The sequence
demonstrates the framework ability to maintain a well-balanced steady
state without triggering spurious numerical oscillations.
\begin{remark}\label{rem:dissipation}
  Once the initial waves departing from the center reach the outflow
  boundaries small reflected waves are formed and move back towards
  the center of the domain.  These waves dissipate quickly and are not
  propagated further in the domain.  Indeed, the appearance of the
  waves is visible only because of the concentration of the contour
  lines, and their depth is always negligible. We would like to note
  here that the dissipation properties of both the ISWE model and the
  Lagrangian-Eulerian numerical approach are still open questions that
  are currently under investigation. Hence, here we have to limit
  ourselves with the observation of the obtained numerical results.
\end{remark}
By $T=20$~s (bottom-right panel), the transient waves have exited the
domain, and the system has successfully recovered the steady-state of
rest, confirming the long-term stability and geometric consistency of
the intrinsic Lagrangian-Eulerian framework.

\paragraph{Test case 1.3: evolution towards a wet/dry steady-state.}

This benchmark, adapted from the one-dimensional test proposed
in~\cite{hernandez2011shallow}, evaluates transitions from dynamic
flow to non-trivial equilibrium involving dry areas. The test
specifically challenges the scheme handling of the
\emph{wetting-and-drying} interface and the \emph{well-balance} in the
presence of topographic obstacles.

The channel domain is defined by $\SubsetU = [-4,4] \times [-10,10]$,
with a topographic obstacle modeled as a local perturbation,
$\obstacle(x) = \frac{1}{2} \cos(\pi x)$ within the subset
$ [-4,4] \times [-0.4, 0.4]$. Initially, the fluid is at rest ($u=0$)
with a leveled total water height of $H=0.8$~m. The water is allowed
to drain through the right boundary (outflow), while the remaining
boundaries are treated with solid wall conditions. The mesh is formed
by $256\times 640$ cells and CFL=0.15.

\begin{figure}
  \centering
  \includegraphics[width=\textwidth]{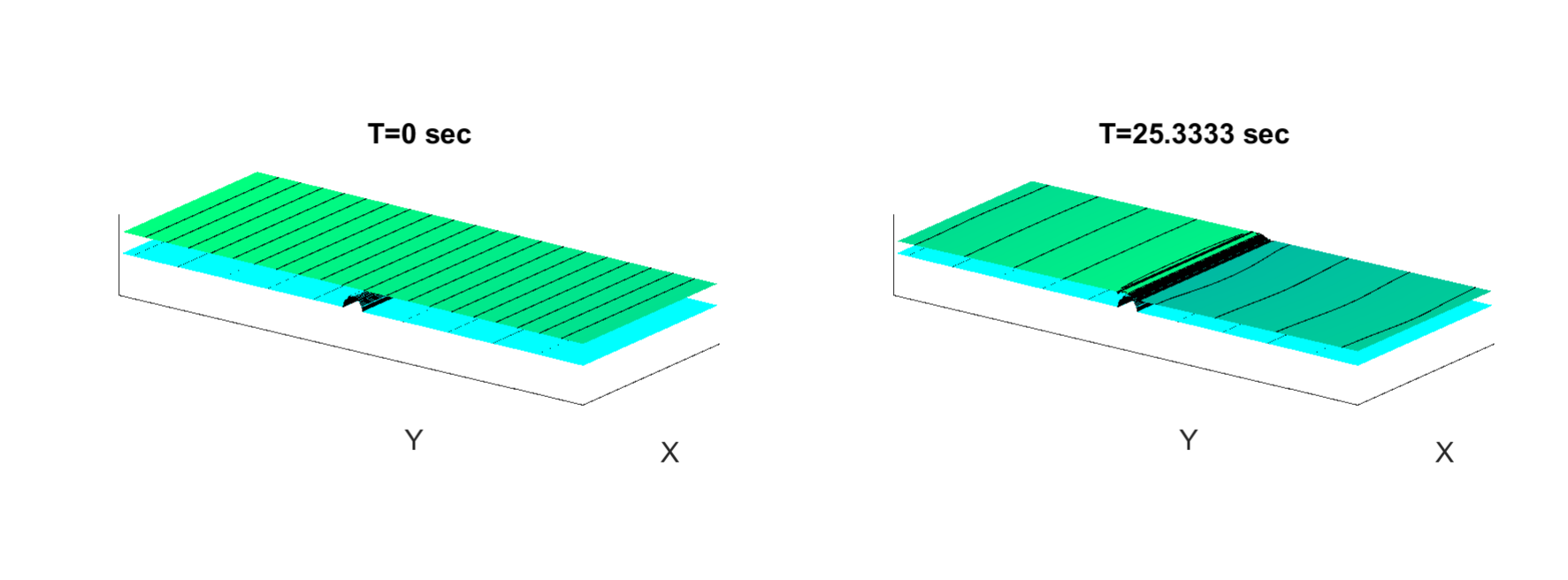}\\[-2.5em]
  \includegraphics[width=\textwidth]{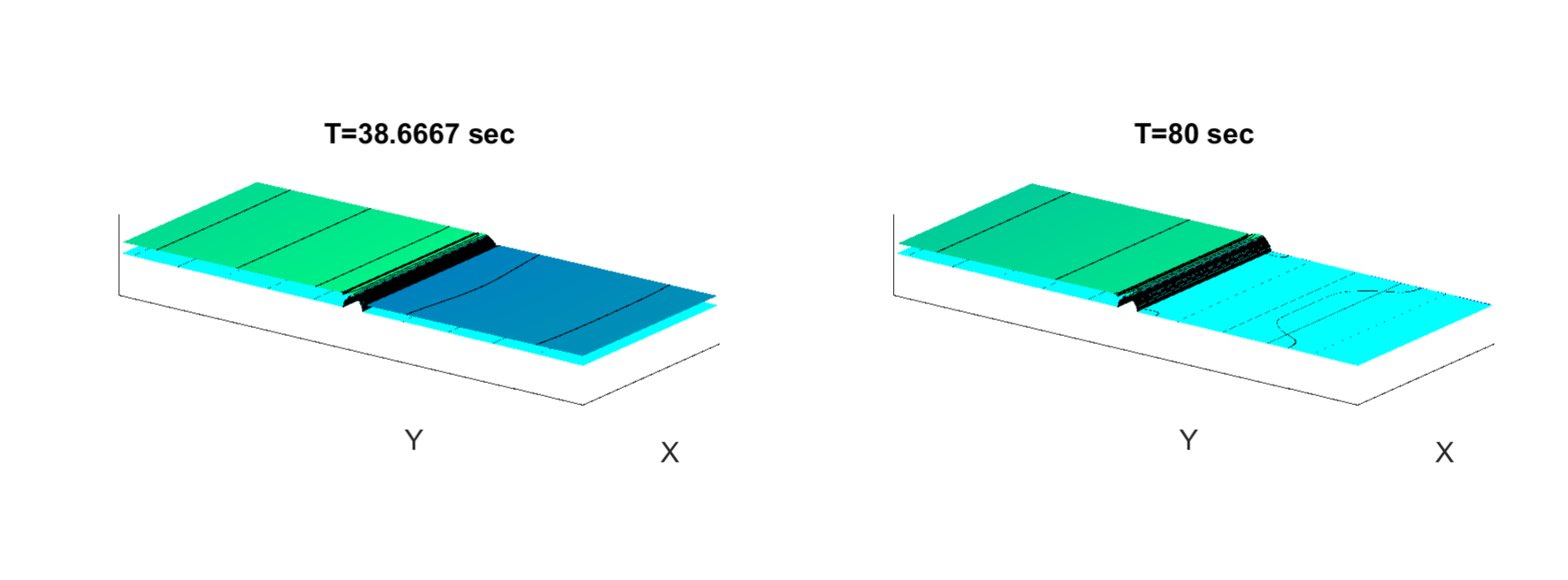}\\[-3em]
  \caption{Test case 1.3 - reservoir drainage. Top-left panel: initial
    data ($T = 0$~s). Top-right ($T = 25.3333$~s) and bottom-left ($T
    = 38.6667$~s) panels: the water starts exiting the domain from the
    outflow boundary, converging to a hybrid wet/dry steady-state in
    the two different regions. Bottom-right panel ($T = 80$~s): final
    equilibrium configuration.}
  \label{WB33}
\end{figure}

\autoref{WB33} shows the simulation progress: the fluid drains rightwards
and the obstacle traps the water on the left portion of the domain.
We reach the hybrid steady-state: a wet-equilibrium state in the left
portion, and a dry-equilibrium state in the right portion of the
domain where the outlet is located.
\autoref{WB33}, bottom-right panel, shows the right portion of the
reservoir essentially fully drained at $T=80$~s. The
Lagrangian-Eulerian method maintains \emph{hyperbolicity} as the water
depth vanishes, preventing the appearance of numerical instabilities
or non-physical velocities, often encountered in the presence of
moving fronts.

\subsection{Local obstacles over general bottom topographies}

Here, we consider four different bottom surfaces with an increasing
complexity in the geometric features (see~\cite{art:ABREU2023iswe}):
\begin{enumerate*}[label={\roman*)}]
\item a simple sloping plane (i.e., constant metric);
\item a parabola
  (i.e., one-dimensional curvature effects);
\item a hyperboloid with a
central bump (i.e., two-dimensional dependence of the geometry with a
central radial symmetry); and,
\item a fully three-dimensional
  surface.
\end{enumerate*}
In all the cases, a localized discontinuous obstacle is
added.
All the surfaces are defined by the use of a (single) global
parametrization $\SurfDomain = \{\xcg,\ycg,\zcg=\BSM(\xcg,\ycg)\}$,
where $\BSM(\xcg,\ycg)$ is the height function. The localized
discontinuous obstacle is defined by the equation
\begin{equation}\label{eq:bump}
  \obstacle(\xcg,\ycg) = \bar{z} - (\xcg-a)^2 - (\ycg-b)^2\,,
\end{equation}
where $\bar{z}$ is the maximum height of the obstacle attained at
point $(a,b)$.  In all test cases, the support of the obstacle is the
square $|\xcg-a|<0.8$ and $|\ycg-b|<0.8$.
This kind of bump is representative of localized obstacles on top of a
regular bottom. The fact that there is only one obstacle is not a
limitation of the proposed approach and a test case considering
multiple bumps is presented in the next subsection.
The working grid discretizing the bottom surface is obtained from a
regular Cartesian grid of a subset $\SubsetU\subset\REAL^2$, i.e. the
chart domain, and projecting vertically to the surface the vertices by
the height function.
The initial conditions describe dam-break problems, with prescribed
normal water depth with zero velocity.

All the tests are developed to highlight that the well-balance
property is preserved at the discrete level also in the presence of
non-trivial bed topography and local discontinuities.
In these cases, the verification of the absence of spurious
oscillations is performed by looking at the long-time behavior of the
$\Lspace[\infty]$- and $\Lspace[1]$-norms of $I({\bSource})$.

\paragraph{Test case 2.1: sloping plane with obstacle.} 
The chart domain is the rectangle
$\SubsetU=[-4,4]\times[-10,10]\subset\REAL^2$ and the height function
is given by $\BSM(x,y)=0.1 y+1.0$. The discontinuous obstacle is
defined by $\obstacle$ in \autoref{eq:bump} with $\bar{z}=1.2$, $a=0$,
and $b=4$.
The initial value for the water depth is set to 2~m for cells upstream
of the dam position, located at $y=2.0$~m, and a depth of 1~m
downstream. A CFL value of $0.20$ is considered to run the simulation.

\begin{figure}
  \centering
  \includegraphics[width=\textwidth]{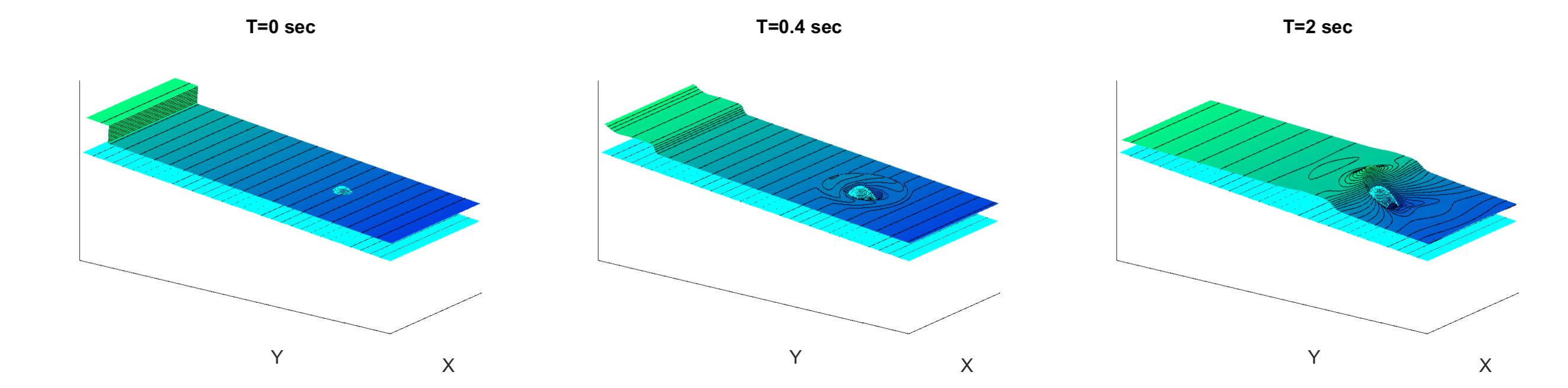}\\
  \includegraphics[width=\textwidth]{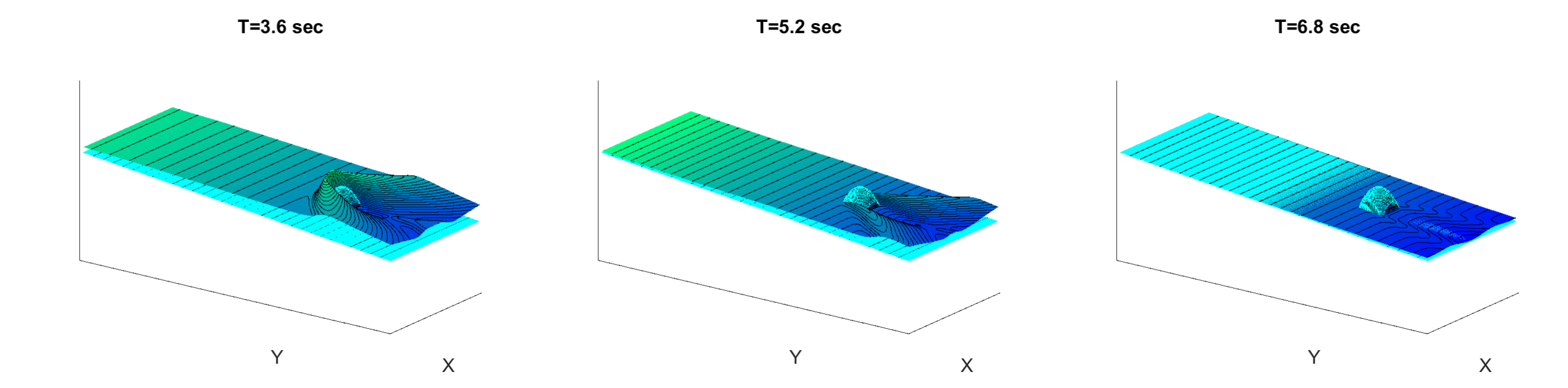}\\
  \caption{Test case 2.1 - sloping plane with obstacle: evolution of
    the normal water depth (drawn above the bottom surface) at
    different times. The interaction of the waves with the obstacle
    are clearly visible. The results showcase the robustness of the
    numerical scheme for long-time simulations also when the water
    depth tends to zero.}
  \label{fig:Slopeplane1}
\end{figure}

\begin{figure}
  \centering 
  \includegraphics[width=0.3\textwidth]{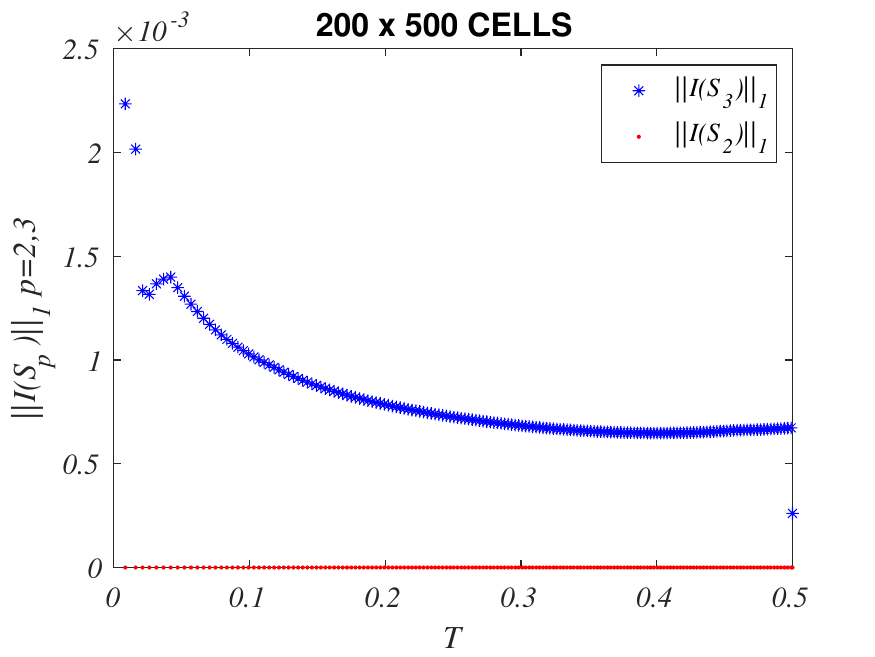}
  \includegraphics[width=0.3\textwidth]{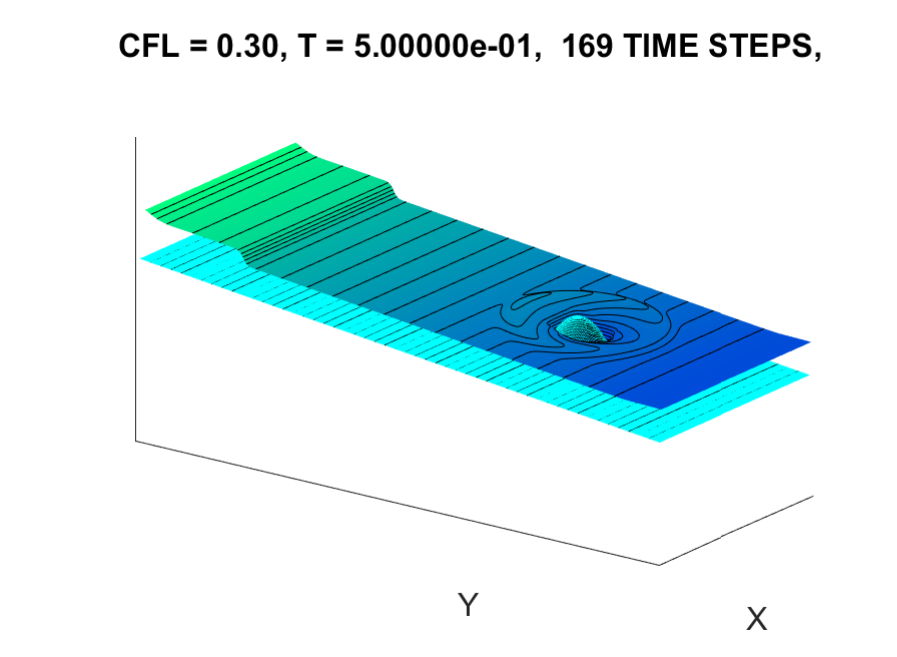}
  \includegraphics[width=0.3\textwidth]{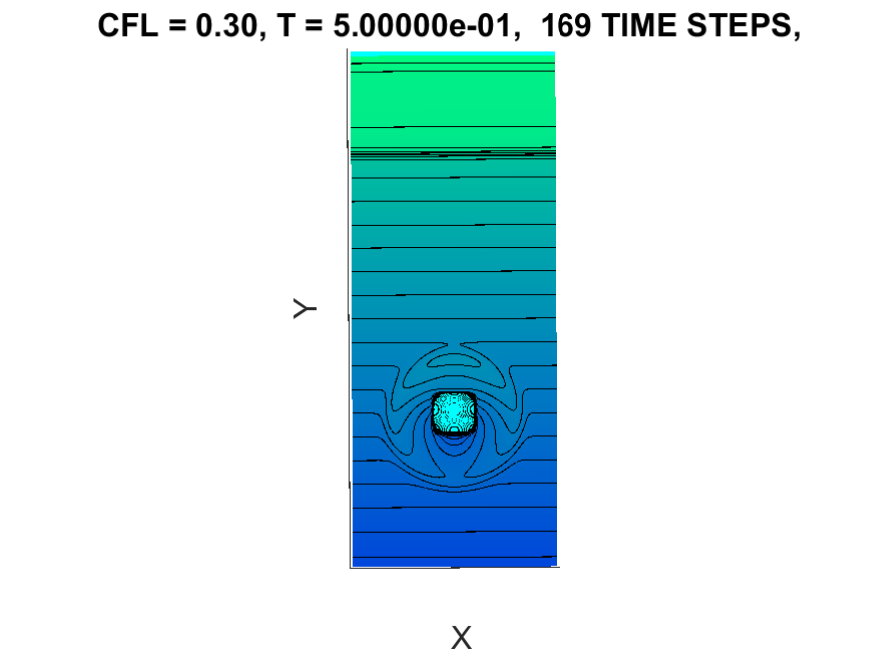}

  \includegraphics[width=0.3\textwidth]{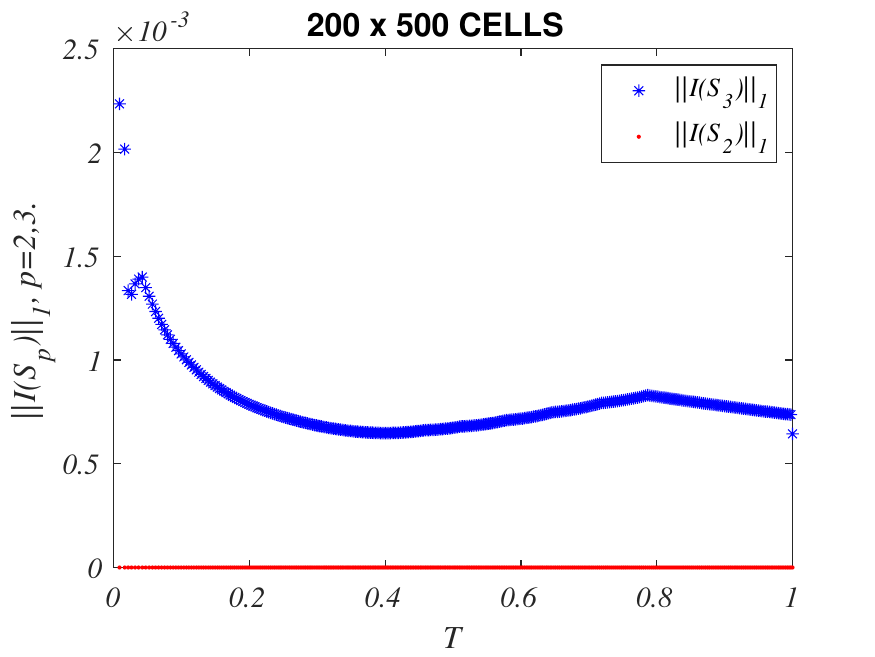}
  \includegraphics[width=0.3\textwidth]{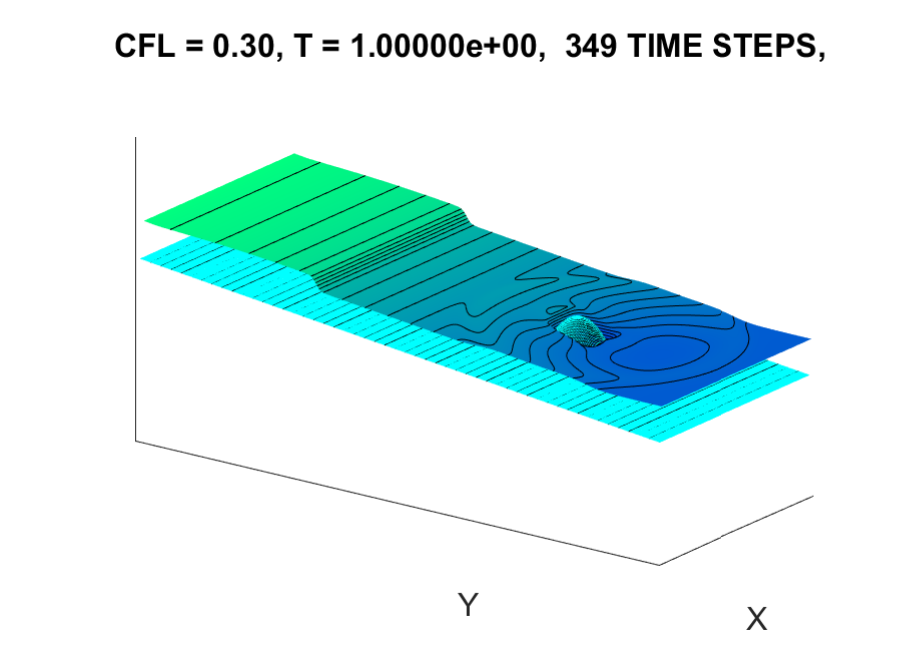}
  \includegraphics[width=0.3\textwidth]{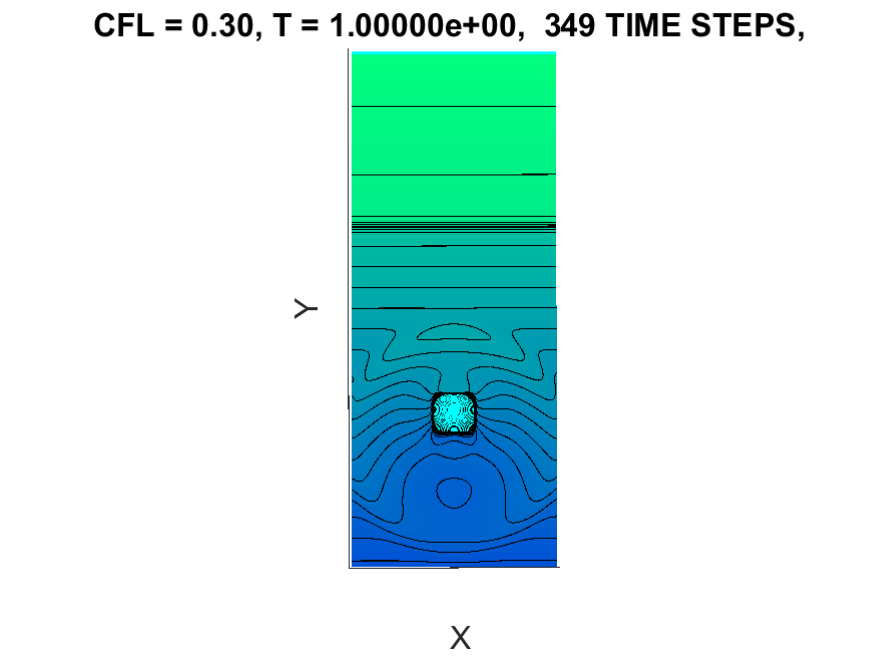}

  \includegraphics[width=0.3\textwidth]{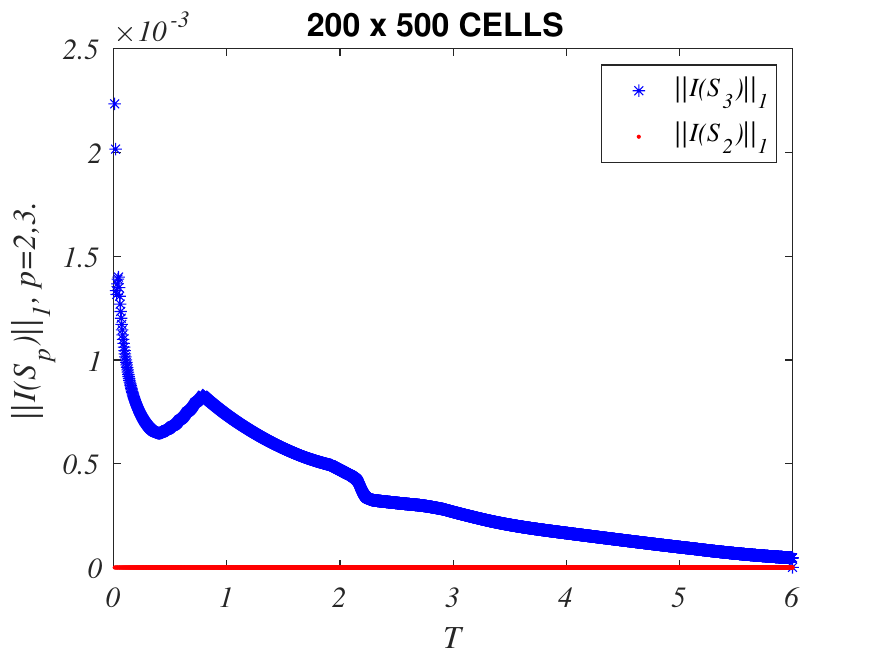}
  \includegraphics[width=0.3\textwidth]{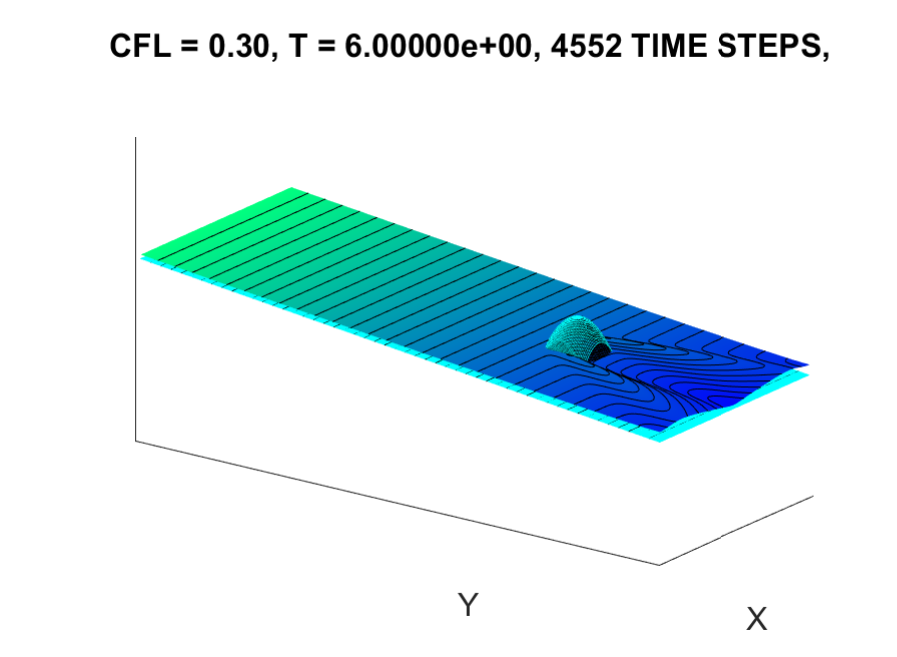}
  \includegraphics[width=0.3\textwidth]{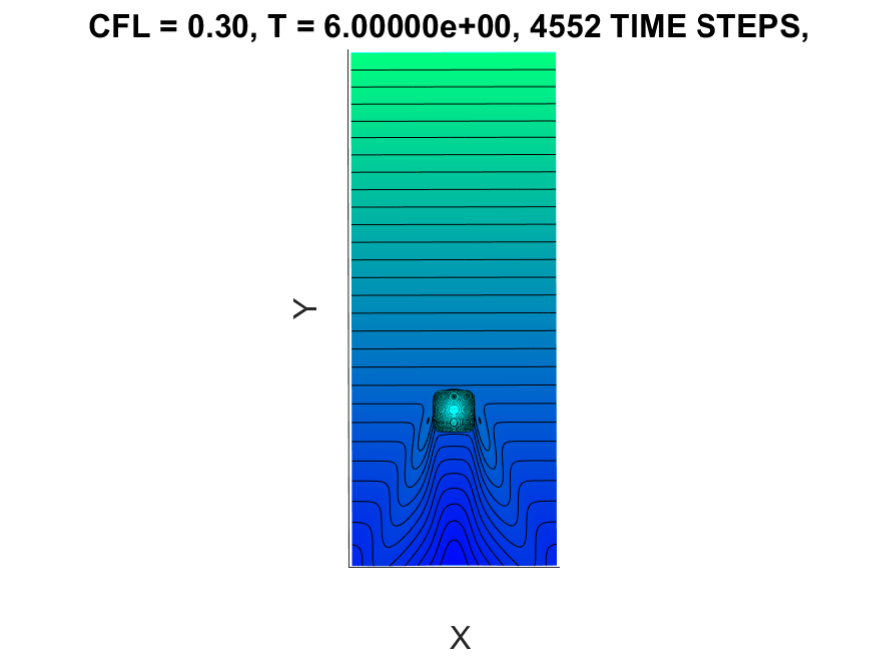}
  \caption{Test case 2.1 - sloping plane with obstacle: well-balance
    simulation. The left column shows the time behavior of the $L_1$
    norms of $I(\bSource)$ in the interval $[0,T]$ for
    $T=0.5, 1.0, 6.0$~s. The central and right columns show the
    solution in terms of water depth in the same time intervals. The
    contours in the right column emphasize the symmetry of the
    numerical solution with respect to the flow direction.}
  \label{fig:Slopeplane2}
\end{figure}

\begin{figure}
  \centering
  \includegraphics[width=0.3\textwidth]{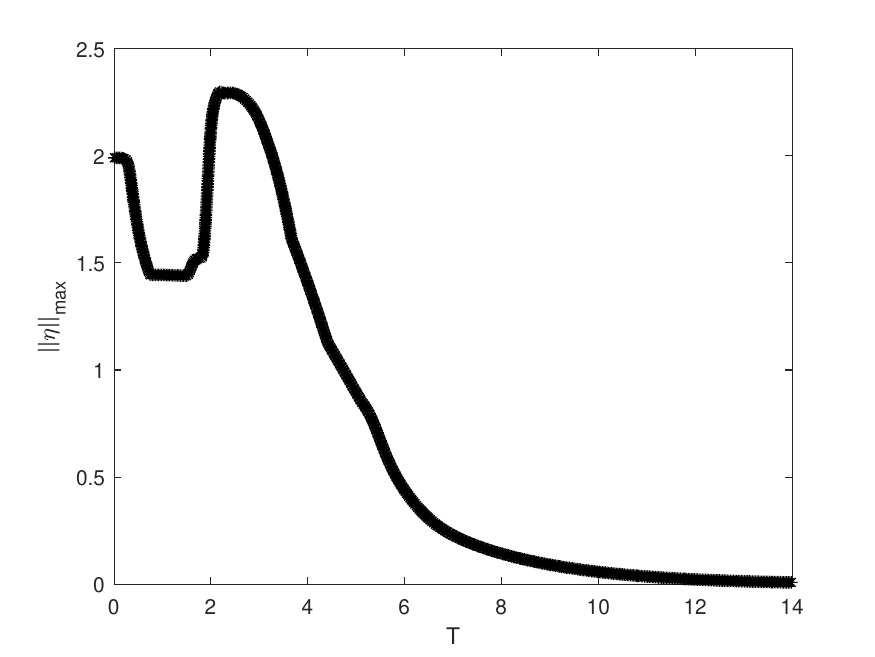}
  \includegraphics[width=0.3\textwidth]{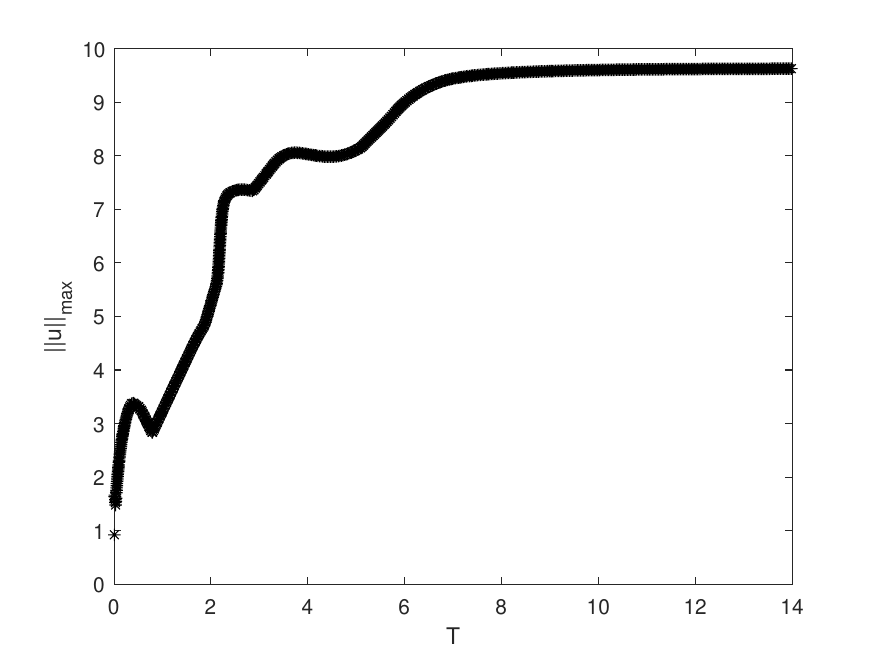}
  \includegraphics[width=0.3\textwidth]{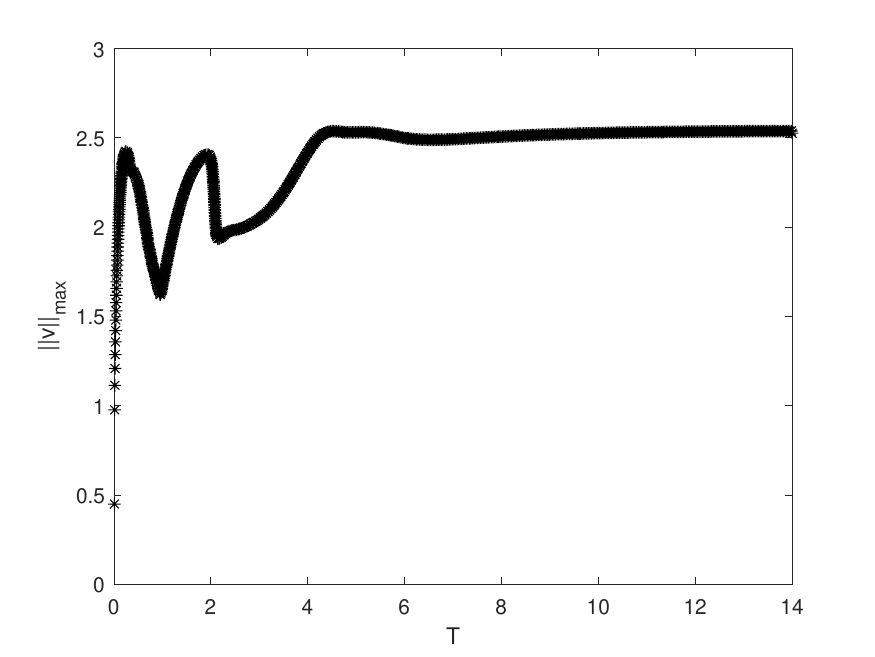}
  
  \includegraphics[width=0.3\textwidth]{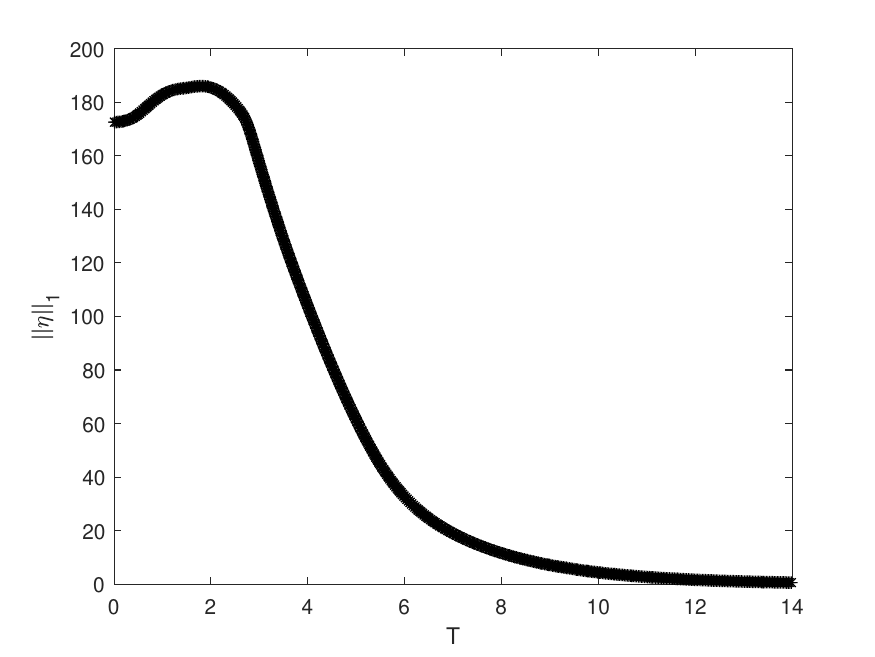}
  \includegraphics[width=0.3\textwidth]{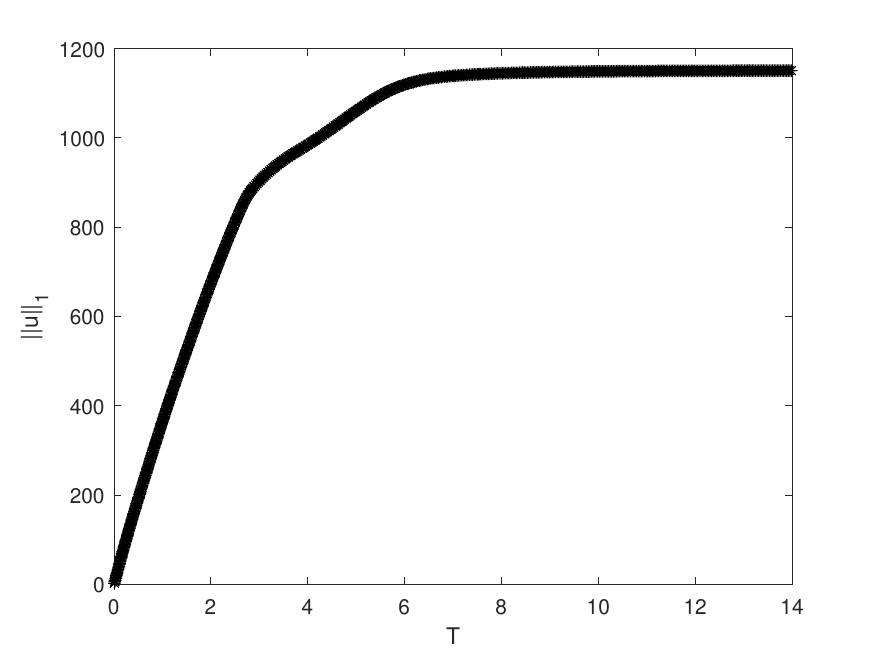}
  \includegraphics[width=0.3\textwidth]{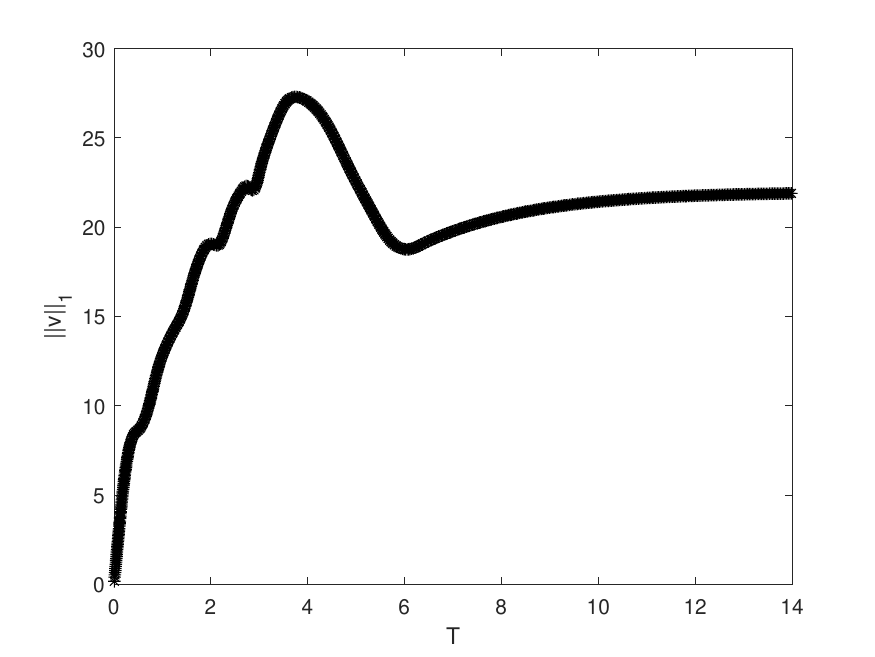}
  \caption{Test case 2.1 - sloping plane with obstacle:
    wetting-and-drying simulations. Values of $\Lspace[\infty]$- and
    $\Lspace[1]$-norms are shown respectively in the top and bottom
    rows for normal water depth $\Depth$ (left column) and the $x$ and
    $y$ components $u$ and $v$ of the tangent velocities (central and
    right columns).  }
  \label{Hplane}
\end{figure}

\autoref{fig:Slopeplane1} shows the evolution of the normal depth at
different times in the interval $T\in [0,8]$~s, on a grid formed by
$256\times 640$ cells along the $\xcl$ and $\ycl$-directions,
respectively.
It is possible to observe that the wave front is captured by the
numerical solution as time evolves.  The method describes well the
effects of the obstacle on the fluid in both directions. The
discretization is robust also in the transition between wet and dry
areas, as shown in the panels in the bottom row, where we reach an
almost dry state for the entire domain.
At $T=0.4$~s the solution displays some dissipation evidenced by
the somewhat smeared shock profile moving downstream. Again, as
stated in \autoref{rem:dissipation}, the amount of numerical
dissipation added by the scheme has not been fully tested and is
currently under investigation.

\autoref{fig:Slopeplane2} shows the evolution of the norm of the forcing
$I(\bSource)$ and the water depth distribution at $T=0.5, 1.0$~s and
$T=6$~s.  In particular, this figure is reported to underline that the
well-balance property of the system is preserved by the numerical
scheme, together with a good symmetrical behavior of the solution
around the obstacle.
\autoref{Hplane} shows both $\Lspace[\infty]$ (top row) and
$\Lspace[1]$ (bottom row) norms of $\Depth$ (left panels) and
$\Qdisch/\Depth=[u,v]$ (central and right panels) for
a long-time simulation of $T=14$~s, which in this case represents the
steady-state solution.
We note that as soon as the value of the depth tends to zero, the
velocities reach their limit values, showing, as mentioned at the end
of the previous section, that our formulation is well-balanced also in
the discrete setting.

\paragraph{Test case 2.2: parabola with obstacle.}
This is a simple one-dimensional bottom surface in which the effect of
the curvature is present. The parabola is defined by the height
function $\BSM(x,y)=\frac{1}{25}(y-10)^2$ on the rectangular domain
$\SubsetU=[-4,4]\times[-10,10]\subset\REAL^2$, and here the
obstacle is defined by $\obstacle$ in \autoref{eq:bump} with
$\bar{z} = 1.2$, $a=0$, and $b=4$.
The domain is discretized by a grid formed by $256\times640$ cells.
The initial depth is 2~m of water upstream and 1~m downstream of the
dam located at $y=2.0$~m and a CFL value of $0.30$ is considered.

\begin{figure}
  \centering
  \includegraphics[width=\textwidth]{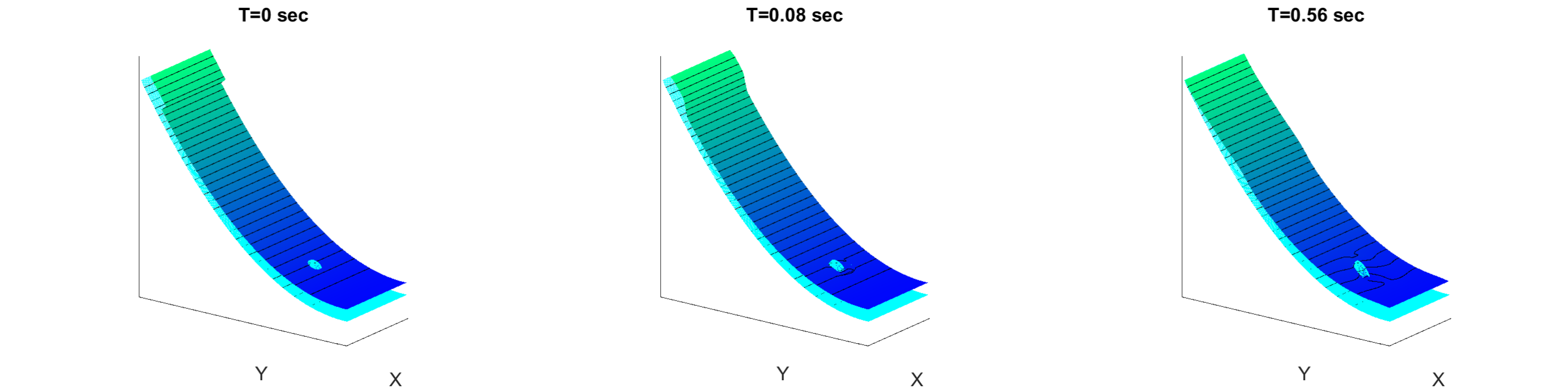}\\
  \includegraphics[width=\textwidth]{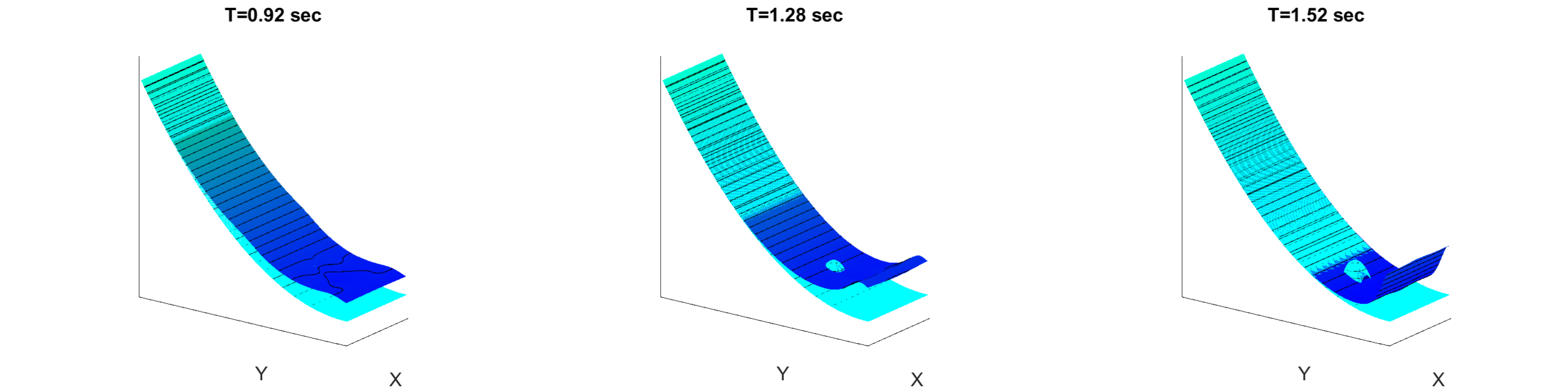}\\
  \caption{Test case 2.2 - parabola with obstacle: evolution of the
    water depth at different times. The grid is formed by 256$\times$640
    cells.}
  \label{fig:parabola1}
\end{figure}

\begin{figure}
  \centering
  \includegraphics[width=0.3\textwidth]{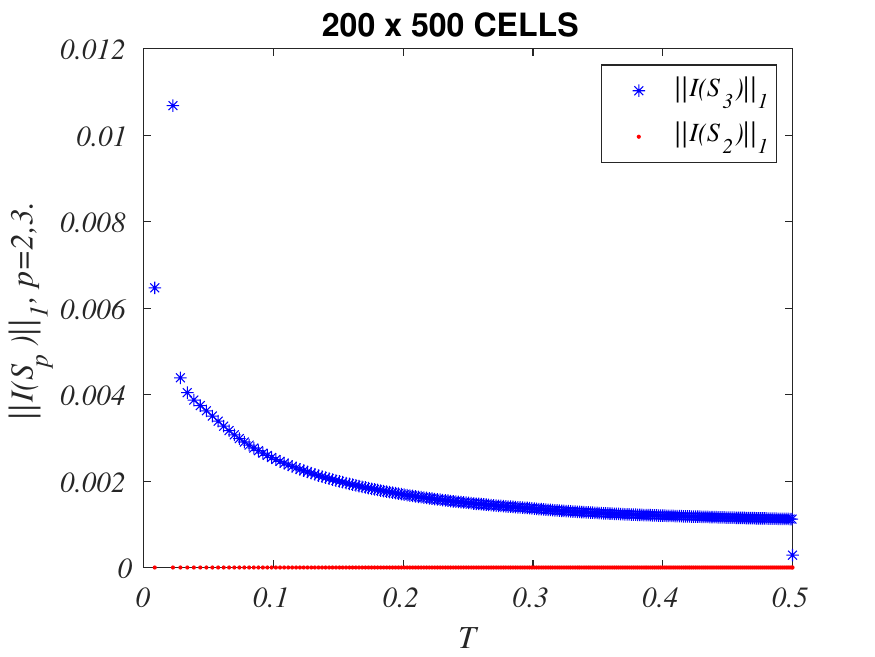}
  \includegraphics[width=0.3\textwidth]{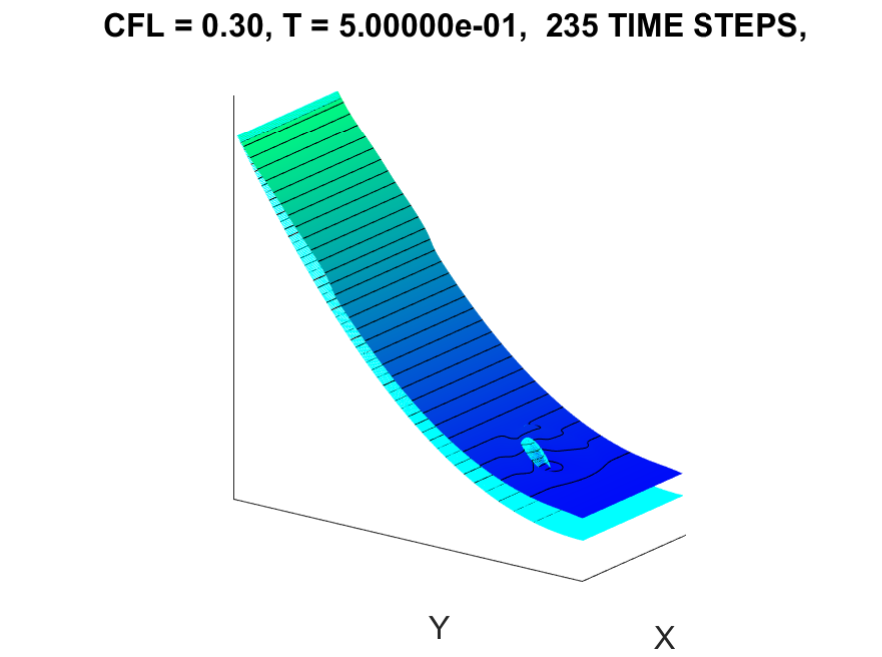}
  \includegraphics[width=0.3\textwidth]{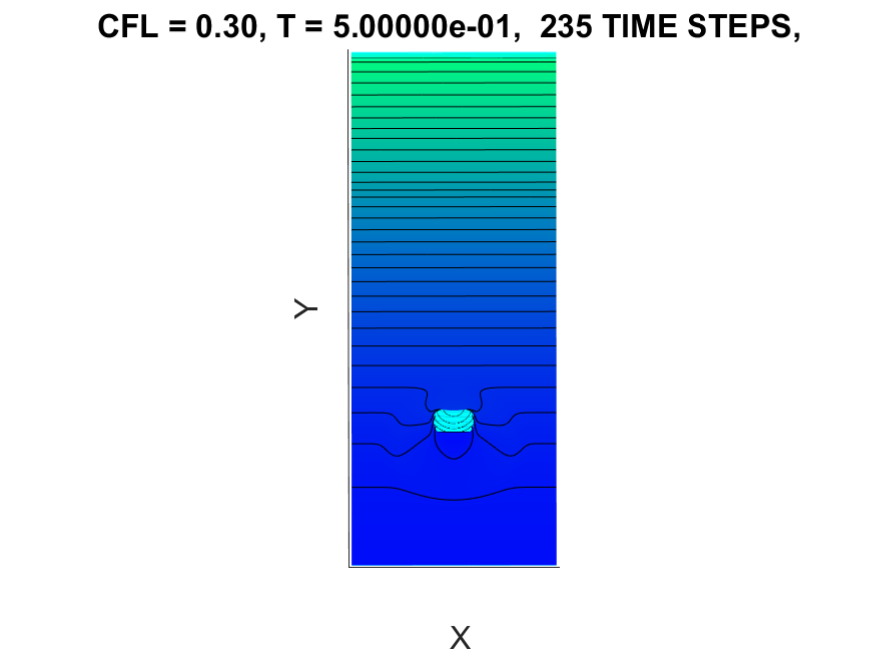}
  
  \includegraphics[width=0.3\textwidth]{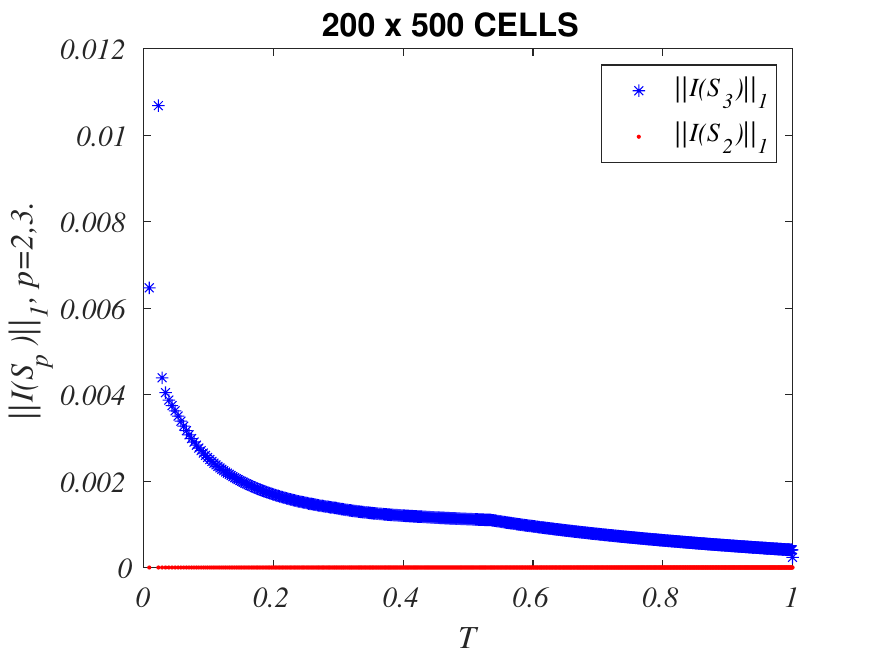}
  \includegraphics[width=0.3\textwidth]{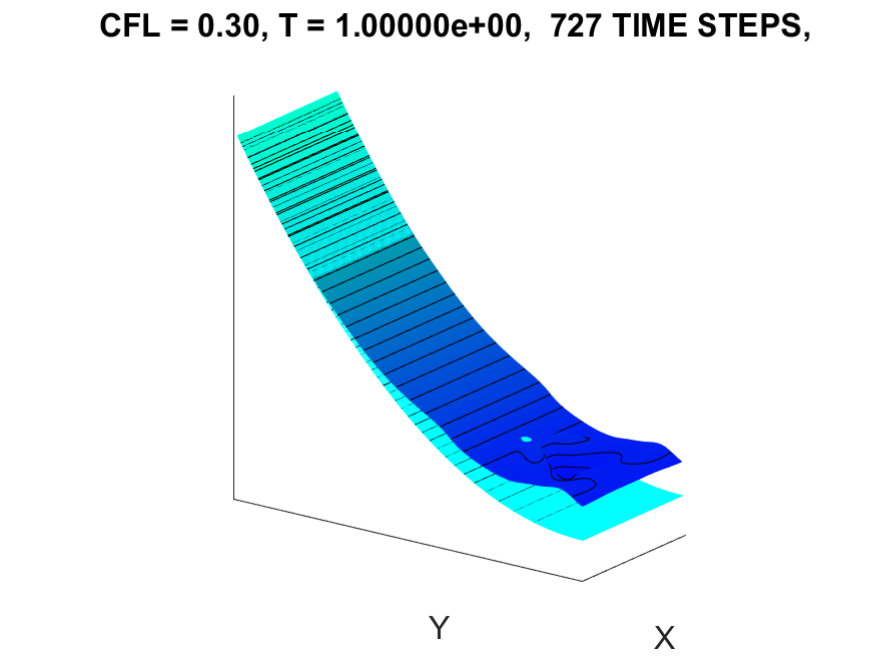}
  \includegraphics[width=0.3\textwidth]{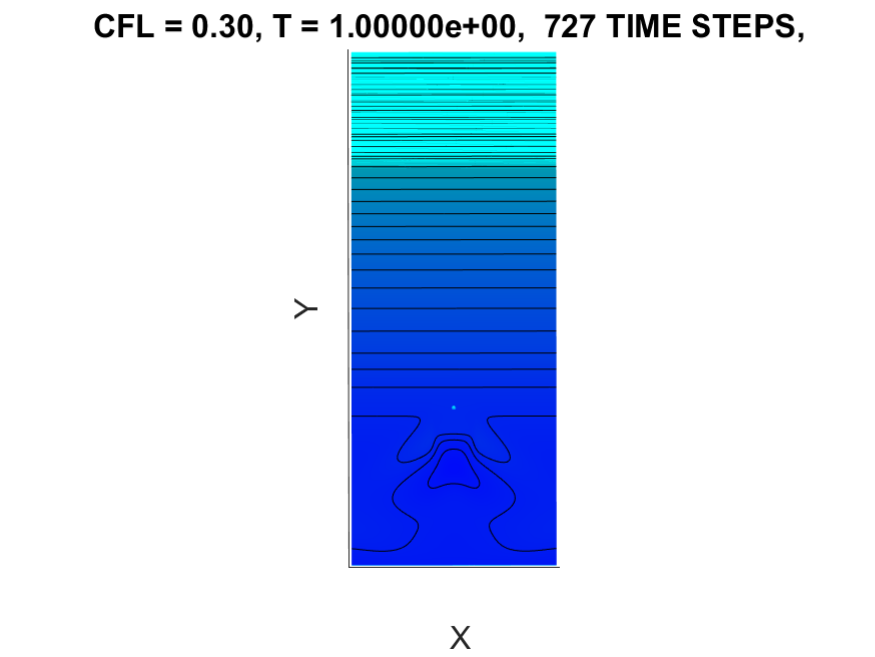}
  
  \includegraphics[width=0.3\textwidth]{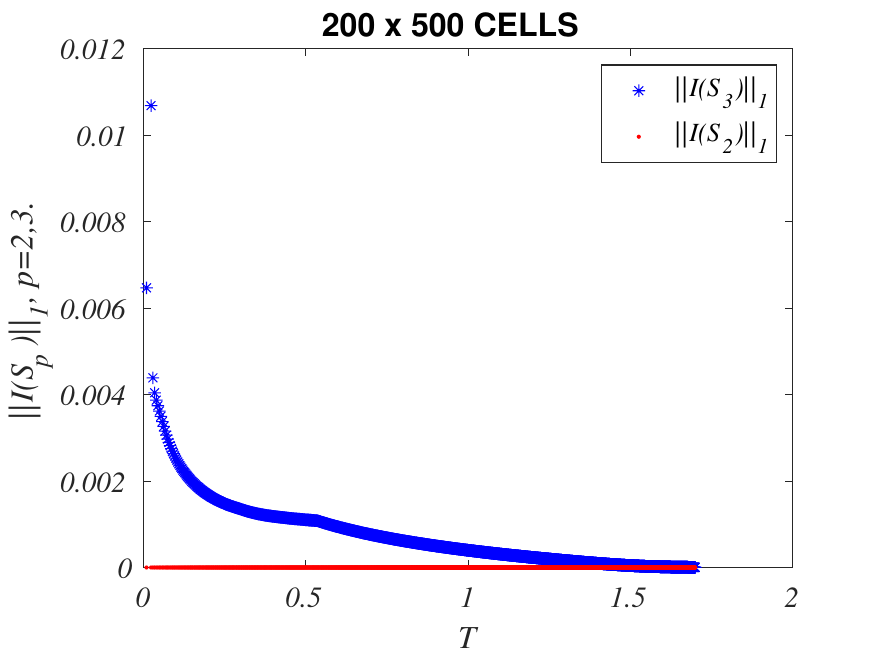}
  \includegraphics[width=0.3\textwidth]{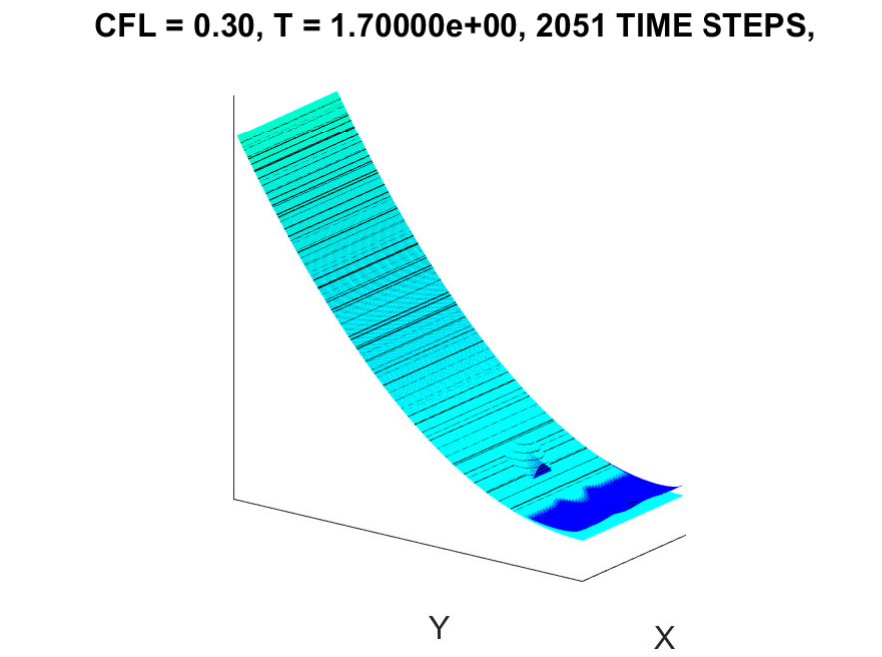}
  \includegraphics[width=0.3\textwidth]{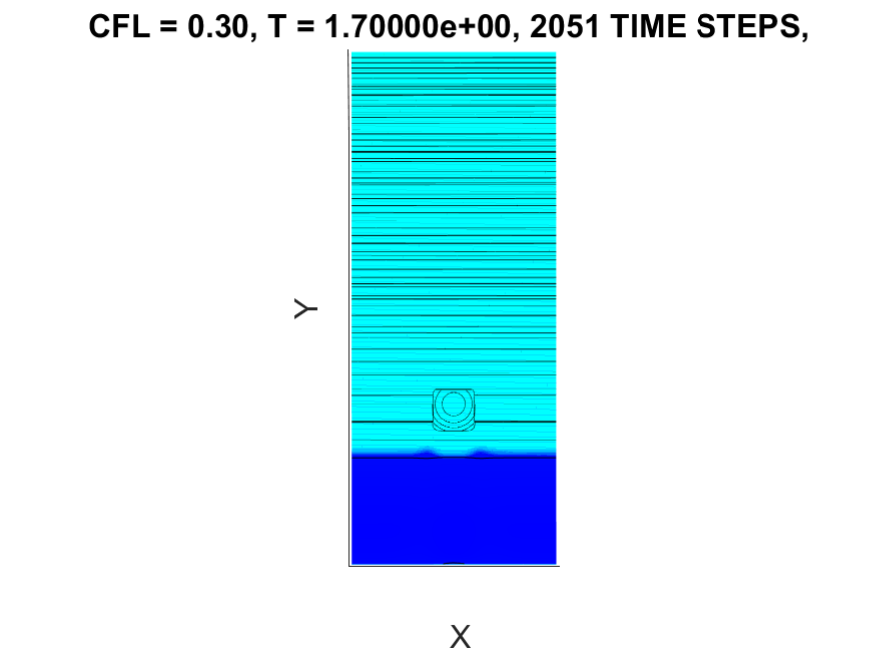}

  \caption{Test case 2.2 - parabola with obstacle: Well-balance
    simulation. The left column shows the behavior of the norms of
    $I(\bSource)$ in the interval $[0,T]$ for $T=0.5, 1.0, 1.7$~s. The
    central and right columns show the solution in terms of depth in
    the same time intervals. The contours on the right column
    emphasize the symmetry of the numerical solution with respect to
    the flow direction.}
  \label{fig:parabola2}
\end{figure}

\begin{figure}
  \centering
  \includegraphics[width=0.3\textwidth]{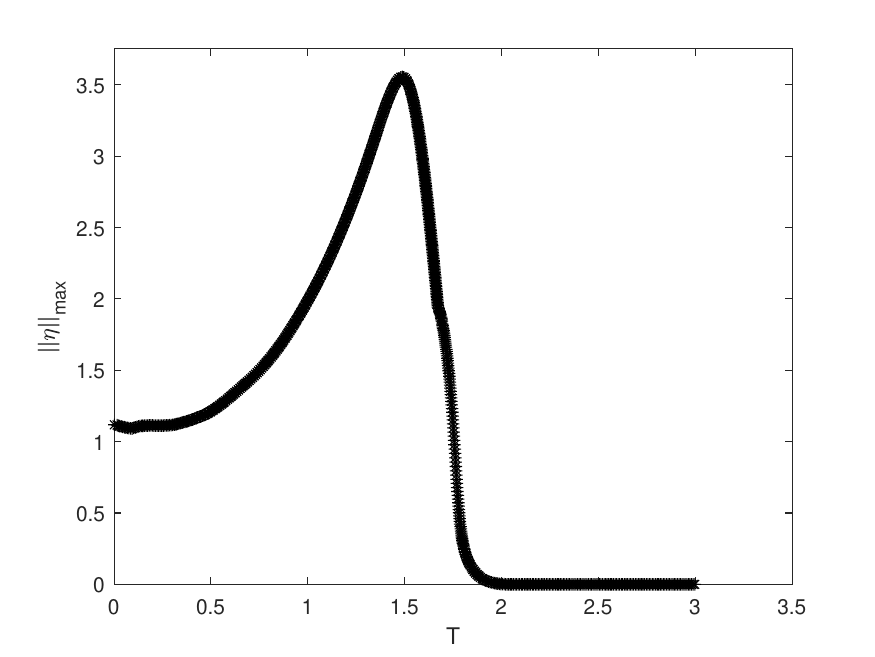}
  \includegraphics[width=0.3\textwidth]{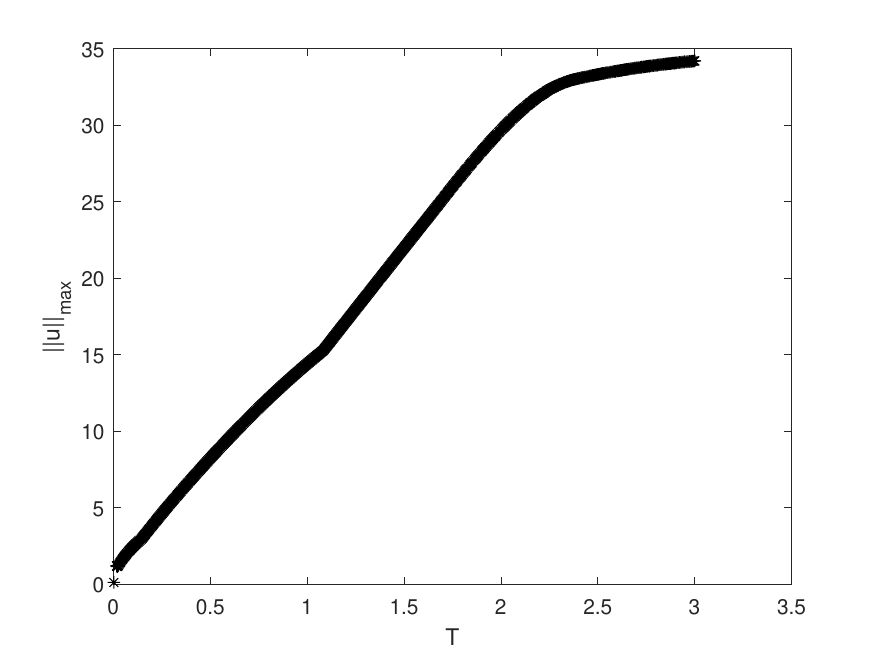}
  \includegraphics[width=0.3\textwidth]{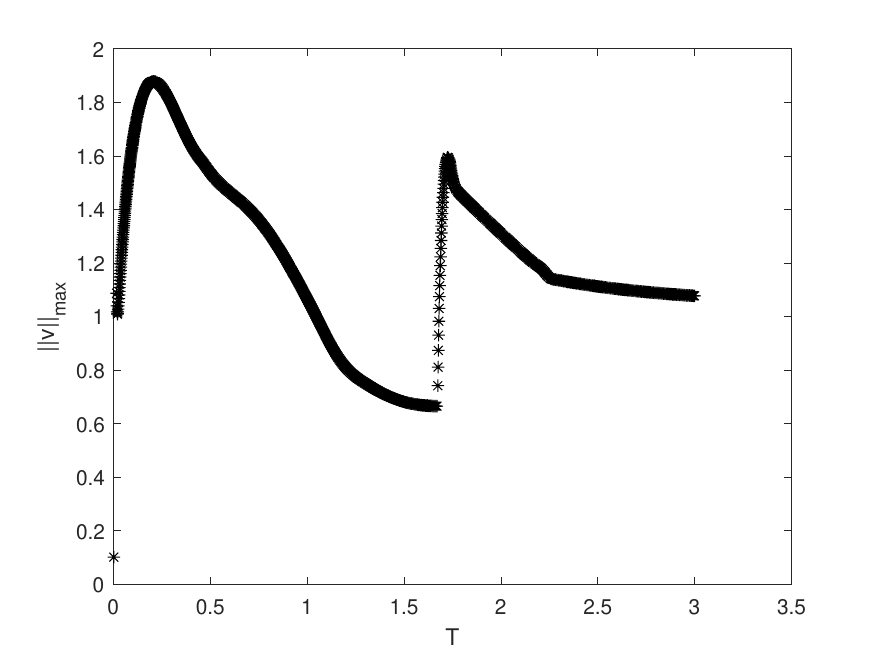}
  \includegraphics[width=0.3\textwidth]{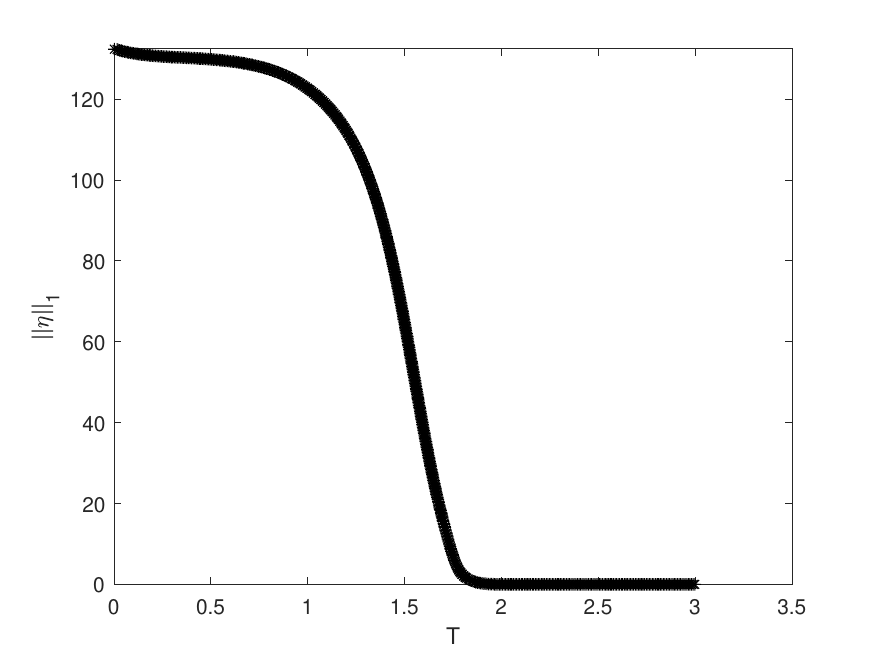}
  \includegraphics[width=0.3\textwidth]{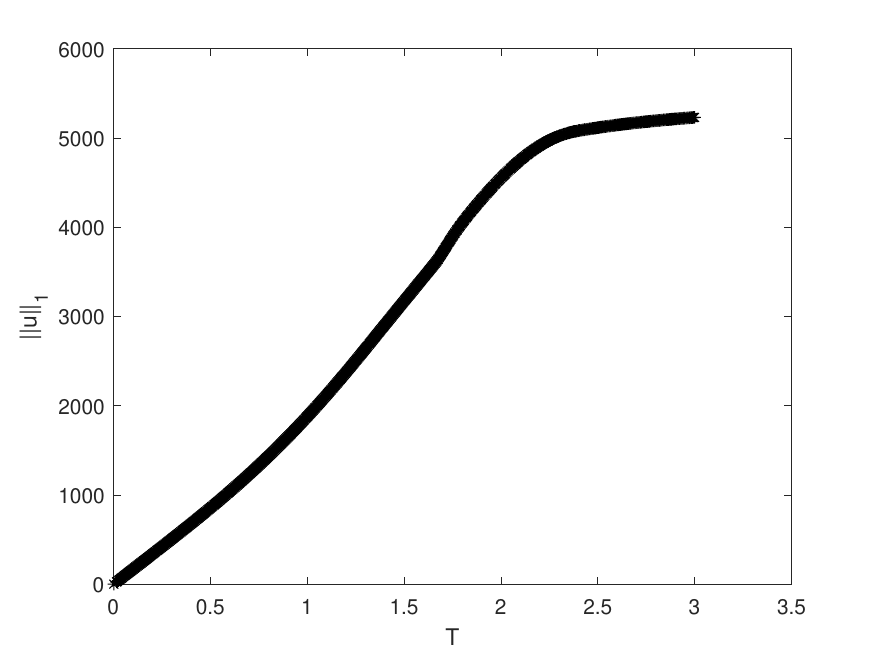}
  \includegraphics[width=0.3\textwidth]{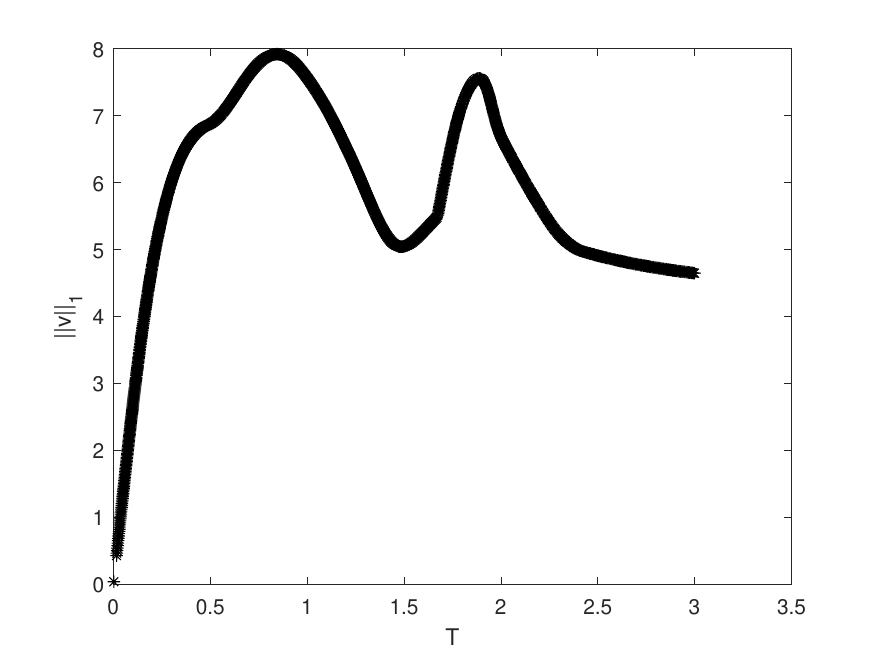}
  \caption{Test case 2.2 - parabola with localized obstacle:
    Wetting-and-drying simulation. Time behavior of $\Lspace[\infty]-$
    and $\Lspace[1]-$norms of the conserved variables.}
  \label{Hparabola}
\end{figure}

\autoref{fig:parabola1} shows the water depth distributions in the time
interval $T\in[0,1.76]$~s.
Similarly as in the previous case, the fluid flow is well captured
along the bottom surface.
Due to the geometric features of the parabola, as time evolves, dry
regions start to appear already after 0.6~s, practically emptying the
domain at the end of the simulation.
Again, the interaction with the obstacle is evident. Looking at the
right panels of \autoref{fig:parabola2}, which reports the same
measurements of \autoref{fig:Slopeplane2}, we observe the symmetric
behavior of the water flowing around the obstacle.
The behavior at long-times is shown in \autoref{Hparabola}, where the
$\Lspace[\infty]$ (top row) and $\Lspace[1]$ (bottom row) spatial
norms of $\eta$, $u$ and $v$ are reported in the time interval
$[0,3]$~s.  As seen in test case 1, as soon as the depth tends to zero
(left panels) the velocities reach their limit value (central and
right panels), indicating an emptying domain.

\paragraph{Test case 2.3: hyperboloid-central-bump with obstacle.}
The third example considers a surface defined by the height function:
$
  \BSM(x,y)=-4\sqrt{x^2+y^2+1}/5
$
on a two-dimensional square domain
$\SubsetU=[-3,3]\times[-3,3]\subset\REAL^2$ discretized using
$256\times256$ cells with the discontinuous obstacle given by
$\obstacle$ in  \autoref{eq:bump} with $\bar{z} = 1.2$, $a=0$, and $b=0$.
In this case, the underling bottom surface presents curvatures in both
directions and the test is devised to underline the ability of the
method is able to capture completely the symmetries in the solution.
We tested two different kinds of uniform boundary conditions: no flow
(first simulation) and outflow (second simulation).
The initial condition is the same for both runs and is defined in
order to capture the radial symmetric feature of the bottom surface in
the long time behavior of the solution: 2~m of (normal) water depth
are considered in the central part of the domain, i.e.  inside of a
circle of radius $0.5$~m, and 1~m outside of central radius. A maximum
CFL number of $0.30$ is used.

\begin{figure}
  \centering
  \includegraphics[width=\textwidth]{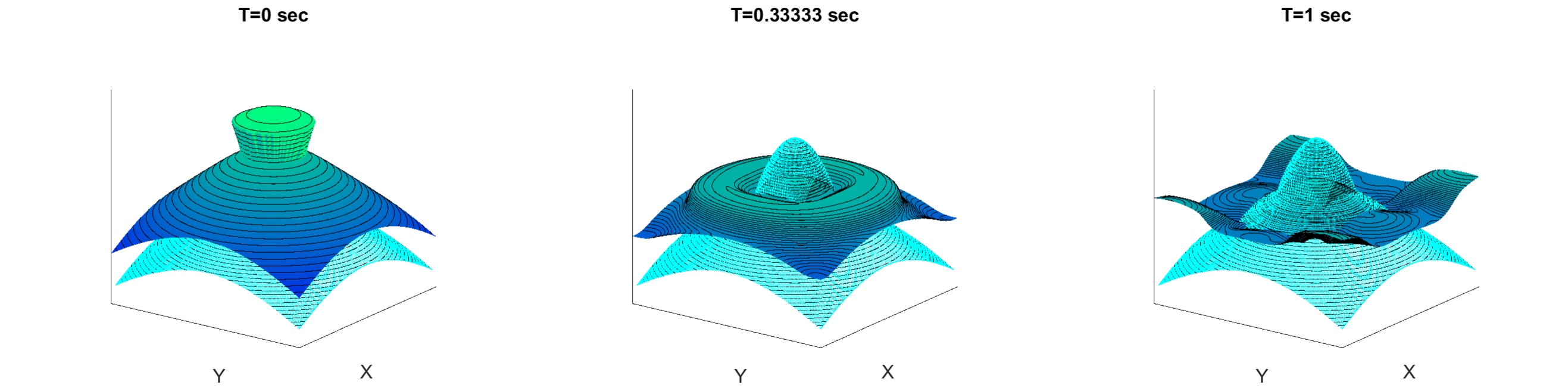}\\
  \includegraphics[width=\textwidth]{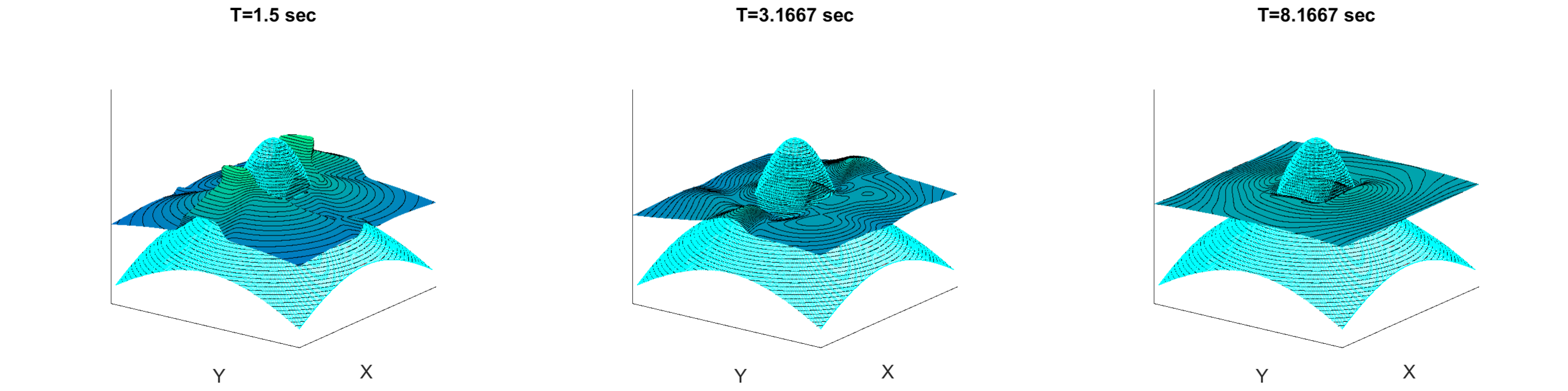}\\
  \caption{Test case 2.3a - hyperboloid-central-bump with obstacle:
    evolution of the normal water depth at different times. No-flow
    conditions are imposed on all boundaries. The computational grid
    is formed by 256$\times$256 cells. }
  \label{bottonfullyacum}
\end{figure}

\begin{figure}
  \centering
  \includegraphics[width=0.3\textwidth]{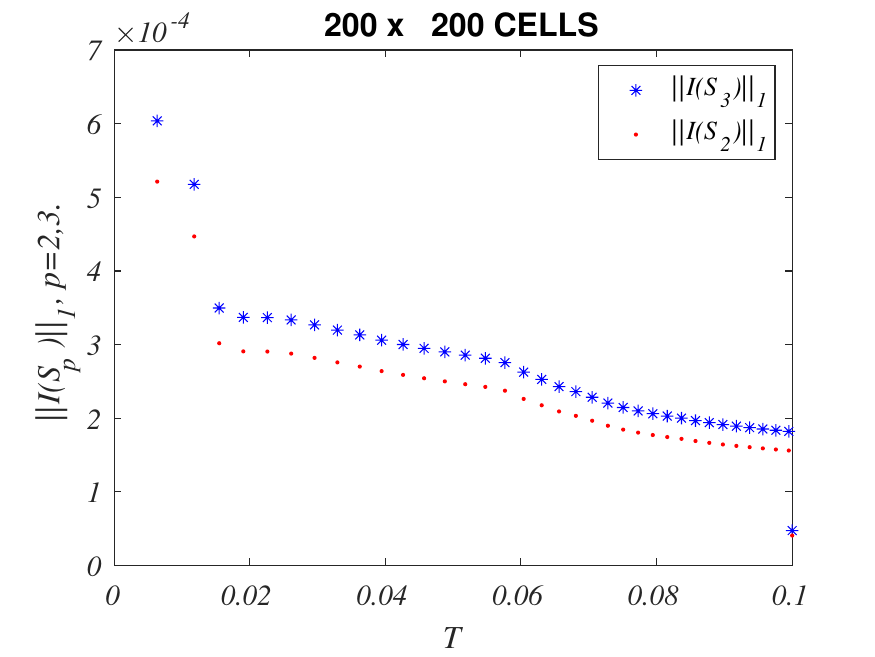}
  \includegraphics[width=0.3\textwidth]{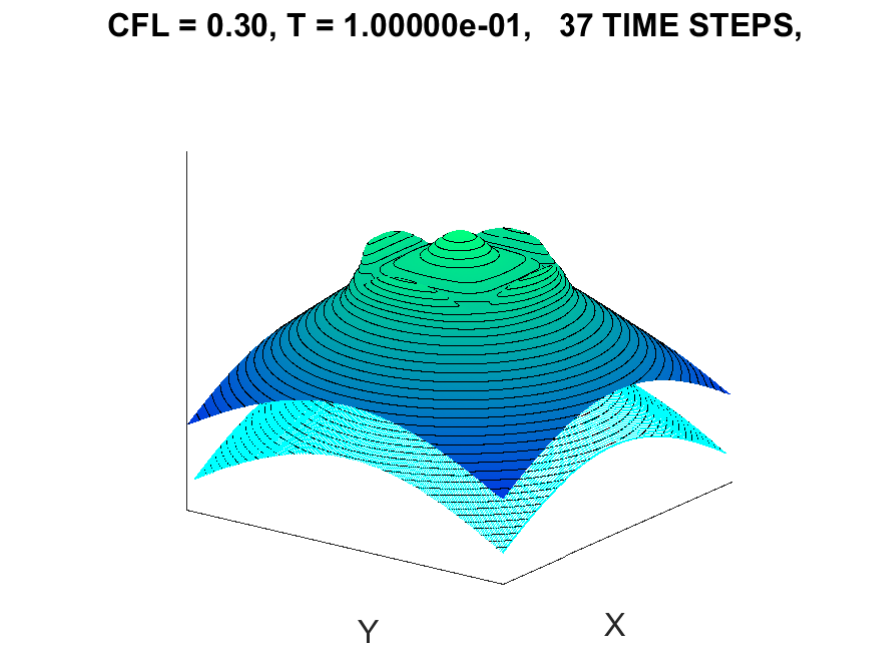}
  \includegraphics[width=0.3\textwidth]{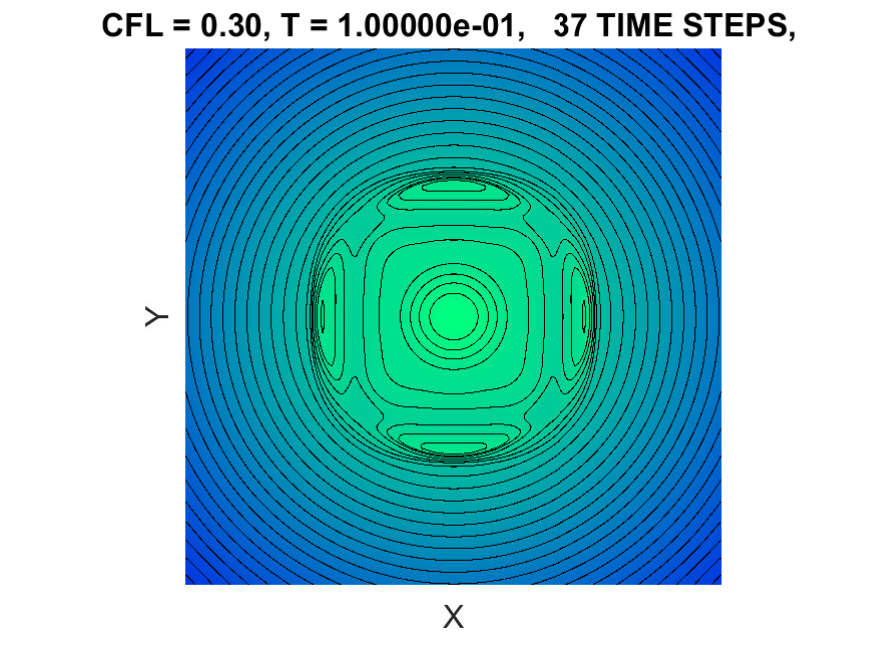}
  \includegraphics[width=0.3\textwidth]{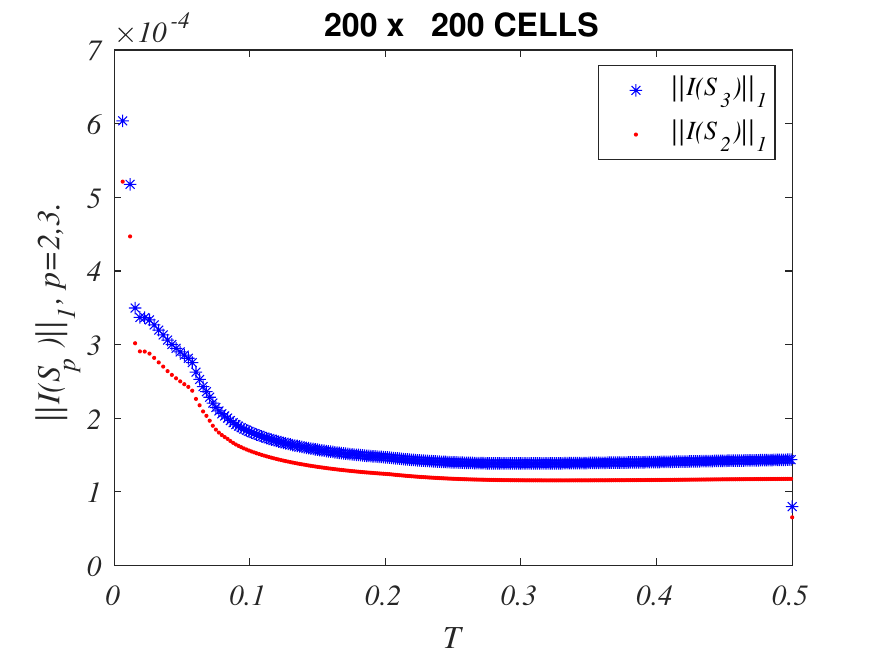}
  \includegraphics[width=0.3\textwidth]{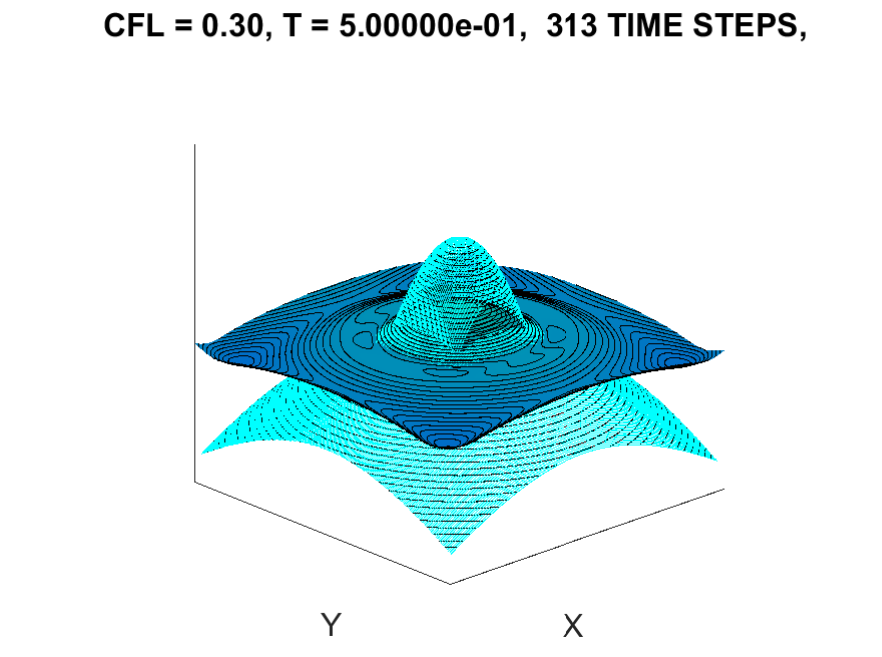}
  \includegraphics[width=0.3\textwidth]{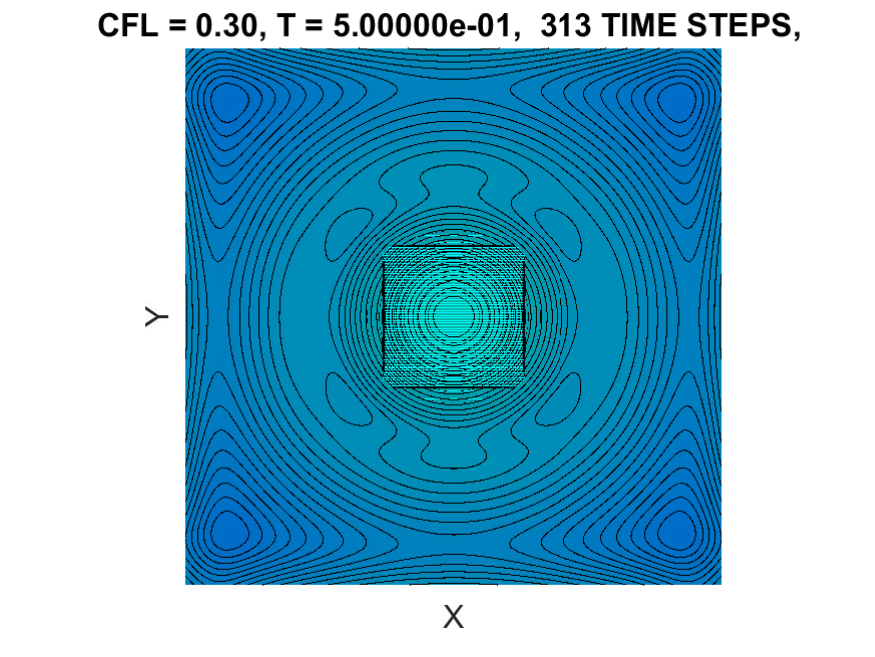}
  \includegraphics[width=0.3\textwidth]{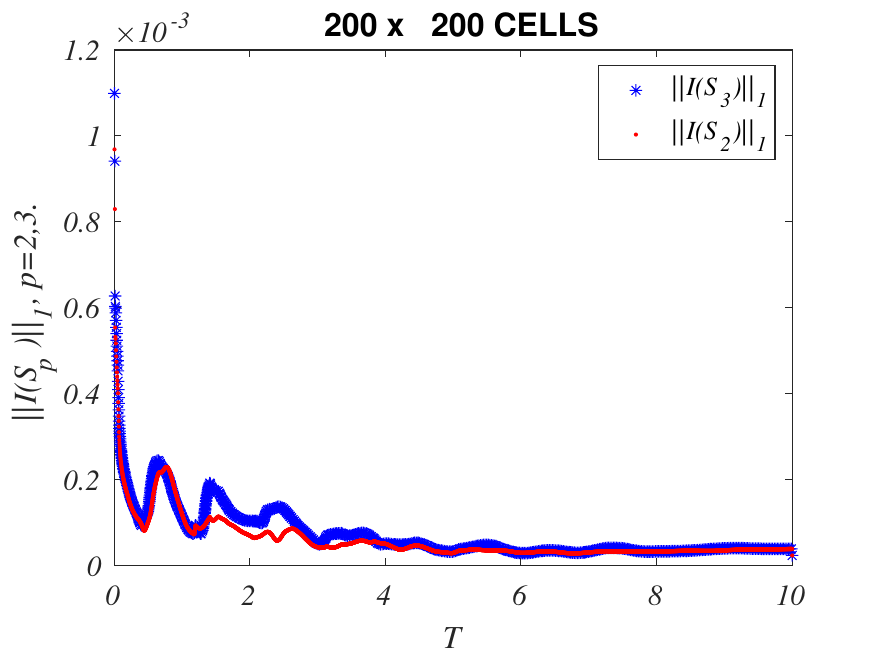}
  \includegraphics[width=0.3\textwidth]{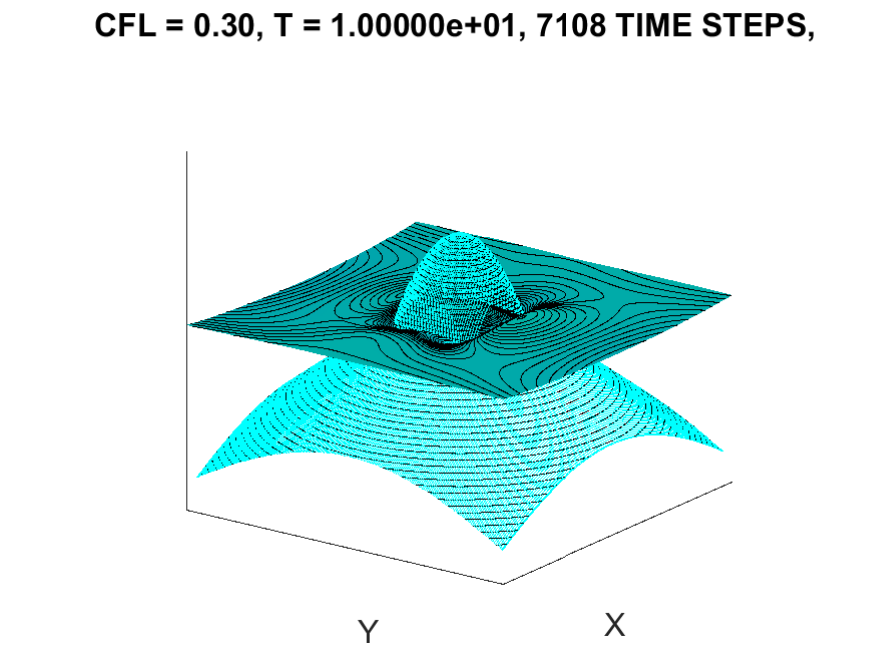}
  \includegraphics[width=0.3\textwidth]{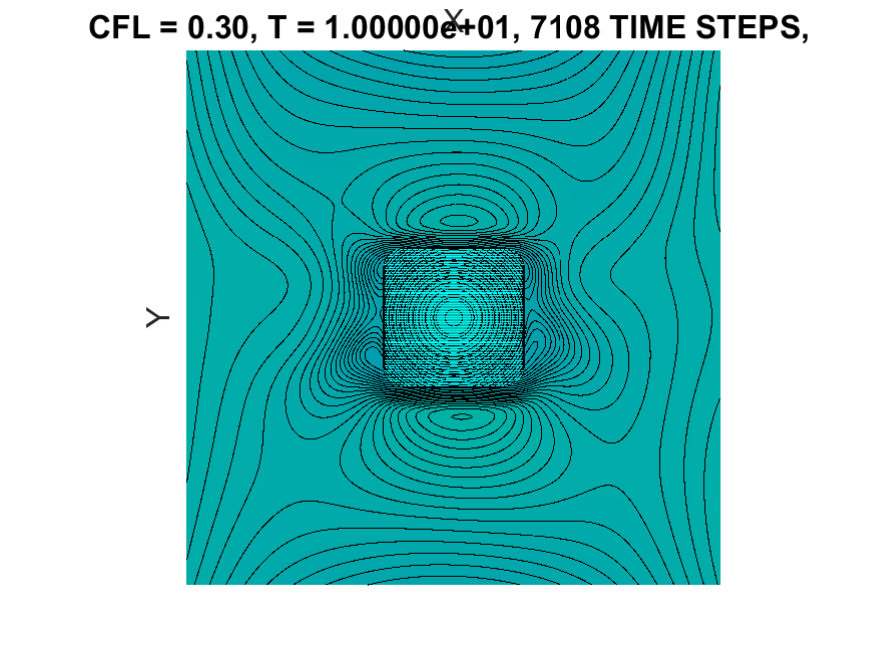}
  
  \caption{Test case 2.3a - hyperboloid-central-bump with obstacle:
    well-balance simulation. The left column shows the time behavior
    of the $L_1$ norms of $I(\bSource)$ in the interval $[0,T]$ for
    $T=0.1, 0.5, 10.0$~s. The central and right columns show the
    solution in terms of depth at the same times. The contours on the
    right column emphasize the symmetry of the numerical solution with
    respect to the flow direction.}
  \label{DbottonfullyWellBalacum}
\end{figure}

\autoref{bottonfullyacum} shows the distribution of the normal water
depth at different times in the time interval $T\in[0,10]$~s,
\autoref{DbottonfullyWellBalacum} shows the state at times $T=0.1, 0.5$~s
and $T=10$~s. This is the case where no-flow boundary conditions have
been imposed at all boundaries The results show that the long-term
steady-state condition is achieved and maintained without the
formation of spurious waves.
Due to the geometric features of the bottom surface, it is possible to
appreciate the central radial symmetry of the waves as they
  dissipate in time due possibly to the dissipative character of ISWE
  or to numerical viscosity or both (see also the discussion in
  \autoref{rem:dissipation}).

\begin{figure}
  \centering
  \includegraphics[width=\textwidth]{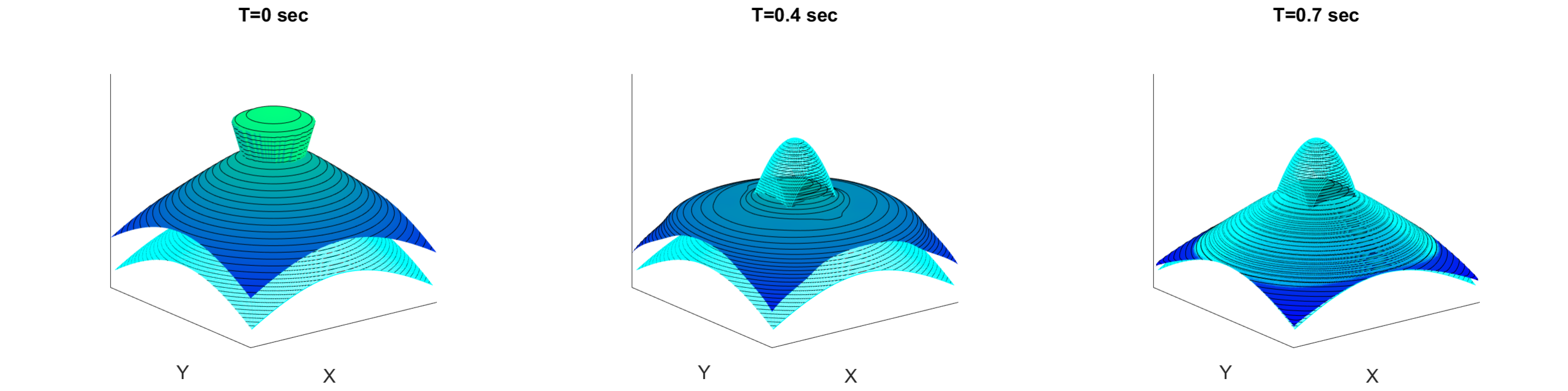}\\
  \includegraphics[width=\textwidth]{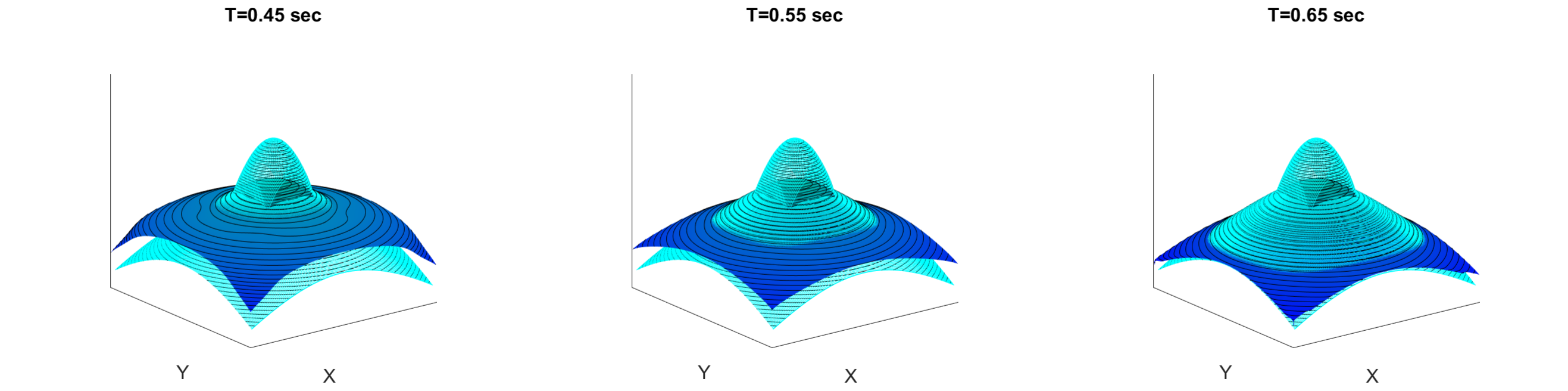}\\
  \caption{Test case 2.3b - hyperboloid-central-bump with obstacle:
    evolution of the water depth at different times. Outflow boundary
    conditions are considered for all boundary. The grid is formed by
    256$\times$256 cells.}
  \label{HCB1}
\end{figure}

\begin{figure}
  \centering
  \includegraphics[width=0.3\textwidth]{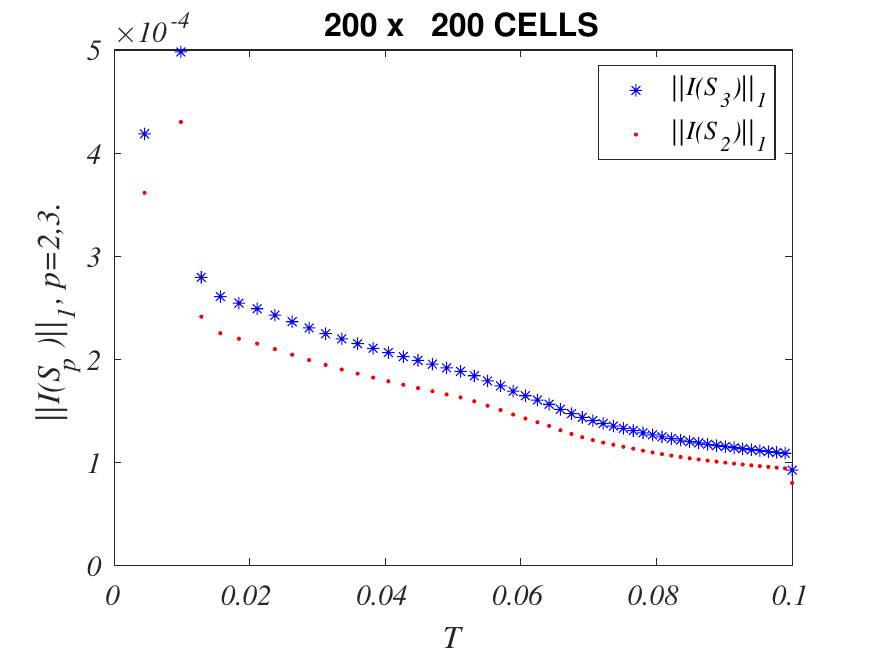}
  \includegraphics[width=0.3\textwidth]{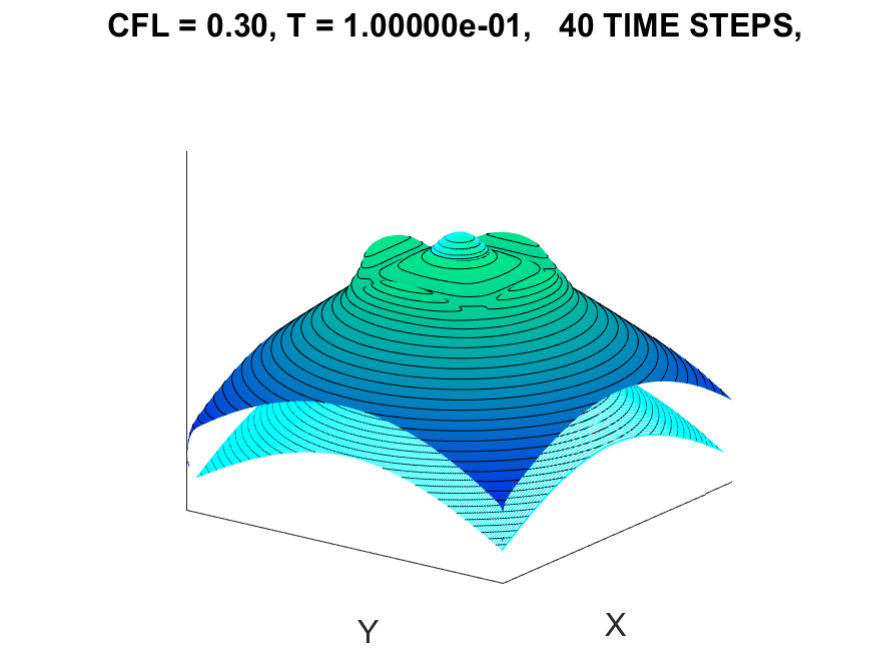}
  \includegraphics[width=0.3\textwidth]{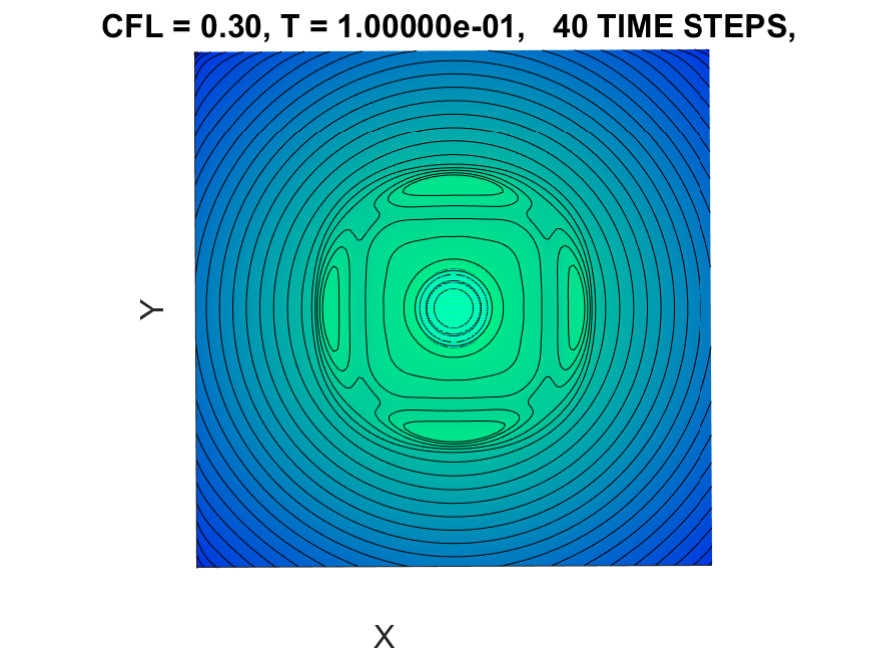}
  \includegraphics[width=0.3\textwidth]{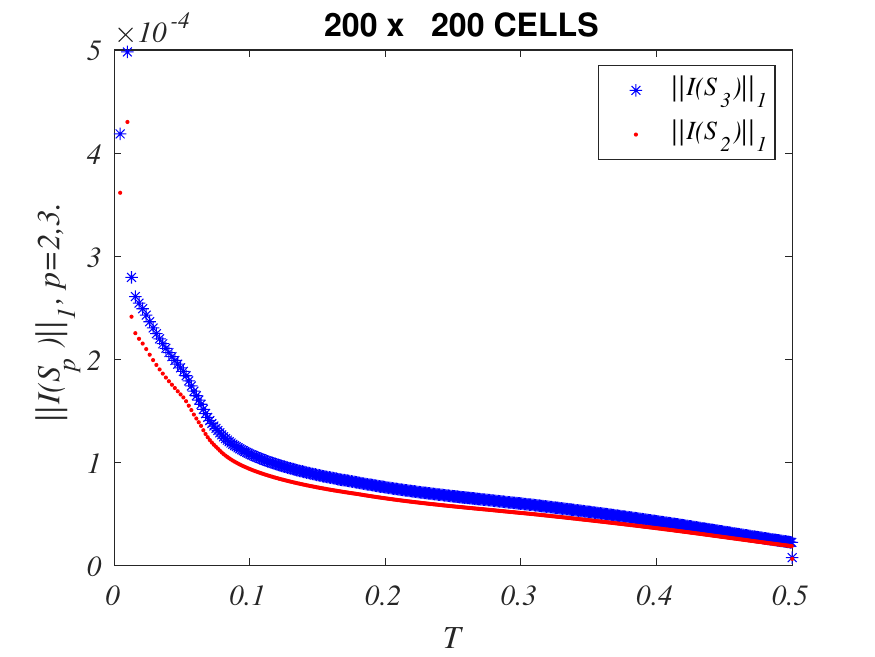}
  \includegraphics[width=0.3\textwidth]{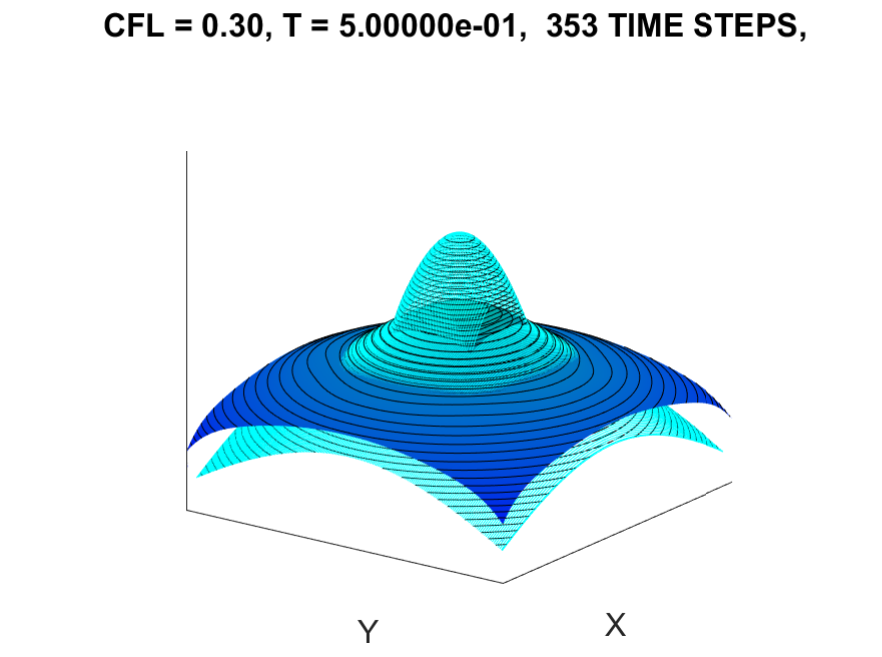}
  \includegraphics[width=0.3\textwidth]{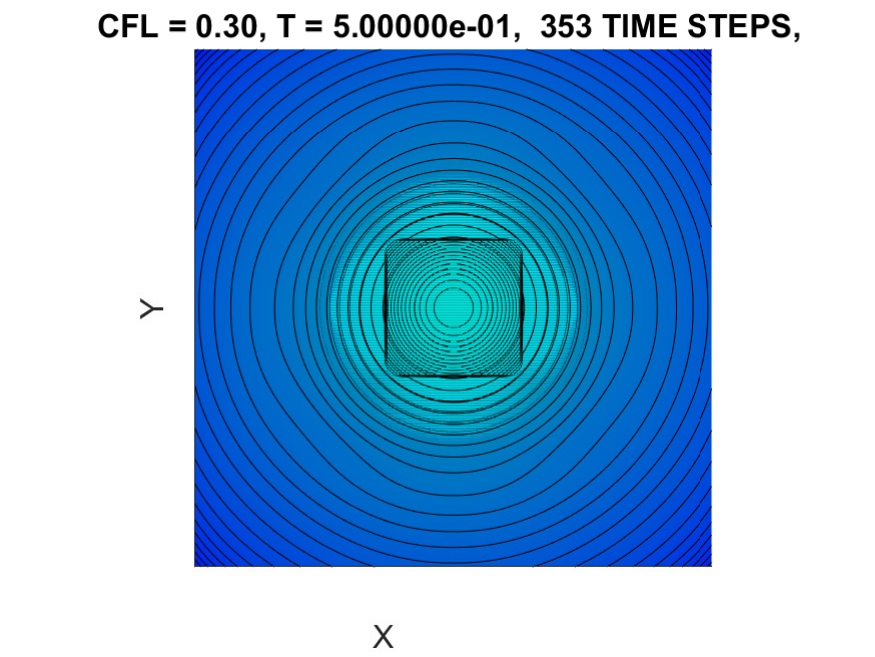}
  \includegraphics[width=0.3\textwidth]{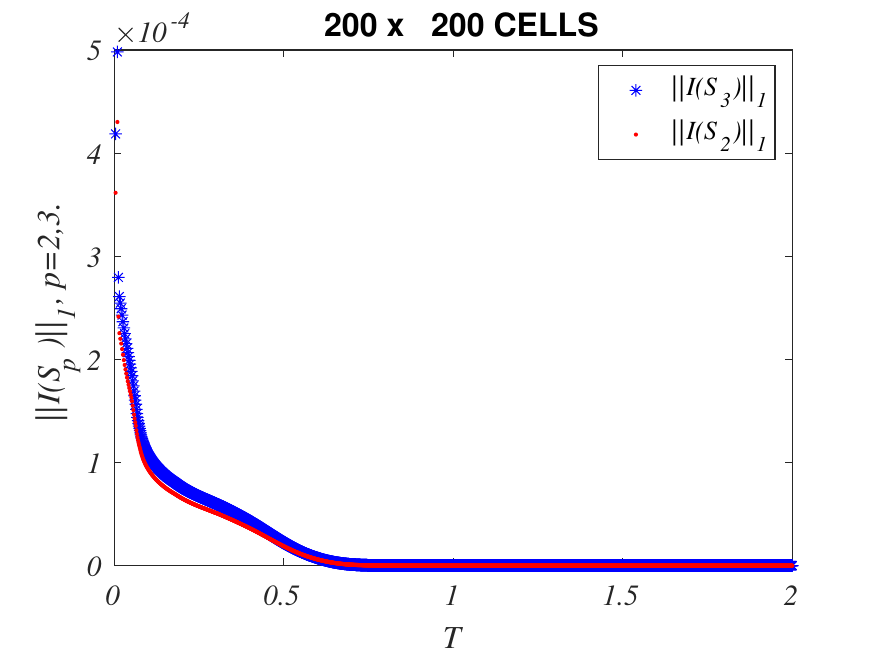}
  \includegraphics[width=0.3\textwidth]{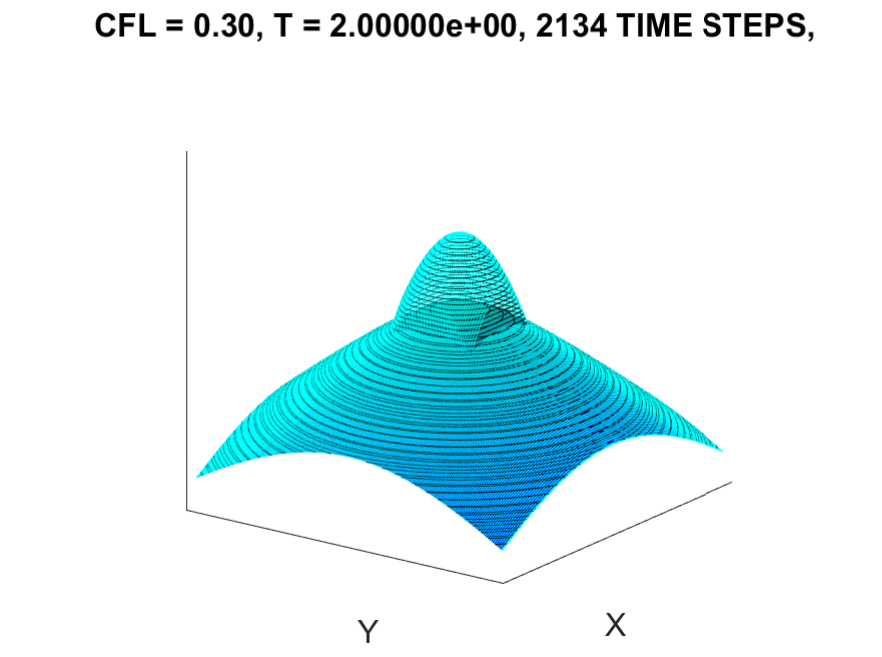}
  \includegraphics[width=0.3\textwidth]{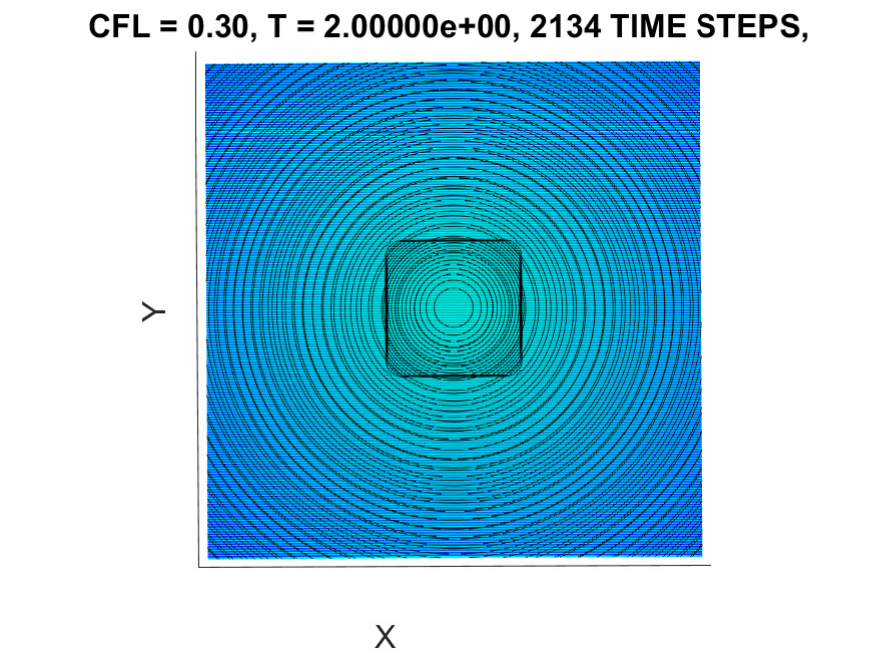}
  \caption{Test case 2.3b - hyperboloid-central-bump with obstacle:
    Well-balance simulation. The left column shows the behavior of the
    norms of $I(\bSource)$ in the interval $[0,T]$ for
    $T=0.1, 0.5, 2.0$~s. The central and right columns show the
    solution in terms of depth in the same time intervals. The
    contours in the right panels emphasize the symmetry of the
    solution with respect to the flow direction.  }
  \label{DbottonfullyWellBaldry1}
\end{figure}

\begin{figure}
  \centering
  \includegraphics[width=0.3\textwidth]{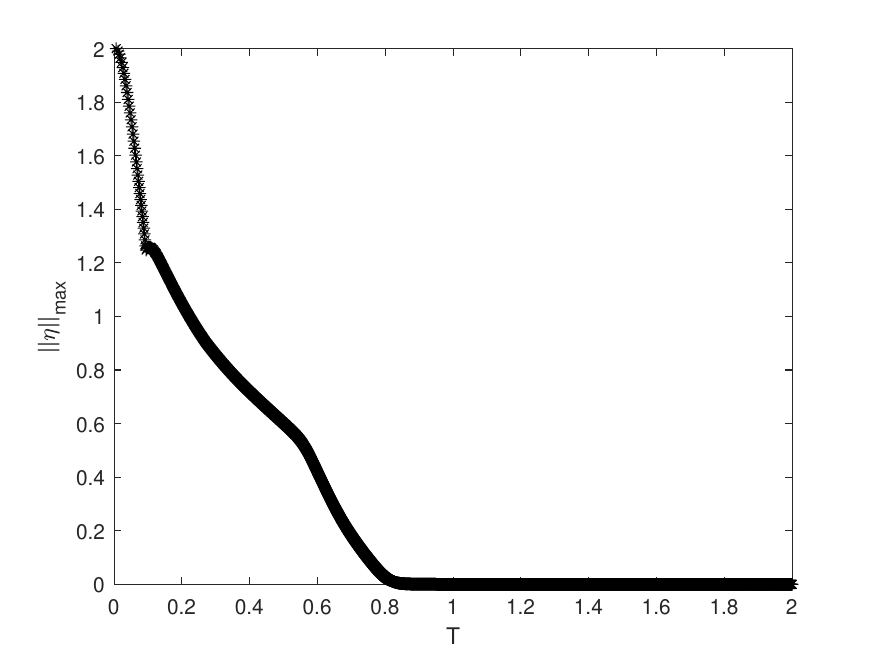}
  \includegraphics[width=0.3\textwidth]{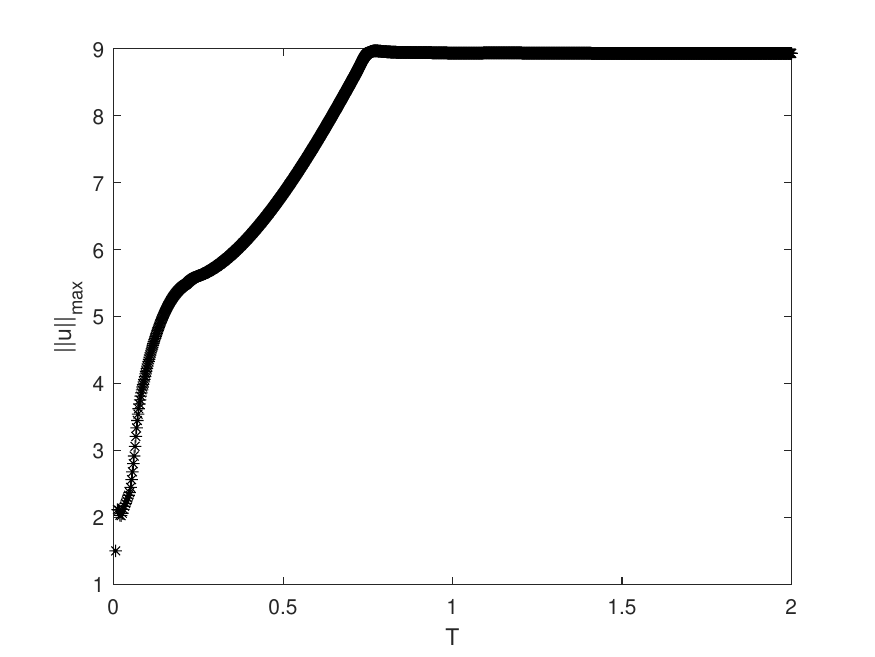}
  \includegraphics[width=0.3\textwidth]{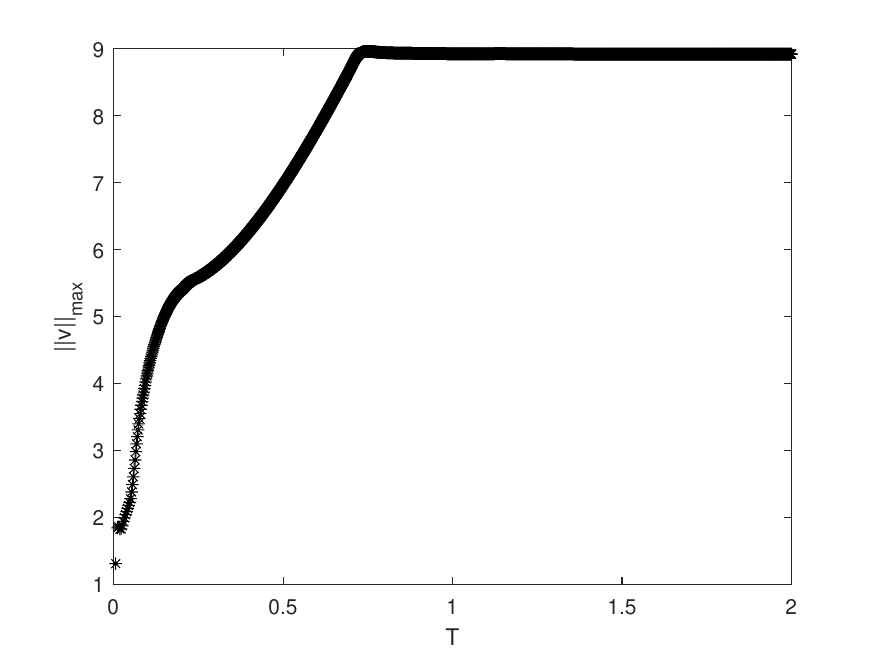}
  \includegraphics[width=0.3\textwidth]{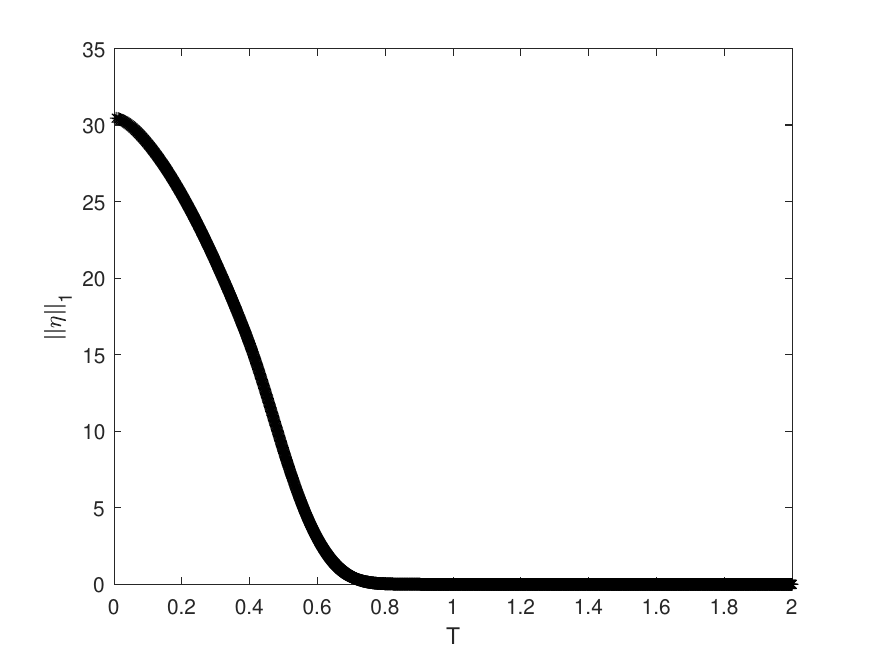}
  \includegraphics[width=0.3\textwidth]{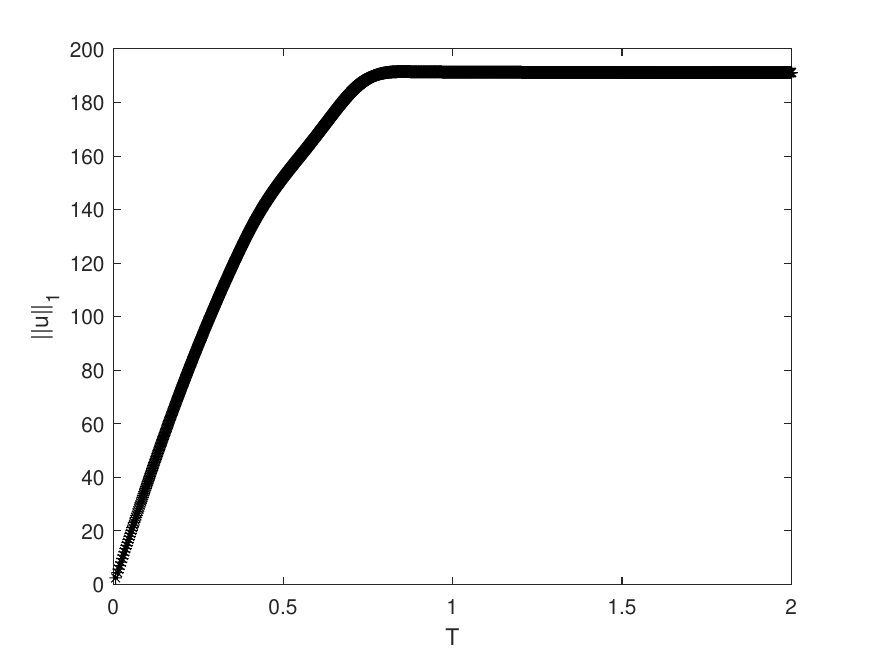}
  \includegraphics[width=0.3\textwidth]{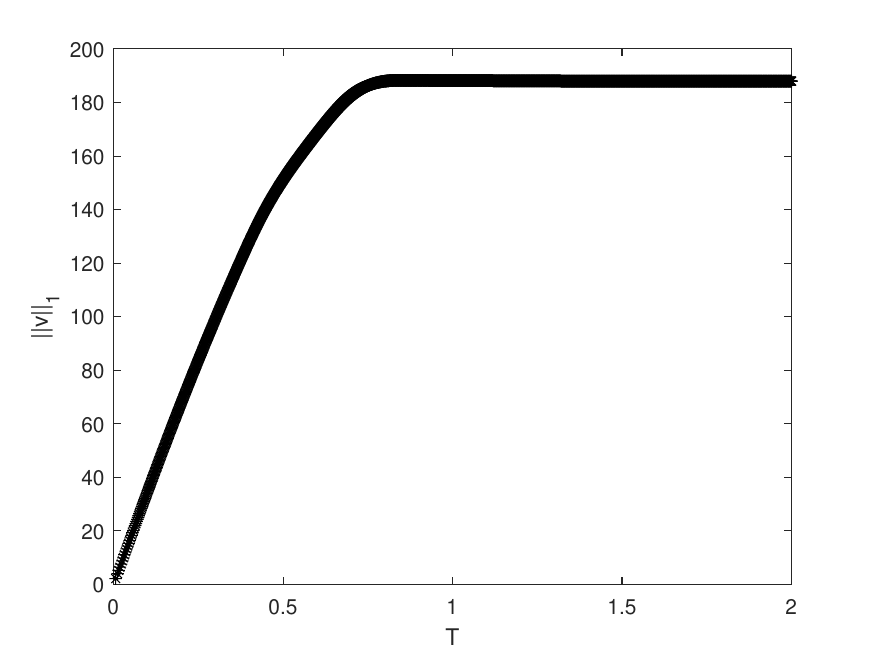}

  \caption{
    Test case 2.3b - hyperboloid-central-bump with obstacle:
    wetting-and-drying simulation. Time behavior of $\Lspace[\infty]-$
    and $\Lspace[1]-$norms of the conserved variables.
    We note that, as soon as the depth tends to zero (see left panels),
    the velocities reach their vanishing limits (see central and right
    panels). Moreover, a symmetric behavior is noted also in the
    values of the limit for the velocities in the $\xcl-$ and
    $\ycl-$directions.}
  \label{HCB-wetdry}
\end{figure}

The results for the outflow run are shown in
 \autoref{HCB1}, \autoref{DbottonfullyWellBaldry1}, and \autoref{HCB-wetdry}. In this case, the
fluid leaves the domain very quickly and at $T=2$~s the domain has
been empty for some time. The results show that the expected
symmetries dictated by the domain with the obstacle are respected and
that the scheme continues to function without generating spurious
oscillations well beyond the time when the fluid fully exits the
domain.
Finally, \autoref{HCB-wetdry} shows both $\Lspace[\infty]$ (top panels)
and $\Lspace[1]$ norms (bottom panels) of $\eta$, $u$ and $v$ in the
long-time behavior until $T=2$~s, certifying the robustness of the
proposed scheme.

\paragraph{Test case 2.4: fully 3D surface with a bump.}
The final test of this section is designed to highlight the behavior
of the proposed scheme when positive and negative curvatures are
present and to test its effectiveness in handling complex
wetting-and-drying patterns.  The height function for the fourth
bottom surface is given by:
$
  \BSM(x,y) = -\frac{1}{500}(y)^3 - \frac{1}{100}y(x)^2
$
on a domain $\SubsetU=[-4,4]\times[-10,10]\subset\REAL^2$ discretized
on a $256\times640$ grid with obstacle given by $\obstacle$ in
\autoref{eq:bump} with $\bar{z}=1.2$, $a=0$, and $b=4$.
It is a dam-break problem with the dam located at $y=-8.5$~m and
initial condition for the normal depth set to $2.0$~m upstream and
$1.0$~m of downstream.

\begin{figure}
  \centering
  \includegraphics[width=\textwidth]{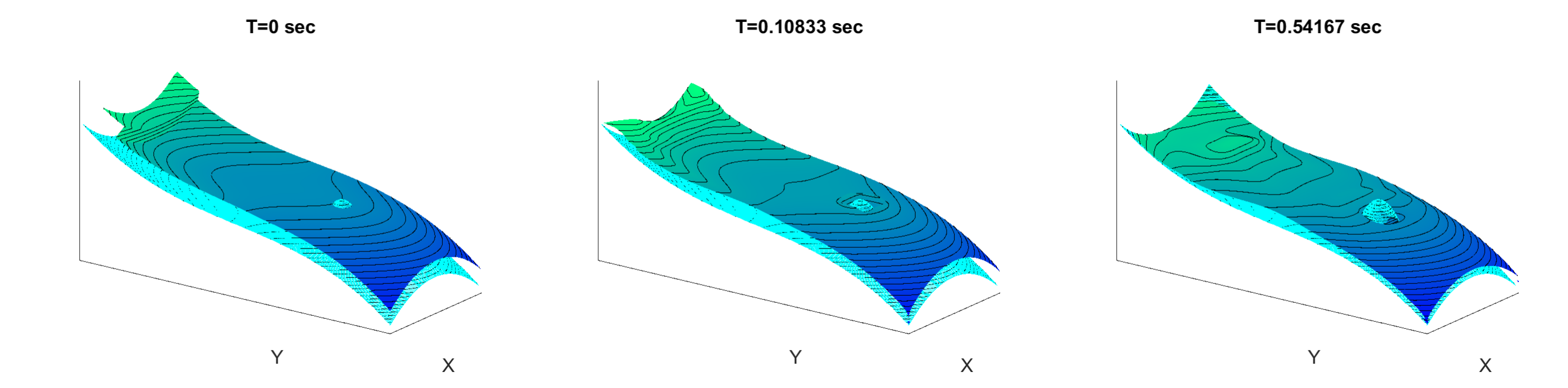}\\
  \includegraphics[width=\textwidth]{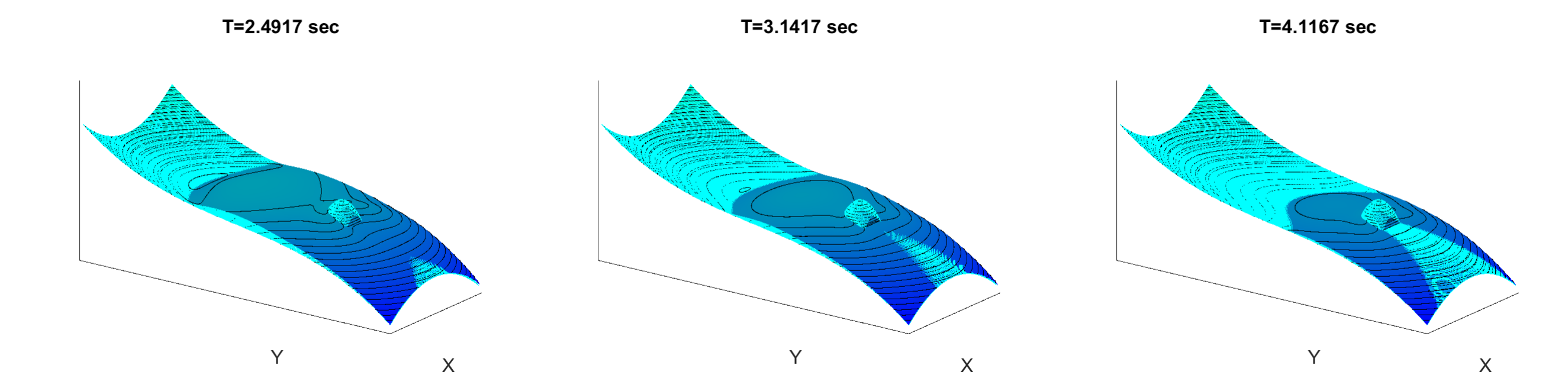}\\
  \caption{Test case 2.4 - fully 3D surface with obstacle: evolution
    of the water depth at different times. The grid is formed by
    256$\times$640 cells.}
  \label{Fully3Dsurface2}
\end{figure}

\begin{figure}
  \centering
  \includegraphics[width=0.3\textwidth]{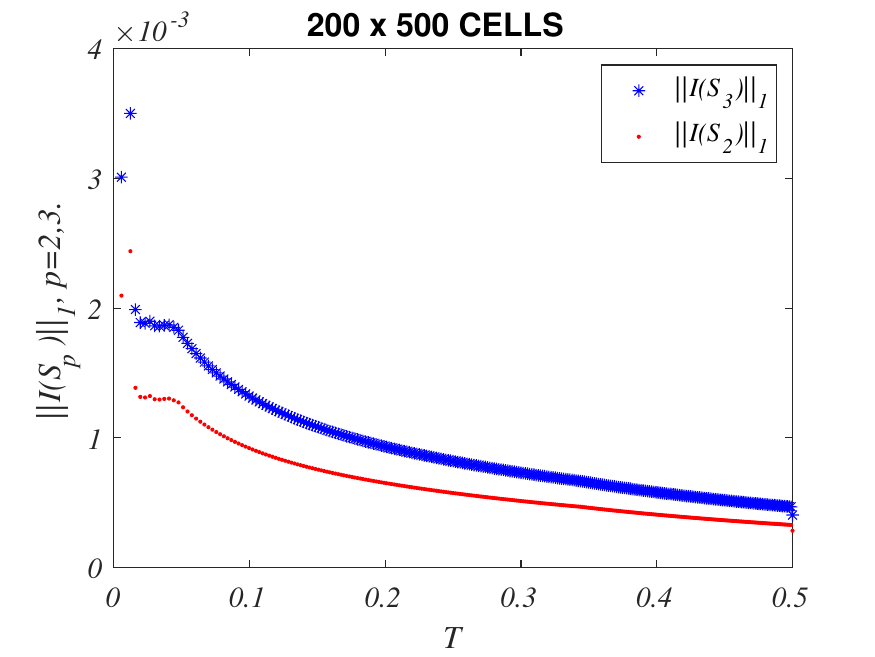}
  \includegraphics[width=0.3\textwidth]{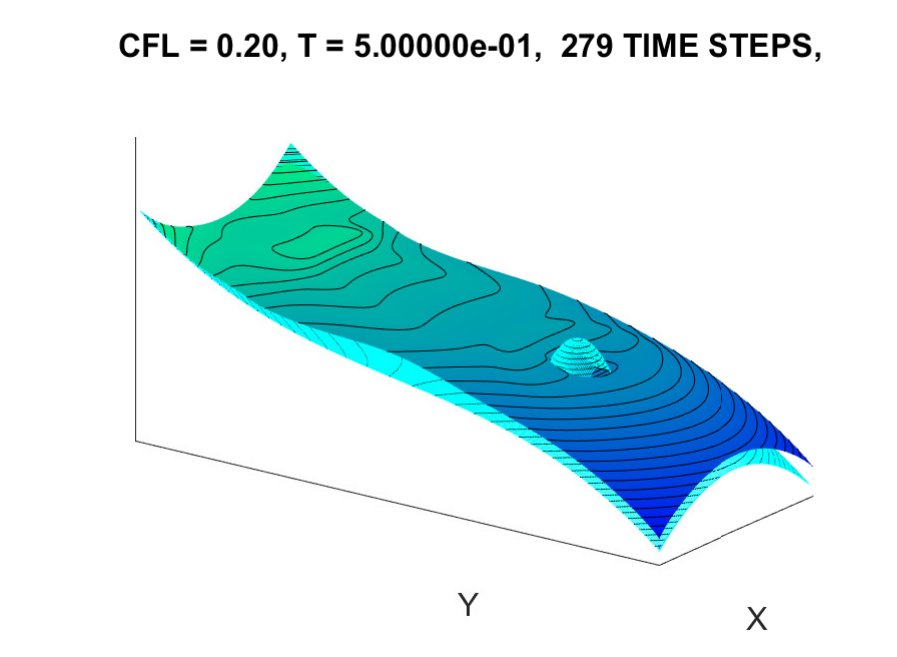}
  \includegraphics[width=0.3\textwidth]{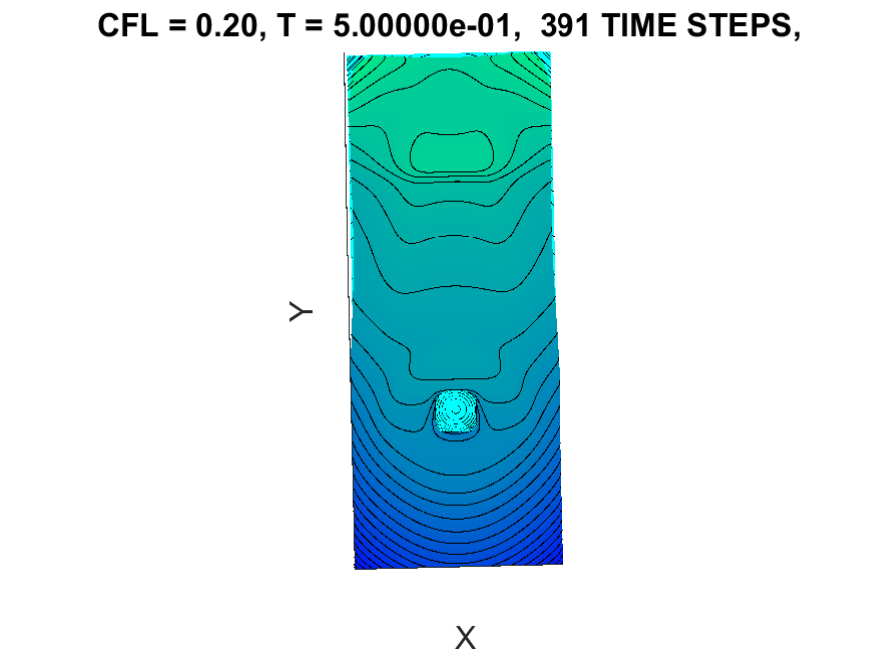}
  \includegraphics[width=0.3\textwidth]{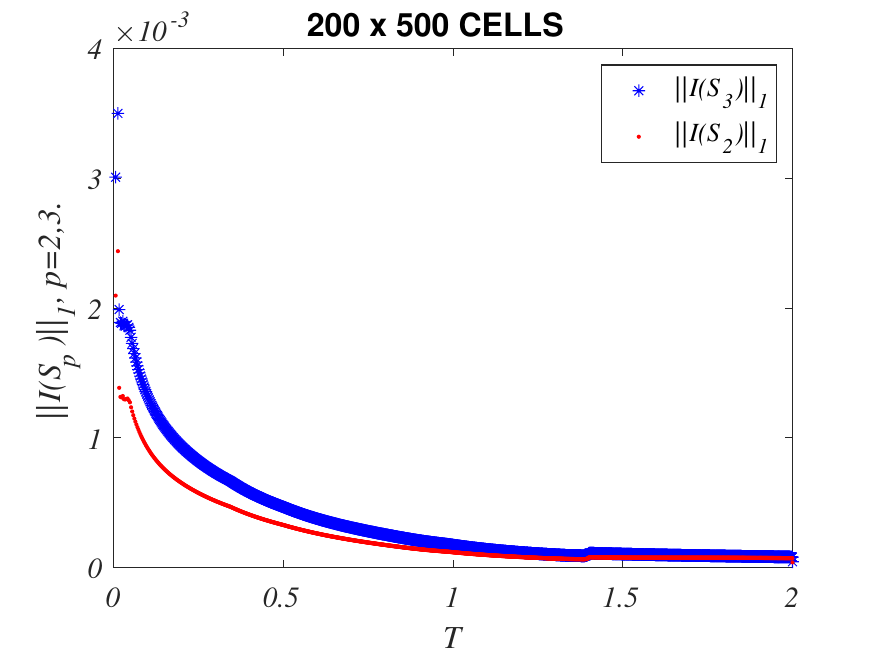}
  \includegraphics[width=0.3\textwidth]{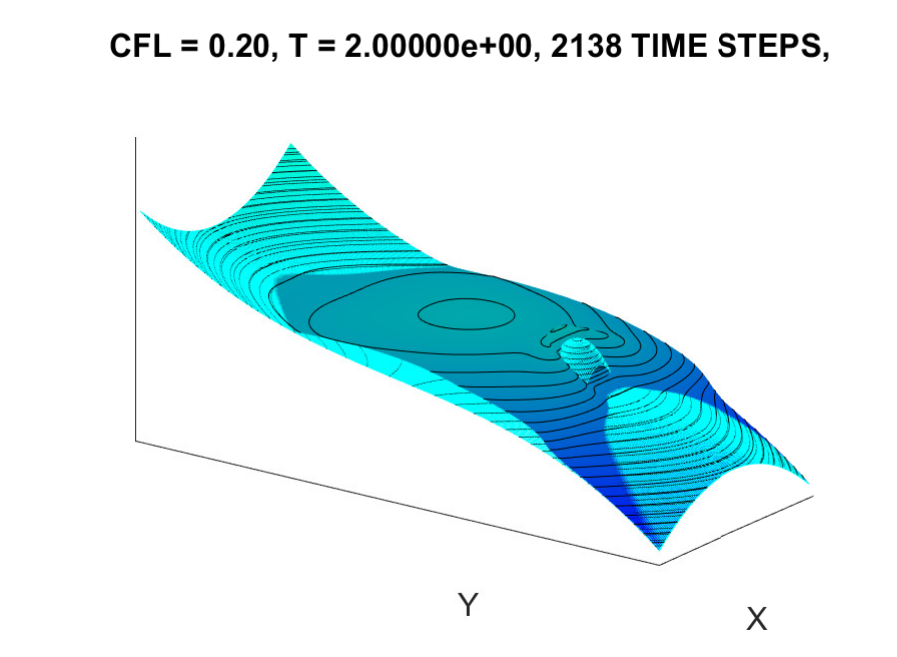}
  \includegraphics[width=0.3\textwidth]{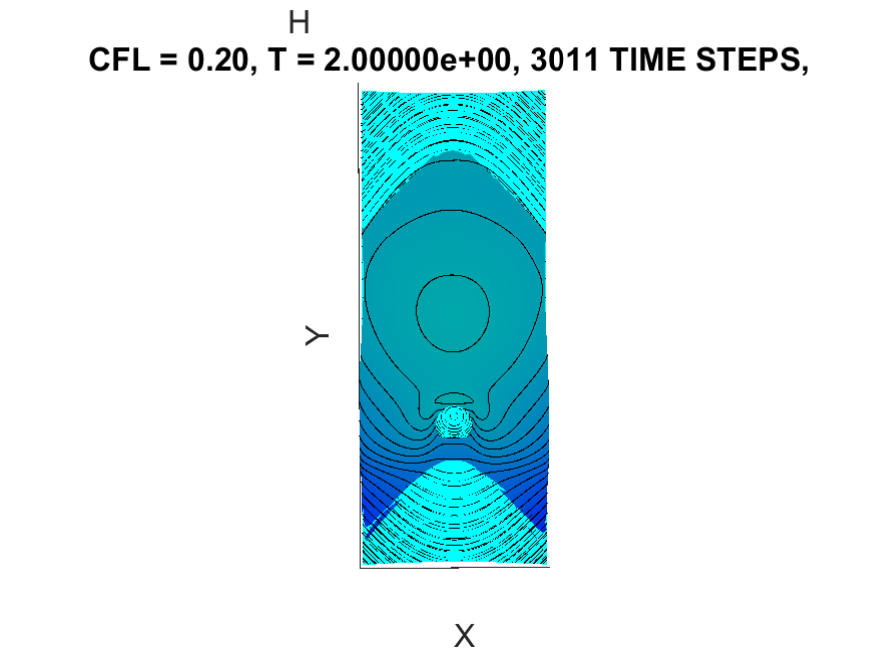}
  
  \includegraphics[width=0.3\textwidth]{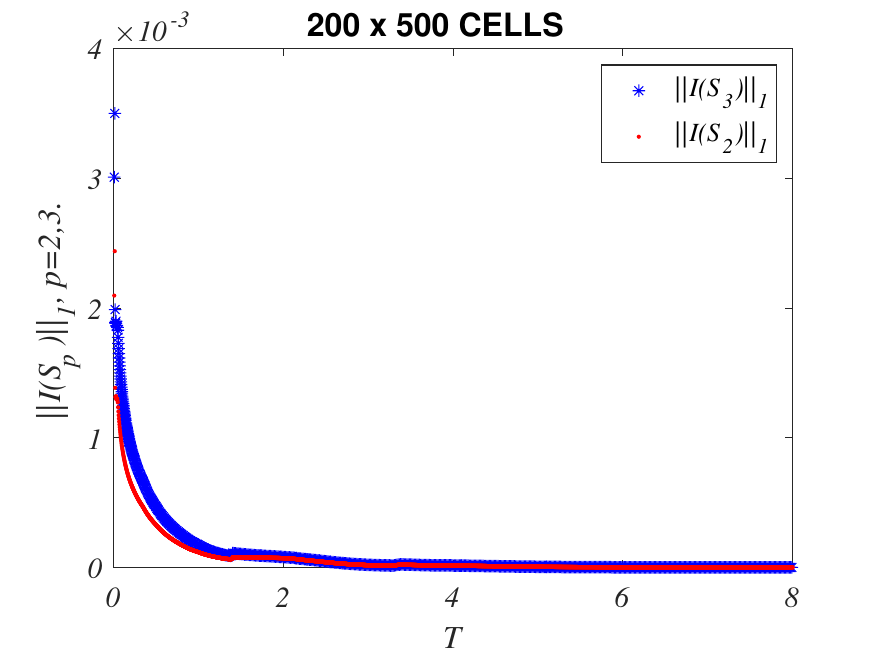}
  \includegraphics[width=0.3\textwidth]{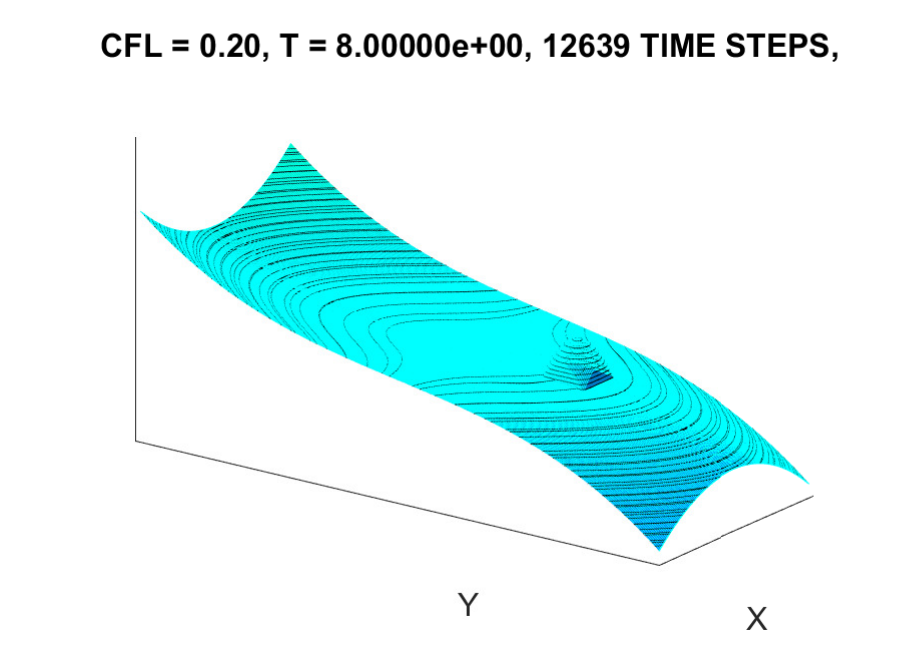}
  \includegraphics[width=0.3\textwidth]{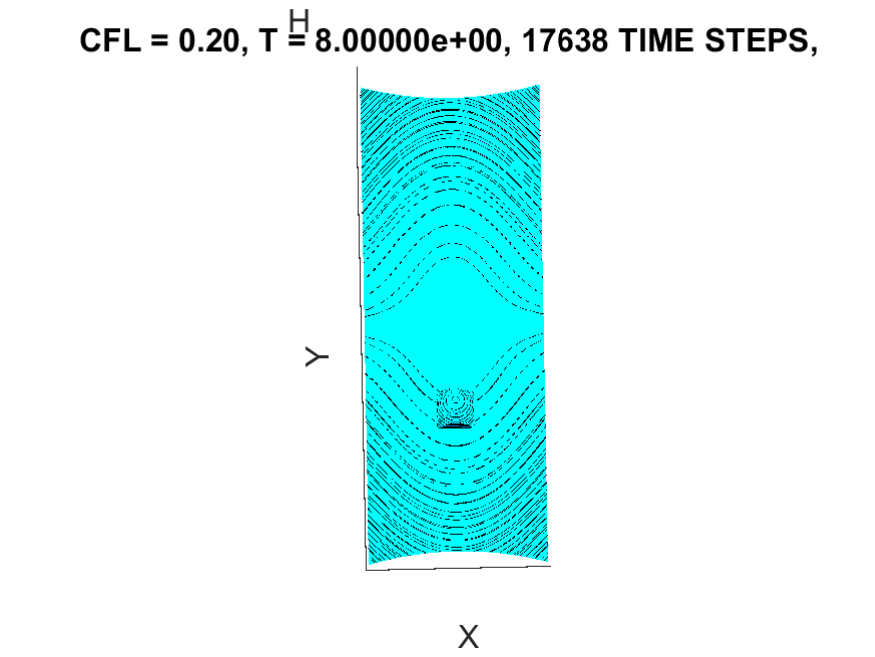}
  \caption{Test case 2.4 - fully 3D surface with obstacle:
    well-balance simulation. The left column shows the behavior of the
    norms of $I(\bSource)$ in the interval $[0,T]$ for
    $T=0.5, 2.0, 8.0$~s. The central and right columns show the solution
    in terms of depth at the same times. The contours on the right column
    emphasize the symmetry of the numerical solution with respect to
    the flow direction.}
  \label{Fully3DsurfacWB1}
\end{figure}

\begin{figure}
  \centering
  \includegraphics[width=0.3\textwidth]{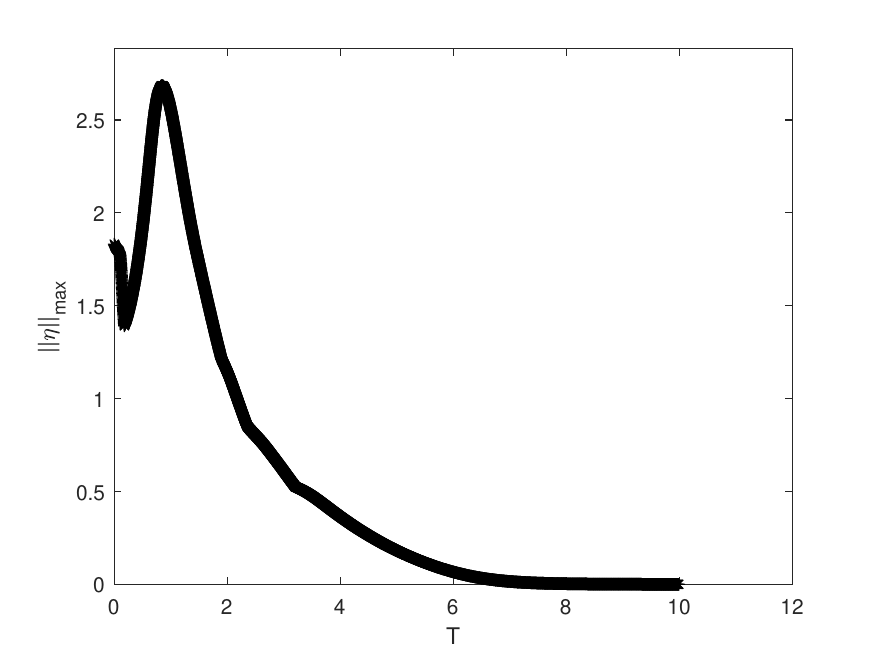}
  \includegraphics[width=0.3\textwidth]{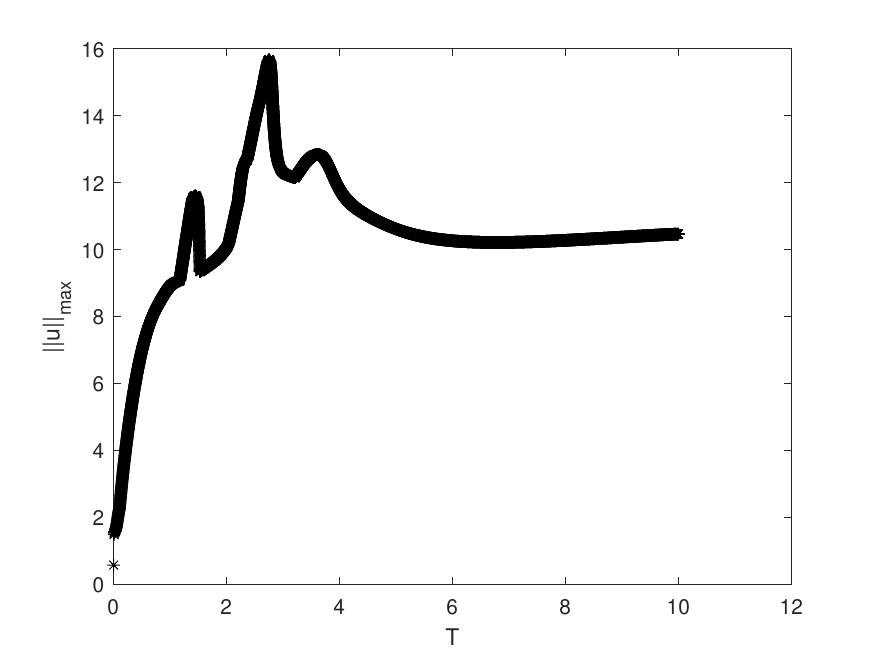}
  \includegraphics[width=0.3\textwidth]{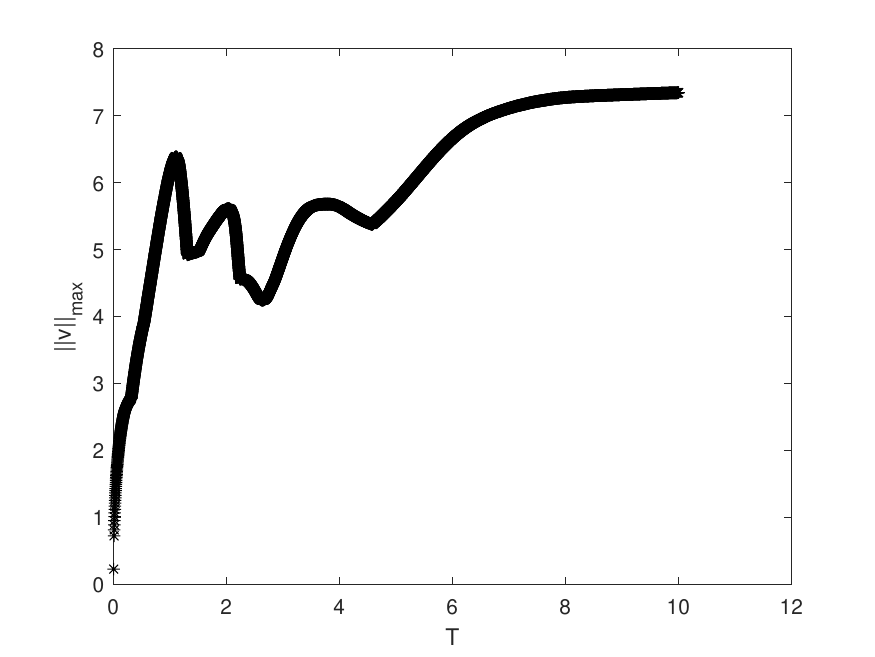}
  \includegraphics[width=0.3\textwidth]{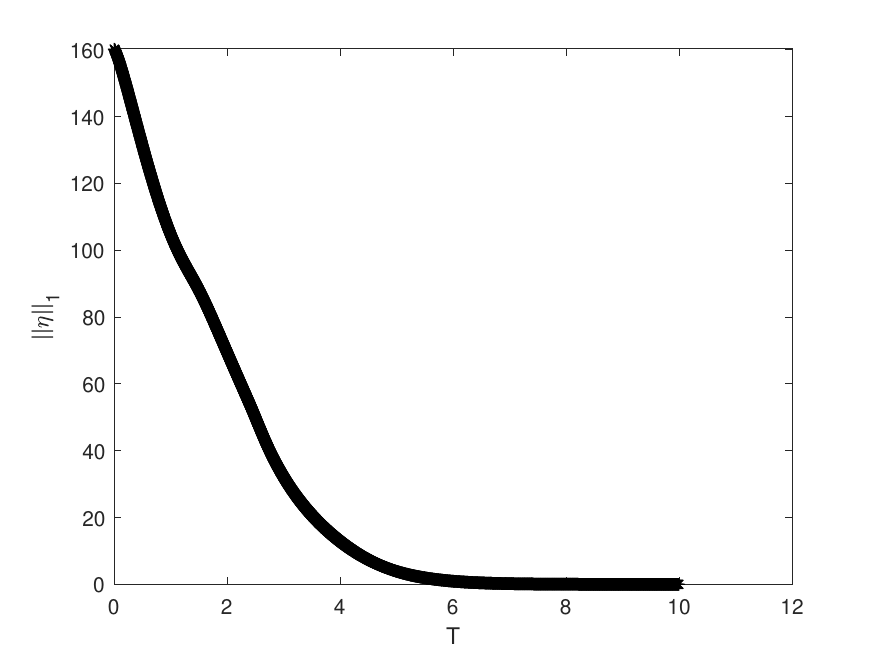}
  \includegraphics[width=0.3\textwidth]{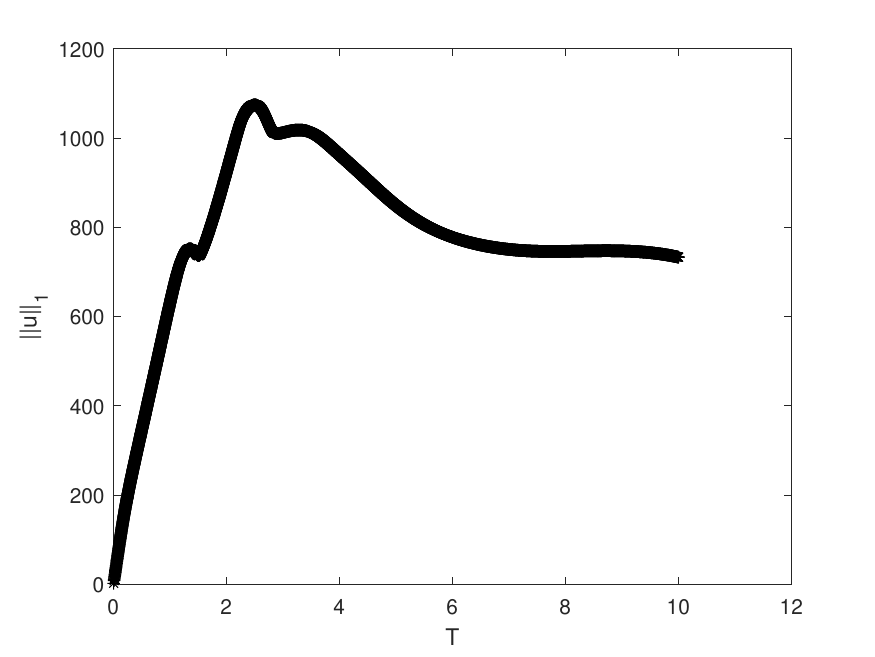}
  \includegraphics[width=0.3\textwidth]{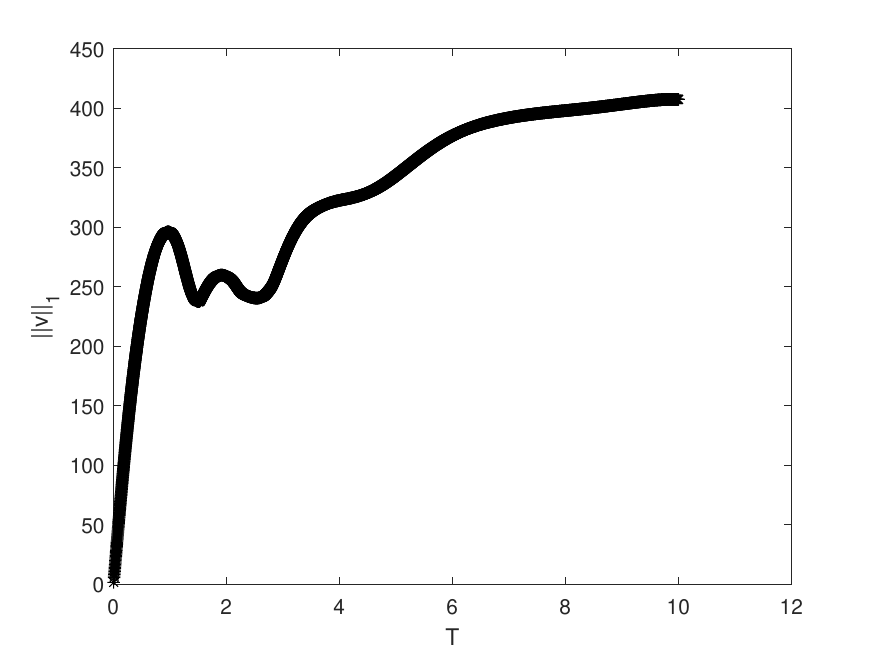}
  \caption{Test case 2.4 - fully 3D surface with obstacle:
    wetting-and-drying simulation. Time behavior of $\Lspace[\infty]-$
    and $\Lspace[1]-$norms of the conserved variables.}
  \label{TC4-wetdry}
\end{figure}

\autoref{Fully3Dsurface2} shows the evolution of the fluid depth at
different times in the interval $T\in[0,5]$~s.
\autoref{Fully3DsurfacWB1}, simlar to \autoref{fig:Slopeplane2} reports the
results at times $T=0.5, 2$~s and $T=8$~s. As in previous test cases,
this figure underlines that the well-balance property of
the system is preserved by the numerical scheme also in the presence
of non-trivial domain curvatures and that the expected symmetries of
the fluid flow around the obstacle are respected.
\autoref{TC4-wetdry} shows the time behavior of both $\Lspace[\infty]$
(top panels) and $\Lspace[1]$ (bottom panels) norms of $\eta$, $u$ and
$v$ in the long-time up to $T=8$~s.

\subsection{Test case 3: multiple obstacles}

This final test case evaluates the qualitative performance of the
solver in a complex environment characterized by multiple localized
discontinuities. The topography consists of an inclined plane, where
$\BSM(x,y) = -0.1y + 1$ within the domain
$\SubsetU = [-4,4] \times [-10,10]$ discretized with a grid of
$256\times640$ cells over which four parabolic localized obstacles are
superimposed given by
$
\obstacle_i(x,y) = 1.2 - (x-a_i)^2 - (y-b_i)^2,
$
with centers $\textbf{c}_i=(a_i,b_i)$, with $\textbf{c}_1=(-2,4)$,
$\textbf{c}_2=(2,-4)$, $\textbf{c}_3=(2,2)$, and
$\textbf{c}_4=(-2,-2)$. The initial condition defines a dam-break
problem with $\Depth=2$~m for $y\leq 9$ and $\Depth=1$~m otherwise. A
CFL~$=0.2$ is used to ensure stability across the sharp gradients
introduced by the presence of multiple obstacles.

\begin{figure}
  \centering
  \includegraphics[width=\textwidth]{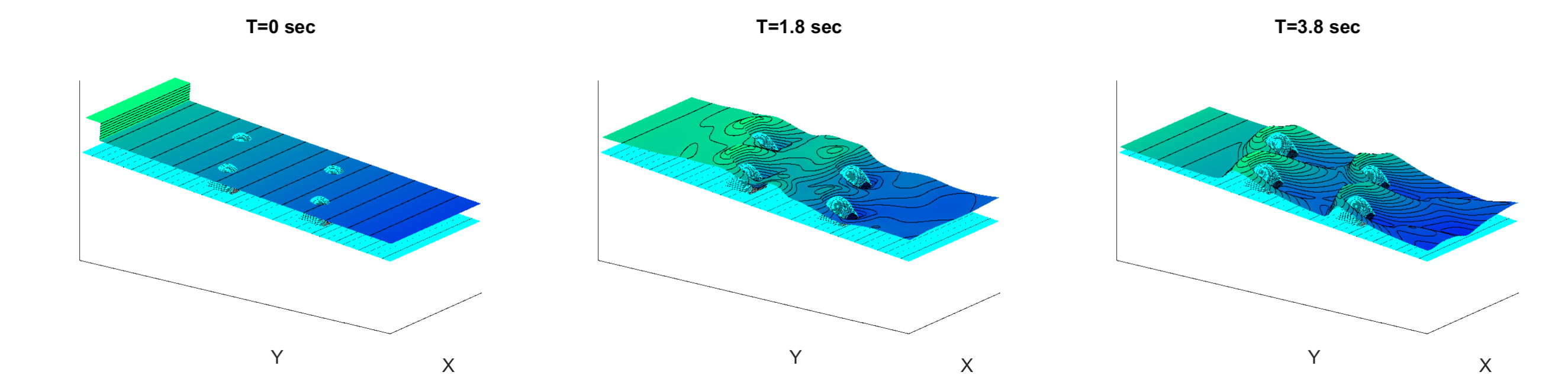}\\
  \includegraphics[width=\textwidth]{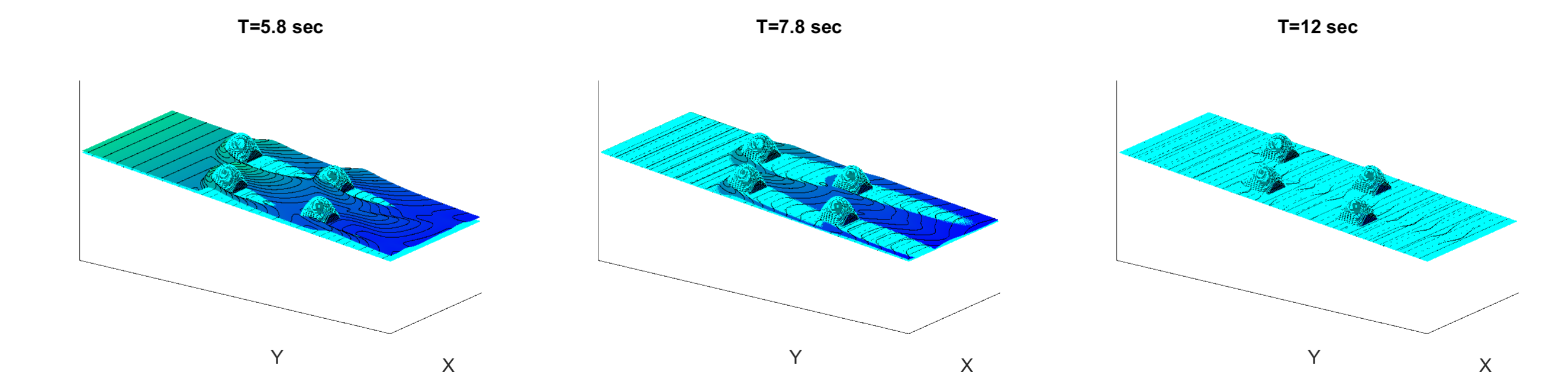}\\
  \caption{Test case 3 - evolution of a drainage process over multiple
    obstacles.  Top-left panel: initial configuration at
    $T=0$~s. Between $T=1.8$~s (top-central panel) and $T=7.8$~s
    (bottom-central panel), the simulation captures the emergence of
    distinct wet regions with non-zero velocity and dry regions where
    the velocity vanishes. It is clear that the drainage patterns are
    strongly influenced by the presence of the topographic
    bumps. Bottom-right panel: final drained steady-state at $T=12$~s.}
  \label{TC33}
\end{figure}

\autoref{TC33} shows the fluid evolution and the convergence to a
non-trivial fully drained steady-state. It is possible to observe
simultaneous dynamic wet-regions, i.e. areas with zero depth and
non-zero velocity corresponding to non-zero eigenvalues in the
hyperbolic system, and completely dry steady-regions, i.e. regions
where the depth vanishes or remains stagnant downstream of the
obstacles and characterized by zero velocity and vanishing
eigenvalues.
We want to underline that our intrinsic Lagrangian-Eulerian scheme is
capable to handle these complex situations without additional
conditions or artificial thresholding.

\section{Concluding remarks and future directions}

In this work, we have extended the Intrinsic Lagrangian-Eulerian
Finite Volume method for Intrinsic Shallow Water Equations on curved
bottoms to take into account localized topographic discontinuities. By
reformulating the discontinuities as local perturbations within our
intrinsic geometric framework, we successfully integrated complex
surface features without necessitating global modifications to the
system structure.

The core strength of our formulation resides in avoiding traditional
Riemann-solver-based techniques. Without requiring local Riemann
problems, the scheme overcomes limitations in handling non-strictly
hyperbolic systems and relaxes some mathematical assumptions, such
as diagonalizability of the Jacobian matrix. To achieve this result,
our approach utilizes the concept of no-flow curves and ensures
stability via a novel ``weak'' CFL condition. The framework allows the
incorporation of discontinuous forcing functions directly into the
Eulerian step, without the need for auxiliary equations to manage the
loss of hyperbolicity typically triggered by topographic jumps.

The numerical results confirm that the Intrinsic Lagrangian-Eulerian
framework effectively circumvents the historical difficulties
associated with gravity-driven flows:

\begin{itemize}

\item Well-Balance: ``lake-at-rest'' benchmarks demonstrate natural
  exact discrete balance between flux gradients and geometric source
  terms, eliminating external root-finding algorithms or complex
  hydrostatic reconstructions.

\item Wetting-and-Drying: velocity stability in long-time simulations
  proves the scheme robustness for $\Depth \rightarrow 0$. By ensuring
  weak positivity \cite{zbMATH07987763,art:abreu22} (in the sense of
  Lax and Liu), our method maintains physically consistent moving
  fronts without non-physical oscillations.

\item Hyperbolicity and Curvature: by respecting the no-flow curve
  stability limits, the formulation ensures the system remains
  well-posed, even in the presence of surface curvatures. The wave
  propagation speeds are captured with high fidelity, regardless of
  the complexity of the underlying metric or topographic obstacles.

\end{itemize}

In summary, the proposed intrinsic formulation provides a versatile
and robust tool for simulating real-life geophysical and engineering
flows - such as debris flows, river floods, and coastal inundation -
where terrains display large curvatures or steep slopes. Our
Riemann-solver-free approach offers greater simplicity and
mathematical flexibility for shallow water equations with complex
topography and friction terms.

The primary objective of this paper is to demonstrate that no-flow
curves delimit a well-balanced region in space-time. Despite employing
an explicit forward Euler time-stepping scheme with linear spatial
reconstruction, the proposed approach accurately captures the
fundamental flow characteristics. Note that a number of questions are
still open about the detailed properties of both the ISWE model and
our Lagrangian-Eulerian numerical approach. They include studying the
structure of nonlinear wave interaction, the dissipative character of
the ISWE approach, the energy and entropy conservation properties and
the amount of artificial dissipation introduced in the
scheme. Additionally, future research will extend this framework to
high-resolution spatio-temporal discretizations and provide a rigorous
analysis of the physical dissipation when friction effects induced by
the Manning formulation are present.

Besides the above future objectives, we also plan to extend this
intrinsic framework to topographies with global discontinuities and
investigate multi-layer systems. In such contexts, the ability to
handle non-diagonalizable interactions between layers will be
critical. Our intrinsic geometric approach will remain fundamental to
accurately taking into account the evolution of layered interfaces
while further pushing the stability limits of the Lagrangian-Eulerian
solver.

\section*{Statements and declarations}

\subsection*{Acknowledgements and fundings} 
E. Abreu acknowledges the financial support of the National Council
for Scientific and Technological Development (CNPq) and São Paulo
Research Foundation (FAPESP), both located in Brazil, under grants
number 307641/2023-6 and 2025/07662-6, respectively.  E. Bachini and
M. Putti received funding from the European Union - Next Generation EU
National Recovery and Resilience Plan [NRRP], Mission 4, Component 2,
Investment 1.3–D.D. 1243 2/082022, PE0000005 Extended Partnership
“RETURN: Multi- Risk sciEnce for resilienT commUnities underR a
changiNg climate” and Mission 4, Component 2, Investment 1.5 - Call
for tender No. 3277 of 30 dicembre 2021 ECS00000043, No. 1058 of June
23, 2023, CUP C43C22000340006, “iNEST: Interconnected Nord-Est
Innovation Ecosystem”. E. Bachini and M. Putti are members of the
Gruppo Nazionale Calcolo Scientifico - Istituto Nazionale di Alta
Matematica (GNCS-INdAM).

\subsection*{Author Contributions} 

E. Abreu: Conceptualization, Methodology, Formal analysis, Writing—review and editing, Implementation, Software, Validation, Data curation.
E. Bachini: Conceptualization, Methodology, Formal analysis, Writing—review and editing. 
J. Perez: Conceptualization, Methodology, Formal analysis, Writing—review and editing, Implementation, Software, Validation, Data curation.
M. Putti: Conceptualization, Methodology, Formal analysis, Writing—review and editing. 

\subsection*{Conflict of interest} The authors have
no financial interests to disclose.

\subsection*{Code availability} Available upon request to the authors.

\appendix
\section{Calculation details: numerical integration of the source term}
\label{app:quadrule}

In this section we report the explicit calculation of the quadrature
approximation for the term $I(\bSource)$ on a Cartesian grid.

\paragraph{Midpoint rule.}

Consider the region $D^n$, and recall that
$\subVolume{}{n}:=\IT[n]\times\Region\subset\REAL\times\SurfDomain$,
with boundary
$\partial \subVolume{}{} = \Region[n] \cup \mathcal{LS}^n\cup
\bRegion[n+1]$ where $\Region[n]$ is the region (cell) on the bottom
surface at time $\tempo^n$ and $\bRegion[n+1]$ is the same region
evolved at times $\tempo^{n+1}$, and $\mathcal{LS}^n$ is the lateral
no-flow surface in space-time.  The latter is then discretized into
the four local no-flow curves associated with the four cell edges,
which in the special case of a Cartesian grid are given by:
\begin{align*}
  \sigma_{i-\frac{1}{2},j}^{\beta}(\tempo)
  &\approx {x_{i-\frac{1}{2}}} + (\tempo -
  \tempo^n){\frac{\FluxLagrangian_{W}^{(\beta,2)}\metrcoefH{1}^2}{\ConservVar{}{\beta,n}}}\,,
  \quad
  \sigma_{i+\frac{1}{2},j}^{\beta}(\tempo) \approx x_{i+\frac{1}{2}} + (\tempo -
  \tempo^n){\frac{\FluxLagrangian_{E}^{(\beta,2)}\metrcoefH{1}^2}{\ConservVar{}{\beta,n}}}\,,\\ \sigma_{i,j-\frac{1}{2}}^{\beta}(\tempo)
  &\approx y_{j-\frac{1}{2}} + (\tempo -
  \tempo^n){\frac{\FluxLagrangian_{S}^{(\beta,3)}\metrcoefH{2}^2}{\ConservVar{}{\beta,n}}}\,,
  \quad
  \sigma_{i,j+\frac{1}{2}}^{\beta}(\tempo)
  \approx y_{j+\frac{1}{2}} + (\tempo -
  \tempo^n){\frac{\FluxLagrangian_{N}^{(\beta,3)}\metrcoefH{2}^2}{\ConservVar{}{\beta,n}}}\,.
\end{align*}
where $i,j$ are the cell vertices and $i\pm1/2, j\pm1/2$ are the edge
midpoints, and $\beta=1,2,3$. 

In this setting, the quadrature approximation of
$I({\bSource})$ can be computed by:
\begin{equation}
  \label{eq:midpoint}
  \begin{aligned}
    I(\bSource)&=\int_{D^n} {\bSource}(\sv,\ConservVar{}{}) 
    = \int_{\tempo^n}^{\tempo^{n+1}}
    \int_{\sigma_{i,j-\frac{1}{2}}^{\beta}(\tempo)}^{\sigma_{i,j+\frac{1}{2}}^{\beta}(\tempo)}
    \int_{\sigma_{i-\frac{1}{2},j}^{\beta}(\tempo)}^{\sigma_{i+\frac{1}{2},j}^{\beta}(\tempo)} {\bSource}(\xcl,\ycl,\ConservVar{}{}(\tempo))\,\Diff{\xcl}\,\Diff{\ycl}\,\Diff{\tempo} \\ 
    & \approx \int_{\tempo^n}^{\tempo^{n+1}}
    \left[\sigma_{i+\frac{1}{2},j}^{\beta}(\tempo)-\sigma_{i-\frac{1}{2},j}^{\beta}(\tempo)\right]
    \left[\sigma_{i,j+\frac{1}{2}}^{\beta}(\tempo)-\sigma_{i,j-\frac{1}{2}}^{\beta}(\tempo)\right]\\
    &\phantom{ \approx \int_{\tempo^n}^{\tempo^{n+1}}\sigma_{i,j+\frac{1}{2}}} {\bSource}\left(\frac{\sigma_{i+\frac{1}{2},j}^{\beta}(\tempo)+\sigma_{i-\frac{1}{2},j}^{\beta}(\tempo)}{2}, \frac{\sigma_{i,j+\frac{1}{2}}^{\beta}(\tempo)+\sigma_{i,j-\frac{1}{2}}^{\beta}(\tempo)}{2},\ConservVar{}{}(\tempo)\right)\,\Diff{\tempo}
    \\
    &\approx \int_{\tempo^n}^{\tempo^{n+1}}
    \meshparam \meshparam\; \bSource\left(\frac{\sigma_{i+\frac{1}{2},j}^{\beta}(\tempo)+\sigma_{i-\frac{1}{2},j}^{\beta}(\tempo)}{2}, \frac{\sigma_{i,j+\frac{1}{2}}^{\beta}(\tempo)+\sigma_{i,j-\frac{1}{2}}^{\beta}(\tempo)}{2},\ConservVar{}{}(\tempo)\right)\,\Diff{\tempo}
    \\
    &\approx \meshparam^2 \Dt\;
    \bSource\left(\frac{\sigma_{i+\frac{1}{2},j}^{\beta}
      \left(\tempo^{n}+\frac{\Dt}{2}\right)+\sigma_{i-\frac{1}{2},j}^{\beta}
      \left(\tempo^{n}+\frac{\Dt}{2}\right)}{2},\right.\\
    &\phantom{\approx \meshparam^2 \Dt\;
      \bSource(\frac{\sigma_{i+\frac{1}{2},j}^{\beta}}{2})}\left.
    \frac{\sigma_{i,j+\frac{1}{2}}^{\beta}
      \left(\tempo^{n}+\frac{\Dt}{2}\right)+\sigma_{i,j-\frac{1}{2}}^{\beta}
      \left(\tempo^{n}+\frac{\Dt}{2}\right)}{2},
    \ConservVar{}{}(\tempo^{n}+{\Dt/2})\right) \\
    &\approx \meshparam^2\,\Dt\;
    \bSource\left(A,B,\ConservVar{}{}(\tempo^n+\Dt/2)\right)\;.
  \end{aligned}
\end{equation}
where we defined $A$ and $B$ by:
\begin{align*}
  A &= \frac{\Dt}{2}\frac{\FluxLagrangian_{W}^{(\beta,2)}\metrcoefH{1}^2}{\ConservVar{}{\beta,n}}
  + x_{i+\frac{1}{2},j}^n
  +\frac{\Dt}{2}\frac{\FluxLagrangian_{E}^{(\beta,2)}\metrcoefH{1}^2}{\ConservVar{}{\beta,n}} + x_{i-\frac{1}{2},j}^n \\
  B &= \frac{\Dt}{2}\frac{\FluxLagrangian_{N}^{(\beta,3)}\metrcoefH{2}^2}{\ConservVar{}{\beta,n}} 
  + x_{i,j+\frac{1}{2}}^n
  +\frac{\Dt}{2}\frac{\FluxLagrangian_{S}^{(\beta,3)}\metrcoefH{2}^2}{\ConservVar{}{\beta,n}} + x_{{i,j-\frac{1}{2}}}^{2}
\end{align*}
Finally, \autoref{eq:midpoint} can be re-written as follows:
\begin{multline*}
  I( {\bSource} )
  \approx \Dt\,\meshparam^2 \,
  \bSource\left(x_j^n+\left[\frac{\FluxLagrangian_{W}^{(\beta,2)}\metrcoefH{1}^2}{\ConservVar{}{\beta,n}}
    -
    \frac{\FluxLagrangian_{E}^{(\beta,2)}\metrcoefH{1}^2}{\ConservVar{}{\beta,n}}\right]\frac{\Dt}{2} + \frac{\meshparam}{2},\right.\\
  \left. y_j^n+\left[\frac{\FluxLagrangian_{N}^{(\beta,2)}\metrcoefH{1}^2}{\ConservVar{}{\beta,n}} - \frac{\FluxLagrangian_{S}^{(\beta,2)}\metrcoefH{1}^2}{\ConservVar{}{\beta,n}}\right]\frac{\Dt}{2} + \frac{\meshparam}{2},
  \ConservVar{}{}(\tempo^n + \Dt/2)\right)\,.
\end{multline*}

\paragraph{The trapezoidal rule.}
Extending the previous computations also in the case of the
trapezoidal rule, the integral $I(\bSource)$ over $D^n$ can be approximated by:

\begin{equation*}
  \label{eq:trap1}
  \begin{aligned}
    I(\bSource) &= \int_{\tempo^n}^{\tempo^{n+1}} \int_{\sigma_{i,j-\frac{1}{2}}^{\beta}(\tempo)}^{\sigma_{i,j+\frac{1}{2}}^{\beta}(\tempo)}
    \int_{\sigma_{i-\frac{1}{2},j}^{\beta}(\tempo)}^{\sigma_{i+\frac{1}{2},j}^{\beta}(\tempo)}
    \bSource(\xcl,\ycl,\ConservVar{}{}(\tempo))\,\Diff{\xcl}\,\Diff{\ycl}\Diff{\tempo} \\[0.7em]
    &\approx \frac{1}{2} \int_{\tempo^n}^{\tempo^{n+1}}
    \int_{\sigma_{i,j-\frac{1}{2}}^{\beta}(\tempo)}^{\sigma_{i,j+\frac{1}{2}}^{\beta}(\tempo)}
    \left[ \sigma_{i+\frac{1}{2},j}^{\beta}(\tempo) -
      \sigma_{i-\frac{1}{2},j}^{\beta}(\tempo)\right]\\
    &\phantom{\approx \frac{1}{2} \int\int}
    \left(
    \bSource(\sigma_{i+\frac{1}{2},j}^{\beta}(\tempo),\ycl,\ConservVar{}{}(\tempo))
    + \bSource(\sigma_{i-\frac{1}{2},j}^{\beta}(\tempo),\ycl,\ConservVar{}{}(\tempo)) \right)\,\Diff{\ycl}\,\Diff{\tempo} \\[0.7em] 
    &\approx \frac{1}{4} \int_{\tempo^n}^{\tempo^{n+1}}
    \left[\sigma_{i+\frac{1}{2},j}^{\beta}(\tempo) -
      \sigma_{i-\frac{1}{2},j}^{\beta}(\tempo)\right]
    \left[ \sigma_{i,j+\frac{1}{2}}^{\beta}(\tempo) -
      \sigma_{i,j-\frac{1}{2}}^{\beta}(\tempo)\right]\\
    &\phantom{\approx \frac{1}{4} \int
      [\sigma]}
    \left[
      \bSource(\sigma_{i+\frac{1}{2},j}^{\beta}(\tempo),\sigma_{i,j+\frac{1}{2}}^{\beta}(\tempo),\ConservVar{}{}(\tempo))
      +
      \bSource(\sigma_{i-\frac{1}{2},j}^{\beta}(\tempo),\sigma_{i,j+\frac{1}{2}}^{\beta}(\tempo),\ConservVar{}{}(\tempo))\right.\\
      & \phantom{\approx \frac{1}{4} \int
        [\sigma][\sigma]\bSource}
      \left.
      +  \bSource(\sigma_{i+\frac{1}{2},j}^{\beta}(\tempo),\sigma_{i,j-\frac{1}{2}}^{\beta}(\tempo),\ConservVar{}{}(\tempo)) + \bSource(\sigma_{i-\frac{1}{2},j}^{\beta}(\tempo),\sigma_{i,j-\frac{1}{2}}^{\beta}(\tempo),\ConservVar{}{}(\tempo)) 
      \right]\,\,\Diff{\tempo} \\[0.7em] 
    &\approx  \frac{\Dt}{4}
    \left[\sigma_{i+\frac{1}{2},j}^{\beta}(\tempo^n) -
      \sigma_{i-\frac{1}{2},j}^{\beta}(\tempo^n)\right]
    \left[ \sigma_{i,j+\frac{1}{2}}^{\beta}(\tempo^n) -
      \sigma_{i,j-\frac{1}{2}}^{\beta}(\tempo^n)\right]\\
    &\phantom{\approx  \frac{\Dt}{4} [\sigma]}
    \left[
      \bSource(\sigma_{i+\frac{1}{2},j}^{\beta}(\tempo^n),\sigma_{i,j+\frac{1}{2}}^{\beta}(\tempo^n),\ConservVar{}{n})
      +
      \bSource(\sigma_{i-\frac{1}{2},j}^{\beta}(\tempo^n),\sigma_{i,j+\frac{1}{2}}^{\beta}(\tempo^n),\ConservVar{}{n})\right.\\
      &\phantom{\approx  \frac{\Dt}{4} [\sigma][\sigma]}
      \left.
      +  \bSource(\sigma_{i+\frac{1}{2},j}^{\beta}(\tempo^n),\sigma_{i,j-\frac{1}{2}}^{\beta}(\tempo^n),\ConservVar{}{n}) + \bSource(\sigma_{i-\frac{1}{2},j}^{\beta}(\tempo^n),\sigma_{i,j-\frac{1}{2}}^{\beta}(\tempo^n),\ConservVar{}{n}) 
      \right]\\[0.7em]
    &= \frac{\Dt\,\meshparam^2}{4}\left[
      \bSource(x_{i+\frac{1}{2}}^{\beta},y_{j+\frac{1}{2}}^{\beta},\ConservVar{}{n}) +
      \bSource(x_{i-\frac{1}{2}}^{\beta},y_{j+\frac{1}{2}}^{\beta},\ConservVar{}{n})
      +\right.\\
      &\phantom{= \frac{\Dt\,\meshparam^2}{4}[\bSource]}
      \left.
      \bSource(x_{i+\frac{1}{2}}^{\beta},y_{j-\frac{1}{2}}^{\beta},\ConservVar{}{n})
      + \bSource(x_{i-\frac{1}{2}}^{\beta},y_{j-\frac{1}{2}}^{\beta},\ConservVar{}{n})
      \right]
  \end{aligned}
\end{equation*}

\section{Refinement study}
\label{app:ref}
\paragraph{Continuous bottom topography.}

\begin{table}
  \centering
  \begin{tabular}{ c  c  c  c  c }
    $h_x = h_y$ &  PLANE   &  PARABOLA   &   HCB   &   3D-SURFACE\\
    \hline
    0.06250     & 2.65144e-02  &  3.85305e-02       &   1.94010e-01    &     4.71504e-02   \\
    0.01563     &  2.12006e-02 &  2.59658e-02       &    1.52012e-01   &     3.92747e-02   \\
    0.00391    & 1.93598e-02   & 1.29988e-02         &   1.05041e-01   &     2.51919e-02    \\
    0.00098    &  1.86435e-02  & 8.34601e-03      &    6.98280e-02     &     1.78970e-02 \\
    0.00024     &  1.46213e-02  &  5.16089e-03     &   4.70245e-02     &     1.12399e-02   \\
    0.00006      & 7.21056e-03   &  3.28825e-03    &    3.01720e-02    &     6.42648e-03   \\
    \hline
  \end{tabular}
  \caption{$\Lspace[1]$-norm of the relative differences between
    solutions (depth values)
    at two sequential refinement levels.}
  \label{tab1cont}
\end{table}

\begin{table}
  \centering
  \begin{tabular}{ c  c  c  c  c }
    $h_x = h_y$ &  PLANE   &  PARABOLA   &   HCB  &   3D-SURFACE\\
    \hline
    0.06250     & 5.48628e-00 &  2.22455e-02     &  1.15492e-01      &   1.58723e-01  \\
    0.01563     & 9.69722e-01 &  1.18538e-02     &   9.70787e-02     &   1.01976e-01  \\
    0.00391     & 1.52901e-01  &  7.08104e-03    &    5.95621e-02    &   6.14649e-02  \\
    0.00098     & 1.86435e-02  &  4.17257e-03    &     3.54222e-02   &   3.56352e-02  \\
    0.00024    & 1.36435e-02  &  2.10924e-03     &     1.89544e-02   &   1.76679e-02  \\
    0.00006    & 1.01069e-02 &  8.22064e-04      &    7.54299e-03    &   3.35342e-03  \\
    \hline
  \end{tabular}
  \caption{$\Lspace[1]$-norm of the relative differences between
    the solution (depth values) at a given refinement level and the finest refinement.}
  \label{tab2cont}
\end{table}

\autoref{tab1cont} and \autoref{tab2cont} show the numerical convergence of
the intrinsic Lagrangian-Eulerian FV scheme for the ISWE model on the
test cases in \cite{art:ABREU2023iswe}, i.e. regular (continuous)
bottom surfaces. For each level of refinement
of the grid, we calculate the relative $\Lspace[1]$-errors of
the numerical depth. In \autoref{tab1cont}, we compute the errors on the
difference of the solution between two sequential refinement levels;
while, in \autoref{tab2cont}, the errors are computed on the difference
between a given level of refinement and the finest refinement. We can
note that, the observed order of convergence is consistent with
an a priori error estimates of $\BigO{\meshparam^{\alpha}}$
with  $\alpha\in[0.25,0.5]$.

\paragraph{Locally discontinuous bottom topography.}

\begin{table}
  \centering
  \begin{tabular}{ c  c  c  c  c }
    $h_x = h_y$ &  PLANE   &  PARABOLA   &  HCB  & 3D-SURFACE\\
    \hline
   0.06250     & 1.55728e-02  & 2.37793e-02      &  1.34914e-01   &   5.67494e-02 \\
    0.01563    &  1.00175e-02 &  2.17714e-02    &    7.87147e-02  &   4.27791e-02 \\
    0.00391    & 7.63599e-03   & 1.29988e-02     &   4.96589e-02   &   3.08401e-02 \\
    0.00098    &  4.83082e-02  & 1.27666e-02     &   1.86103e-02   &   2.16371e-02 \\
    0.00024    & 3.04227e-03  &  1.01799e-02     &   4.70245e-02  &   1.41185e-02 \\
    0.00006    &  1.41513e-03  &  1.14275e-03   &    1.67135e-02  &   6.72282e-03 \\
    \hline
  \end{tabular}
  \caption{$\Lspace[1]$-norm of the relative differences between
    solutions (depth values) at two sequential refinement levels.}
  \label{tab1disc}
\end{table}

\begin{table}
  \centering
  \begin{tabular}{ c  c  c  c  c }
    $h_x = h_y$ &  PLANE   &  PARABOLA   &   HCB  &   3D-SURFACE\\
    \hline
    0.06250     & 3.81885e-02 &  7.91877e-02       &  2.85631e-01   &   1.87205e-01  \\
    0.01563     & 2.52503e-02 &  6.14942e-02       &  1.67627e-01   &   1.20504e-01  \\
    0.00391     & 1.72809e-02  &  4.43324e-02      & 1.00541e-01    &   7.50085e-02 \\
    0.00098     & 1.06938e-02  &  2.97128e-02      & 5.64059e-02    &   4.37445e-02 \\
    0.00024    & 6.59381e-02  &  1.76131e-02      &  2.76677e-02    &   2.24702e-02 \\
    0.00006    & 4.02438e-02   &  7.59395e-03     & 1.01721e-02     &   8.52291e-03 \\
    \hline
  \end{tabular}
  \caption{$\Lspace[1]$-norm of the relative differences between
    the solution (depth value) at a given refinement level and the finest refinement.}
  \label{tab2disc}
\end{table}

\autoref{tab1disc} and \autoref{tab2disc} show the numerical convergence of
the intrinsic Lagrangian-Eulerian FV scheme for the ISWE model in the
case of locally discontinuous bottom topographies, as defined
in \autoref{sec:results}. For each level of mesh refinement, we calculate
the relative $\Lspace[1]$-errors of the depth. Similarly to the
continuous case, in \autoref{tab1disc} we show the errors of the
difference of the solution between two sequential level of refinement;
while, in \autoref{tab2disc}, the errors on the differences of the
solution are computed between a given level and the finest refinement
level. Again, it is possible to note that the observed order of
convergence is consistent with the a priori error estimates of
$\BigO{\meshparam^{\alpha}}$ with $\alpha\in[0.25,0.5]$.
%

\bibliographystyle{abbrvnat}
\bibliography{strings,biblio}

\end{document}